\documentclass[12pt]{article}

\newif\ifpreprint
\preprinttrue

\pdfoutput=1
\usepackage{fullpage}

\newenvironment{keywords}{\paragraph{Keywords:}}{}
\def\sep{, }

\usepackage{siunitx}

\usepackage[latin1]{inputenc}
\usepackage{amsfonts}
\usepackage{amssymb}
\usepackage{amsmath}
\usepackage{color}
\usepackage{euscript}
\usepackage{graphicx}
\usepackage{stmaryrd} % used by \ic
\usepackage[colors]{optsys}
\usepackage{comment}

\newcommand{\mathscr}{\EuScript}

\newcommand{\optstate}{\state}
\newcommand{\OptState}{\State}

\newcommand{\icinterval}[2]{\ic{#1,#2}}
\newcommand{\kernelargs}[2]{\nsetp{#2}{#1}}
\newcommand{\nkernelargs}[2]{\nsetp{#2}{#1}}
\newcommand{\bkernelargs}[2]{\bsetp{#2}{#1}}
\newcommand{\Bkernelargs}[2]{\Bsetp{#2}{#1}}
\renewcommand{\bRR}{[0,+\infty]}

\newcommand{\HistoryFeedback}{\gamma}
\newcommand{\HISTORYFEEDBACK}{\Gamma}

\newcommand{\Flow}[3]{\Phi_{#1:#2}^{#3}}

\newcommand{\StochKernel}[2]{\rho_{#1:#2}}
\newcommand{\tildeStochKernel}[2]{{\widetilde\rho}_{#1:#2}}
\newcommand{\StochKernelFeed}[3]{\rho_{#1:#2}^{#3}}
\newcommand{\varStochKernelFeed}[3]{\varrho_{#1:#2}^{#3}}
\newcommand{\tildevarStochKernelFeed}[3]{\tilde{\varrho}_{#1:#2}^{#3}}

\newcommand{\BelOp}[2]{\mathcal{B}_{#1:#2}}
\newcommand{\tildeBelOp}[2]{{\widetilde{\mathcal{B}}}_{#1:#2}}
\newcommand{\overlineBelOp}[2]{{\overline{\mathcal{B}}}_{#1:#2}}
\newcommand{\naturalBelOp}[2]{{\mathcal{B}^\Before_{#1:#2}}}

\newcommand\InstantneousCost{\ell}
\newcommand\FinalCost{\ell}

\newcommand{\volume}{S}
\newcommand{\turbined}{Q}
\newcommand{\spilled}{R}
\newcommand{\inflow}{A}

\newcommand{\tildecriterion}{\tilde \jmath}
\newcommand{\barcriterion}{\overline \jmath}

\newcommand{\Observation}{Y}

\usepackage{xparse}

\NewDocumentCommand{\dsum}{ O{} O{} }{\ds\sum_{#1}^{#2}}
\NewDocumentCommand{\dprod}{ O{} O{} }{\ds\prod_{#1}^{#2}}
\NewDocumentCommand{\dargmin}{ O{} O{} }{\ds\argmin_{#1}^{#2}}
\NewDocumentCommand{\dargmax}{ O{} O{} }{\ds\argmax_{#1}^{#2}}

\newcommand{\SlowSpan}{\mathsf{S}}

\newcommand{\topp}[1]{\bar{#1}}
\newcommand{\bottom}[1]{\underline{#1}}
\newcommand{\tslow}{s}
\newcommand{\FastSpan}{\mathsf{F}}
\newcommand{\tfast}{f}

\newcommand{\nex}[1]{{#1}^+}

\newcommand{\pre}[1]{{#1}^-}

\newcommand{\slow}{\mathsf{s}}
\newcommand{\fast}{\mathsf{sf}}

\newcommand{\bm}[1]{#1}

\newcommand{\tDHD}{t}%{s}
\newcommand{\horizonDHD}{\horizon}
\newcommand{\tDHDIt}{t'}

\newcommand{\Before}{\natural}
\newcommand{\After}{\flat}

\NewDocumentCommand{\conso}{O{} O{\slow}}{u_{#1}^{#2}}
\NewDocumentCommand{\vaconso}{O{} O{\slow}}{\bm{u}_{#1}^{#2}}
\NewDocumentCommand{\CONSO}{O{} O{\slow}}{\mathbb{U}_{#1}^{#2}}

\NewDocumentCommand{\price}{O{}}{w_{#1}^{\slow}}
\NewDocumentCommand{\vaprice}{O{}}{\bm{w}_{#1}^{\slow}}
\NewDocumentCommand{\PRICE}{O{}}{\mathbb{W}_{#1}^{\slow}}

\NewDocumentCommand{\buy}{O{} O{\fast}}{u_{#1}^{#2}}
\NewDocumentCommand{\vabuy}{O{} O{\fast}}{\bm{u}_{#1}^{#2}}
\NewDocumentCommand{\BUY}{O{} O{\fast}}{\mathbb{U}_{#1}^{#2}}

\NewDocumentCommand{\cost}{O{} O{\fast}}{w_{#1}^{#2}}
\NewDocumentCommand{\vacost}{O{} O{\fast}}{\bm{w}_{#1}^{#2}}
\NewDocumentCommand{\COST}{O{} O{\fast}}{\mathbb{W}_{#1}^{#2}}

\NewDocumentCommand{\stock}{O{} O{\slow}}{\optstate_{#1}^{#2}}
\NewDocumentCommand{\vastock}{O{} O{\slow}}{\bm{\optstate}_{#1}^{#2}}
\NewDocumentCommand{\STOCK}{O{} O{\slow}}{\mathbb{\OptState}_{#1}^{#2}}

\NewDocumentCommand{\buffer}{O{} O{\fast}}{\optstate_{#1}^{#2}}
\NewDocumentCommand{\vabuffer}{O{} O{\fast}}{\bm{\optstate}_{#1}^{#2}}
\NewDocumentCommand{\BUFFER}{O{} O{\fast}}{\mathbb{\OptState}_{#1}^{#2}}

\NewDocumentCommand{\instantcost}{O{} O{}}{\Lambda_{#1}^{#2}}
\NewDocumentCommand{\tilValue}{O{} O{} O{}}{#3{V}_{#1}^{#2}}

\newcommand{\dyn}{G} %{\mathcal{F}}

\newcommand{\hatdyn}{\mathcal{\dyn}}

\NewDocumentCommand{\alea}{O{} O{}}{w_{#1}^{#2}}
\NewDocumentCommand{\ALEA}{O{} O{}}{\mathbb{W}_{#1}^{#2}}

\NewDocumentCommand{\hist}{O{} O{}}{h_{#1}^{#2}}
\NewDocumentCommand{\Hist}{O{} O{}}{\mathcal{H}_{#1}^{#2}}
\NewDocumentCommand{\HIST}{O{} O{}}{\mathbb{H}_{#1}^{#2}}

\newcommand{\analytically}[1]{\tribu{A}_{#1}}
\newcommand{\universally}[1]{\tribu{U}_{#1}}
\renewcommand{\borel}[1]{\tribu{B}_{#1}}

\usepackage{fullpage}
\usepackage[pdftex]{hyperref}
\usepackage[square,numbers]{natbib}
\usepackage[resetlabels,labeled]{multibib}
\makeatother
\newcommand{\citeAM}[2][]{\citep[AM\!\!][#1]{#2}}
\makeatletter

\graphicspath{{./}}
\ifpreprint
\newcommand{\preprintstart}{}
\newcommand{\preprintstop}{}
\fi

\title{Time Block Decomposition of\\
  Multistage Stochastic Optimization Problems}

\author{%
  Pierre Carpentier\footnote{UMA, ENSTA Paris, Institut Polytechnique de Paris, Palaiseau, France},
  \and
  Jean-Philippe Chancelier\footnote{CERMICS, Ecole des Ponts ParisTech, Marne-la-Vall\'ee, France},
  \and
  Michel De Lara\footnotemark[2],
  \and
  Thomas Martin\footnotemark[2],
  \and
  Tristan Rigaut\footnote{Efficacity, Marne-la-Vall\'ee, France}
}

\date{\today}

\begin{document}
\maketitle

\begin{abstract}
  Multistage stochastic optimization problems are, by essence, complex
  as their solutions are functions of both stages and uncertainties.
  Their large scale nature makes decomposition methods appealing,
  like dynamic programming which is a sequential decomposition
  using a state variable defined at all stages. By contrast,
  in this paper we introduce the notion of state reduction by time
  blocks, that is, at stages that are not necessarily all the original
  stages. Then, we prove a dynamic programming equation with value
  functions that are functions of a state only at some stages.
  This equation crosses over time blocks, but involves a dynamic
  optimization inside each block.
  We illustrate our contribution by showing its potential
  in three applications in multistage stochastic optimization:
  mixing dynamic programming and stochastic programming,
  two-time-scale optimization problems,
  decision-hazard-decision optimization problems.
\end{abstract}


\section{Introduction}
\label{gdp:Introduction}

Solutions of multistage stochastic optimization problems are functions
of both time and uncertainties. This makes such problems complex.
However, their structure makes decomposition methods appealing
to solve them~\cite{Ruszczynski:2003}.
One of the most common approaches are time decomposition
(state-based resolution methods), like stochastic dynamic programming,
in stochastic optimal control, and scenario decomposition,
like progressive hedging, in stochastic programming.
On the one hand, stochastic programming deals with an
underlying random process taking a finite number of values,
called scenarios \cite{Shapiro-Dentcheva-Ruszczynski:2009}.
Solutions are indexed by a scenario tree, the size of which
increases exponentially with the number of stages
(hence generally a few stages in practice).
However, to overcome this obstacle, stochastic programming
takes advantage of scenario decomposition methods
(progressive hedging \cite{Rockafellar-Wets:1991}).
On the other hand, stochastic control deals with
a state model driven by a white noise, that is,
the noise is made of a sequence of independent random variables.
Under such assumptions, stochastic dynamic programming
is able to handle many stages,
as it offers reduction of the search for a solution
among state feedbacks (instead of functions of the past noise)
\cite{Bellman:1957,Puterman:1994}.

In a word, dynamic programming is good at handling
multiple stages --- but at the price of assuming that noises
are stagewise independent --- whereas stochastic programming
does not require such assumption, but can only handle a few stages.
Could we take advantage of both methods?
Is there a way to apply stochastic dynamic programming at a slow time scale
--- a scale at which noises could be considered statistically independent ---
crossing over fast time scale optimization problems where independence
would not hold? This question is one of the motivations of this paper,
and we indeed provide a method to decompose multistage stochastic
optimization problems by time blocks.
This decomposition method and the main result are, mathematically
speaking, quite natural, but the main difficulty is notational.
Indeed, the rigorous formulation of multistage stochastic optimization
problems on so-called history spaces requires rather heavy notation.

Although specialists in stochastic optimal control and dynamic
programming will find the results as natural and non surprising,
or as part of folklore, the fact is that we have not been able
to find references that treat the case of
\emph{a state defined only at a subset of stages}.
This is why we set out to write this paper, without any real
theoretical ambition, but with the objective that this result
be established and can be used for applications using several
forms of decomposition\footnote{The starting point of our reflections
  on this subject were conversations
  that three of us held with Roger Wets in Bogota in 2013. We discussed
  the interest and the way of mixing the techniques of scenario trees
  (to be able to take into account correlated noises) with the techniques
  of dynamic programming (to have a vision of the optimal future costs
  by means of value functions). \label{ft:Bogota2013}}.
This is also why we present three (theoretical) applications in multistage
stochastic optimization: mixing dynamic programming and stochastic programming,
two-time-scale optimization problems, decision-hazard-decision optimization problems.

As there are several ways to tackle the difficulties of dealing with
a large number of time steps, we compare our approach with other ones.
In this paper, we propose an \emph{exact decomposition} of a multistage
stochastic optimization problem \emph{by time blocks}
using \emph{a state defined only at a subset of stages},
to be distinguished from either \emph{time aggregation}
or \emph{approximate decomposition by timescales},
which both yield \emph{approximate} problems. We discuss both now.

\emph{Time aggregation} consists in grouping the time steps, that
is, in considering a partition of the time steps in time blocks
and ``aggregating'' variables and constraints in each time block.
To our knowledge, this approach was initiated in~\cite{Birge:1985}
for stochastic linear programs.
For such linear programs, it is indeed easily conceived that,
by summing (``aggregating'') linear constraints, one obtains lower
bounds for minimization problems. This approach was generalized in
the paper~\cite{Wright:1994} which puts forward a measure-theoretic
framework with coarser and finer filtrations, and uses linear duality.
Then, this was extended in~\cite{Kuhn:2008} for stochastic convex
programs, using filtrations and convex duality. The main idea can
be sketched as follows:
the coarser filtration is used to reduce the measurability of the
decision variables, whereas the finer filtration is used to enlarge
the measurability of the dual variables associated with the constraints,
so that the optimal value of the problem obtained by using these
two filtrations is an upper bound of the true optimal value;
exchanging the role of the filtrations leads to a lower bound.
Thus, with time aggregation, one obtains simpler problems that are
lower and upper bounds for the original minimization stochastic problems,
hence are approximations.

In \emph{approximate decomposition by timescales}, one identifies
several timescales in the original multistage stochastic optimization
problem and then sets up an optimization problem for each timescale.
It is approximate in that the connexion between the  problems formulated
for each timescale and the whole multistage problem is not explicit.

Approximate decomposition by timescales can be done
in the context of dynamic programming, with the value functions obtained
for a given timescale entering the final cost of the problem
at the finer timescale. This approach gives a cascade of easier
to solve optimization problems, and again corresponds
to approximate the original problem.
An example of this approach can be found in~\cite{Cheng-Powell:2016}
where --- for a problem involving both the control of the storage
of a battery (5 minutes time steps) and the frequency regulation
(2 seconds time steps) --- is introduced a first hourly resource
model whose resolution by dynamic programming leads to value
functions used in a five minute storage model as final costs.
The value functions, obtained by solving by dynamic programming
this second model, are themselves used in a 2~seconds frequency model.
Another possibility arises when the considered optimization problem
displays a \emph{periodical behavior}. In that case, a natural
time block decomposition is given by the period of the system.
In~\cite{Shapiro:2020}, by taking into account
such a periodical pattern in the dynamic programming equations,
one significantly reduces the computational effort to solve
the problem using a fixed point approach. Finally,
\cite{Porteiro-Ferragut-Paganini:2018} presents a preliminary
work on extending the Stochastic Dual Dynamic Programming
approach to two-time-scale problems, such as those encountered
in energy systems involving both long-term hydro storages and
short-term battery storages.

Approximate decomposition by timescales can also be  done
in the framework of stochastic programming.
In~\cite{Kaut-Midthun-Werner-Tomasgard-Hellemo-Fodstad:2014},
the authors introduce a slow scenario tree, that is, a tree
involving only the time stages of the slow time scale;
but at each node of this slow scenario tree are attached
fast time scale scenarios, which do not interfere with
the other nodes of the slow scenario tree.
The structure allows one to model and solve problems
that need to combine strategic (long term) and operational
(short term) uncertainty, without the explosion in the problem
size that would follow from using a standard multistage model.
The special situation where decisions are taken only
at the slow time scale (whereas uncertainties occur at each
time stage) is considered in~\cite{Glanzer-Pflug:2020}.
The authors propose to build a scenario tree branching
at the slow time stages, and designed using the theory
of bridge processes between two consecutive nodes in order
to represent the noise at the fast time scale.

The paper is organized as follows. In
Sect.~\ref{gdp:Stochastic_Dynamic_Programming_and_State_Reduction_by_Time_Blocks},
we present stochastic dynamic programming with histories as a way to solve a
stochastic optimal control problem formulated in discrete time.  In
Sect.~\ref{gdp:State_Reduction_by_Time_Blocks}, we revisit the notion of
``state'' by defining state reduction by time blocks --- that is, at time stages
that are not necessarily all the original ones --- and then we prove a reduced
dynamic programming equation.  This is the central contribution.  Then, we
illustrate our contribution by showing its potential in three cases.
In~Sect.~\ref{sec:Mixing_dynamic_programming_and_stochastic_programming}, we
show how to mix dynamic programming and stochastic programming.  In
Sect.~\ref{Two_time_scale_optimization_problem}, we detail how to handle
problems with two time scales, and illustrate this with the crude oil
procurement problem.
In~Sect.~\ref{gdp:Decision_Hazard_Decision_Dynamic_Programming}, we introduce
what we call the decision-hazard-decision framework, and we provide a dynamic
programming equation.  We relegate technical results and proofs in
Appendix~\ref{Technical_Details_and_Proofs}.

\section{Stochastic dynamic programming with histories}
\label{gdp:Stochastic_Dynamic_Programming_and_State_Reduction_by_Time_Blocks}

In~\S\ref{The_Bertsekas-Shreve_setting},
we present the setting to formulate multistage stochastic optimization problems
over the so-called history space, with history feedbacks. Then, to prepare
the main result in Sect.~\ref{gdp:State_Reduction_by_Time_Blocks},
we establish in~\S\ref{Stochastic_dynamic_programming_equation_with_histories}
a dynamic programming equation when the state is the history, that is,
is made of the uncertainties and the controls prior to the current time (see the
``canonical construction'' in~\cite[p.~15]{Hernandez-Lerma-Lasserre:1996}).
Although quite natural, this equation is generally not written in
the literature, as most frameworks in dynamic programming assume
the \emph{a priori} existence of a state.

We use the notation $ \ic{r,s}=\na{r,r+1,\ldots,s-1,s} $ for any two
natural numbers~$r,s$ such that $r \leq s$. We will also use the shorter
notation $r{:}s=\ic{r,s}$, for example in subscripts as in $\history_{r:s}$.
From now on, time is discrete and runs among the natural numbers
$t\in\ic{0,\horizon}$, where $\horizon \in \NN^*$ is a positive
natural number. Finally, we say that a function is \emph{numerical}
if it takes its values in $ \barRR=[-\infty,+\infty] $.

\subsection{The Bertsekas-Shreve setting \cite{Bertsekas-Shreve:1996}}
\label{The_Bertsekas-Shreve_setting}

To obtain a stochastic dynamic programming with histories
requires technical assumptions.
Indeed, as Bertsekas and Shreve notice at the beginning of
\cite[\S7.6]{Bertsekas-Shreve:1996}:
``The dynamic programming algorithm is centered around infimization of
functions, and this is intimately connected with projections of sets'';
``Unfortunately, the  projection of a Borel-measurable set need not be
Borel-measurable. In Borel spaces, however, the  projection
of a Borel-measurable set is an analytic set''.
They devote \cite[\S7.6]{Bertsekas-Shreve:1996} to the definition
and study of analytic sets, and in \cite[\S7.7]{Bertsekas-Shreve:1996}
define universally measurable functions, as well as lower semianalytic
functions.

We call \emph{Borel space}~$\np{\STATE,\borel{\STATE}}$ a Borel set~$\STATE$
equipped with its Borel $\sigma$-field~$\borel{\STATE}$
\cite[Definition~7.7, p.~118]{Bertsekas-Shreve:1996}.
By abuse of notation, we often speak of the {Borel space}~$\STATE$.
There exist two other interesting $\sigma$-fields:
the \emph{analytic $\sigma$-field}~$\analytically{\STATE}$
\cite[Definition~7.19, p.~171]{Bertsekas-Shreve:1996};
the \emph{universal $\sigma$-field}~$\universally{\STATE}$
\cite[Definition~7.18, p.~167]{Bertsekas-Shreve:1996}.
We have the inclusions
\( \borel{\STATE} \subset \analytically{\STATE} \subset \universally{\STATE} \)
\cite[p.~171]{Bertsekas-Shreve:1996}.

For any Borel space~$\STATE$, subset~\( \State\subset\STATE \)
and numerical function \( \varphi : \State\to\barRR \),
the function \( \varphi \) is said to
be \emph{lower semianalytic} \cite[Definition~7.21]{Bertsekas-Shreve:1996}
if the subset~\( \State\) is analytic (\( \State\in \analytically{\STATE}\))
and if the subset \( \defset{\state\in\State}{ \varphi\np{\state}<c } \)
is analytic for all $c\in\RR$. We denote by $ \espace{L}^{0}_{+}(\STATE)$
the space of lower semianalytic nonnegative numerical functions over~$\STATE$.

For a Borel space~$\mathbb{X}$ (resp.~$\mathbb{Y}$)
equipped with the Borel $\sigma$-field~$\borel{\mathbb{X}}$
(resp.~$\borel{\mathbb{Y}}$),
a mapping~$f : \mathbb{X}\to\mathbb{Y}$ is said to be
\emph{universally measurable} \cite[Definition~7.20,
p.~171]{Bertsekas-Shreve:1996} (resp. \emph{Borel-measurable})
if, for all~$B\in\borel{\mathbb{Y}}$,
$f^{-1}(B)\in\universally{\mathbb{X}}$
(resp. $f^{-1}(B)\in\borel{\mathbb{X}}$).

\subsubsubsection{Histories and history spaces}

For each time~$t\in\ic{0,\horizon{-}1}$, the control~$\control_{t}$
takes its values in a Borel space~$\CONTROL_{t}$.
For each time~$t\in\ic{0,\horizon}$, the uncertainty~$\uncertain_{t}$
takes its values in a Borel space~$\UNCERTAIN_{t}$.
For $t\in\ic{0,\horizon}$,
we define the \emph{history space}~$\HISTORY_{t}$
as the product Borel space \cite[Proposition~7.13,
p.~119]{Bertsekas-Shreve:1996}
\begin{equation}
    \HISTORY_{t}
  = \UNCERTAIN_{0} \times \prod_{\tbis=1}^{t}
    \np{ \CONTROL_{\tbis{-}1} \times \UNCERTAIN_{\tbis} }
    \eqfinv
\end{equation}
with the particular case $ \HISTORY_{0} = \UNCERTAIN_{0} $,
$ \tribu{\History}_{0} = \tribu{\Uncertain}_{0} $.
A generic element $ \history_{t}= \bp{ \uncertain_{0},
  \np{\control_{\tbis{-}1},\uncertain_{\tbis}}_{\tbis\in \ic{1,t} }} =$
% TYPO: SHOULD HAVE BEEN   \np{\control_{\tbis{-}1},\uncertain_{\tbis}}_{\tbis \in \ic{1,t}} } =$
$(\uncertain_{0}$, $\control_{0}$, $\uncertain_{1}$, $\control_{1}$,
$\uncertain_{2}$, $\ldots$, $\control_{t-2}$, $\uncertain_{t{-}1}$, $
\control_{t{-}1}$, $\uncertain_{t}) \in \HISTORY_{t} $
is called a \emph{history} at time~$t$.
For $1 \leq \tun \leq \tbis \leq \tter$,
we introduce the $\interval{\tun}{\tbis}$-\emph{history subpart}
$  \history_{\tun:\tbis} =
\np{ \control_{\tun{-}1},\uncertain_{\tun}, \ldots,
  \control_{\tbis{-}1},\uncertain_{\tbis} } \in \HISTORY_{\tun:\tbis}=
\prod_{\tau=\tun}^{\tbis} \np{ \CONTROL_{\tau{-}1} \times \UNCERTAIN_{\tau}} $,
so that we have
$\history_{\tter} = \np{ \history_{\tun{-}1}, \history_{\tun:\tter} }$.

\subsubsubsection{History feedbacks}

For $ 0 \leq \tun \leq \tter \leq \horizon-1 $,
we define a $\interval{\tun}{\tter}$-\emph{history feedback} as a sequence
$ \bseqa{\HistoryFeedback_{\tbis}}{\tbis\in\ic{\tun,\tter}} $
of universally measurable mappings
$  \HistoryFeedback_{\tbis} :
\HISTORY_{\tbis} \to \CONTROL_{\tbis} $.
We call~$\HISTORYFEEDBACK_{\tun:\tter}$
the set of $\interval{\tun}{\tter}$-history feedbacks.
The history feedbacks reflect the following information
structure. At the end of the time interval~$[t-1,t[$,
an uncertainty variable~$\uncertain_{t}$ is revealed.
Then, at the beginning of the time interval~$[t,t+1[$,
a decision-maker chooses a control~$\control_{t}$ contingent
on no more than the past, giving the chronology
\begin{equation}
  \uncertain_{0} \rightsquigarrow \control_{0} \rightsquigarrow
  \uncertain_{1} \rightsquigarrow \control_{1} \rightsquigarrow
  \cdots \rightsquigarrow
  \uncertain_{t} \rightsquigarrow \control_{t} \rightsquigarrow
  \cdots \rightsquigarrow
  \uncertain_{\horizon{-}1} \rightsquigarrow \control_{\horizon{-}1}
  \rightsquigarrow \uncertain_{\horizon}
  \eqfinp
  \label{gdp:eq:interplay_Noise_Control}
\end{equation}

\subsubsubsection{Stochastic kernels}

In what follows, given a Borel space~$\mathbb{Y}$,
$\Delta\np{\mathbb{Y}}$ denotes the Borel space of probability measures
over~$\mathbb{Y}$ (see \cite[Corollary 7.25.1]{Bertsekas-Shreve:1996}).
Uncertainty is represented by a sequence $ \bseqa{ \StochKernel{t{-}1}{t}}%
{t\in\ic{1,\horizon}}$ of Borel-measurable stochastic kernels
(see~\cite[Definition~7.12 and Proposition~7.26, p.~134]{Bertsekas-Shreve:1996})
\begin{equation}
  \StochKernel{t{-}1}{t} :
  \HISTORY_{t{-}1} \to \Delta\np{\UNCERTAIN_{t}}
  \eqsepv
  \forall t \in \ic{1,\horizon}
  \eqfinp
  \label{gdp:eq:sequence_of_stochastic_kernels}
\end{equation}
Thus, for any past history~$\history_{t{-}1}\in\HISTORY_{t{-}1}$, we have that
$\StochKernel{t-1}{t}\np{\history_{t{-}1}} \in \Delta\np{\UNCERTAIN_{t}}$,
the space of probability measures over~$\UNCERTAIN_{t}$.
It is common practice (see~\cite[Definition~7.12, p.~134]{Bertsekas-Shreve:1996})
to use the notation
$\StochKernel{t-1}{t}\np{\dd \uncertain_{t} | \history_{t{-}1}}$
to denote this probability distribution, element of~$\Delta\np{\UNCERTAIN_{t}}$.
So, the notation $| h_{t}$ is here to \emph{evoke} a conditional
distribution (of the next uncertainty knowing the past history),
but it is not introduced as a conditional distribution, but simply
as a way to express a \emph{parametric dependence} (as explicitely
said in \cite[Definition~7.12, p.~134]{Bertsekas-Shreve:1996}).
We could have indifferently written
$\StochKernel{t-1}{t}\np{\dd \uncertain_{t}, \history_{t{-}1}}$
or $\StochKernel{t-1}{t}\np{\dd \uncertain_{t}; \history_{t{-}1}}$.

We define, for any feedback
$\nseqa{\HistoryFeedback_{\tbis}}{\tbis\in\ic{\tter,\horizon\!-\!1}}
\in \HISTORYFEEDBACK_{\tter:\horizon\!-\!1} $,
a new sequence of Borel-measurable stochastic kernels
$ \StochKernelFeed{\tter}{\horizon}{\HistoryFeedback} :
\HISTORY_{\tter} \to \Delta\np{ \HISTORY_{\horizon} } $,
that capture the transitions between histories when the dynamics
$\history_{\tbis{+}1} = \bp{\history_{\tbis},\control_{\tbis},\uncertain_{\tbis{+}1}}$
is driven by $ \control_{\tbis} = \HistoryFeedback_{\tbis}\np{\history_{\tbis}} $
for all $\tbis$ in $\ic{\tter,\horizon{-}1}$
(see Definition~\ref{gdp:de:stochastic_kernels_rho} in Appendix~\ref{Technical_Details_and_Proofs}
for the detailed construction
of~$\StochKernelFeed{\tun}{\tter}{\HistoryFeedback}$).
Note that~$\StochKernelFeed{\tter}{\horizon}{\HistoryFeedback}$
generates a probability distribution on the space~$\HISTORY_{\horizon}$
of histories over the whole timespan~$\ic{0,\horizon}$.

\subsubsubsection{Cost function}

The cost criterion to be minimized is a nonnegative\footnote{We could also
  consider a cost criterion~$ \criterion : \HISTORY_{\tter} \to \barRR $,
  either bounded function, or uniformly
  bounded below function. However, for the sake of simplicity,
  we will deal in the sequel with nonnegative numerical
  functions. The case $ \criterion\np{\history_{\horizon}}=+\infty $
  materializes joint constraints between uncertainties
  and controls in~$\history_{\horizon}$.\label{gdp:ft:nonnegative}}
and lower semianalytic numerical function
\begin{equation}
  \criterion : \HISTORY_{\horizon}   \to \bRR \eqfinp
  \label{gdp:eq:criterion}
\end{equation}
Notice that~\eqref{gdp:eq:criterion} does not represent a cost
at final time, but a cost function of the whole history
$\history_{\horizon}=
(\uncertain_{0}, \control_{0}, \uncertain_{1}, \control_{1},
\ldots, \uncertain_{\horizon-1},\control_{\horizon-1},
\uncertain_{\horizon}) \in \HISTORY_{\horizon}$.
As $ \history_{\horizon}$ contains \emph{all} past controls and
uncertainties, a function $\criterion : \HISTORY_{\horizon} \to \bRR$
covers the most general case. For instance, the function~$\criterion$
can have the special form of a sum of time block costs, like
in Equation~\eqref{gdp:eq:criterion_time_block_additive}.

\subsection{Stochastic dynamic programming equation with histories}
\label{Stochastic_dynamic_programming_equation_with_histories}

\subsubsubsection{Family of optimization problems}

We consider the following family of optimization problems,
indexed by $t$ in $ \ic{0,\horizon{-}1} $ and parameterized by
the history $ \history_{t} \in \HISTORY_{t} $:
for all $t$ in $\ic{0,\horizon{-}1}$, we define the minimum value
\begin{subequations}
  \begin{align}
    \Value_{t}(\history_{t})
    & = \inf_{\HistoryFeedback_{t:\horizon{-}1} \in \HISTORYFEEDBACK_{t:\horizon{-}1}}
      \int_{\HISTORY_{\horizon}} \criterion\np{\history'_{\horizon}}
      \StochKernelFeed{\tter}{\horizon}{\HistoryFeedback}
      \kernelargs{\history_{\tter}}{\dd \history'_{\horizon}}
      \eqsepv \forall \history_{\tter} \in \HISTORY_{\tter}
      \eqfinv
      \intertext{and we also define}
      \Value_{\horizon}(\history_{\horizon})
    & = \criterion(\history_{\horizon})
      \eqsepv \forall \history_{\horizon} \in \HISTORY_{\horizon}
      \eqfinp
  \end{align}
  \label{gdp:eq:value_functions}
  The numerical function $ \Value_{t} : \HISTORY_{t} \to \bRR $
  is called the \emph{value function} at time~$t$.
\end{subequations}

Next, we show how the sequence
$\nseqa{\Value_{t}}{t\in \ic{0,\horizon}}$ of value functions
can be used to solve, via dynamic programming, the optimization
problem of interest, that is, the one starting at~$t=0$, whose
value is (recall that $\history_{0}=\uncertain_{0}$)
\begin{equation}
  \Value_{0}(\uncertain_{0})
  = \inf_{\HistoryFeedback_{0:\horizon{-}1} \in \HISTORYFEEDBACK_{0:\horizon{-}1}}
  \int_{\HISTORY_{\horizon}} \criterion\np{\history'_{\horizon}}
  \StochKernelFeed{0}{\horizon}{\HistoryFeedback}
  \kernelargs{\uncertain_{0}}{\dd \history'_{\horizon}}
  \eqfinp
  \label{eq:Vzero}
\end{equation}

\subsubsubsection{Bellman operators and dynamic programming}

We show that the value functions in~\eqref{gdp:eq:value_functions}
are \emph{Bellman functions}, in that they are solution
of a Bellman or dynamic programming equation.

\begin{theorem}
  \label{gdp:pr:DP_withoutstate_third}
  We suppose to be in the setting of~\S\ref{The_Bertsekas-Shreve_setting}.
  For $t$ in $\ic{0,\horizon{-}1}$, we define the \emph{Bellman operator}
  $\BelOp{t{+}1}{t}$ by, for all
  $\varphi\in \espace{L}^{0}_{+}(\HISTORY_{t{+}1})$
  and for all~$\history_{t} \in \HISTORY_{t}$,
  \begin{subequations}
    \begin{equation}
      \bp{ \BelOp{t{+}1}{t}\varphi }\np{\history_{t}} =
      \inf_{\control_{t}\in\CONTROL_{t}} \int_{\UNCERTAIN_{t{+}1}}
      \varphi\np{\history_{t},\control_{t},\uncertain_{t{+}1}}
      \StochKernel{t}{t{+}1} \kernelargs{\history_{t}}{d\uncertain_{t{+}1}}
      \eqfinp
      \label{gdp:eq:Bellman_operators_rho}
    \end{equation}
    Then, the Bellman operators are such that
    \begin{equation}
      \BelOp{t{+}1}{t} : \espace{L}^{0}_{+}(\HISTORY_{t{+}1})
      \to \espace{L}^{0}_{+}(\HISTORY_{t})
      \eqfinv
      \label{gdp:eq:Bellman_operators_rho_lower_semianalytic_preservation}
    \end{equation}
  \end{subequations}
  and the value functions~$\Value_{t}$ defined in~\eqref{gdp:eq:value_functions}
  are lower semianalytic and satisfy the \emph{Bellman equation},
  or \emph{(stochastic) dynamic programming equation},
  \begin{equation}
    \Value_{\horizon} =\criterion  \eqsepv
    \Value_{t} =\BelOp{t{+}1}{t}\Value_{t{+}1} \eqsepv
    \mtext{ for } t\in\ic{0,\horizon\!-\!1}
    \eqfinp
    \label{gdp:eq:Bellman_equation}
  \end{equation}
\end{theorem}

The proof is sketched in Appendix~\ref{Technical_Details_and_Proofs}.
This theorem is inspired
by \cite[Chap.~8]{Bertsekas-Shreve:1996},
with the feature that the state~$\optstate_t$ is, in our case,
the canonical history~$\history_t$, with the canonical dynamics
$ \history_{t{+}1} = \bp{\history_{t},\control_{t},\uncertain_{t{+}1}} $.
This quite general dynamic programming result is the basis of
all future developments in this paper.
Although the recalls and statements presented in this
Sect.~\ref{gdp:Stochastic_Dynamic_Programming_and_State_Reduction_by_Time_Blocks}
are mostly  straightforward consequences of results already established
in the literature, the developments are indispensable
to tackle time block decomposition in the forthcoming
Sect.~\ref{gdp:State_Reduction_by_Time_Blocks}.

\section{State reduction by time blocks and dynamic programming}
\label{gdp:State_Reduction_by_Time_Blocks}

In standard approaches to solve, by dynamic programming, a stochastic
optimal control problem formulated in discrete time, either a state
is given for all times (as in \cite{Bertsekas-Shreve:1996},
\cite{Hernandez-Lerma-Lasserre:1996}, \cite{Puterman:1994}
and~\cite{Witsenhausen:1973}), or no state is given
(as in \cite{Evstigneev:1976},\cite{Yuksel:2020}).
In this paper, our approach is intermediate, in that a state is possibly
obtained, but only at certain times. Thus, in this section, we consider
the question of reducing the history using a compressed ``state'' variable.
Differing with traditional practice, such a variable may not be available
at any time~$t\in\ic{0,T}$, but at some specified times
$ 0 = t_{0}< \cdots <t_{N} = \horizon $. We have recalled in
Sect.~\ref{gdp:Stochastic_Dynamic_Programming_and_State_Reduction_by_Time_Blocks}
that the history~$\history_{t}$
is itself a state variable with associated canonical dynamics
$ \history_{t{+}1} = \bp{\history_{t},\control_{t},\uncertain_{t{+}1}} $.
However, the size of this canonical state increases with time~$t$,
which is an unpleasant feature for dynamic programming --- quickly
leading to the well-known curse of dimensionality --- hence
the practical need to introduce a (ideally low dimensional) state
space, at least at some specified times, as done in this paper.
As already said in the introduction, the main difficulty
in achieving this goal is notational.

In~\S\ref{State_reduction_on_a_single_time_block},
we start by introducing the notion of state reduction on a single time block.
In~\S\ref{State_reduction_on_multiple_consecutive_time_blocks_and_dynamic_programming_equations},
we move to state reduction on multiple consecutive time blocks and we give
the corresponding dynamic programming equations across time blocks.
In~\S\ref{State_reduction_on_multiple_consecutive_time_blocks_and_reduced_optimal_feedbacks},
we conclude on how we obtain reduced optimal feedbacks.

\subsection{State reduction on a single time block}
\label{State_reduction_on_a_single_time_block}

We first present the case where the reduction only occurs at two
times denoted by~$\tun$ and~$\tter$, and such
that $  0 \leq \tun < \tter \leq \horizon $.

\begin{definition}
  \label{gdp:de:reduction-dynamics}
  Let~$\STATE_{\tun}$ and $\STATE_{\tter}$ be two Borel \emph{state spaces},
  $\theta_{\tun}$ and~$\theta_{\tter}$ be two Borel-measurable \emph{reduction mappings}
  \begin{subequations}
    \begin{equation}
      \theta_{\tun} : \HISTORY_{\tun} \to \STATE_{\tun} \eqsepv
      \theta_{\tter} : \HISTORY_{\tter} \to \STATE_{\tter} \eqfinv
      \label{gdp:eq:reduction_mappings}
    \end{equation}
    and $\Dynamics{\tun}{\tter}$ be a Borel-measurable \emph{dynamics}
    \begin{equation}
      \Dynamics{\tun}{\tter} : \STATE_{\tun} \times \HISTORY_{\tun{+}1:\tter}
      \to \STATE_{\tter} \eqfinp
    \end{equation}
    The triplet~$\np{\theta_{\tun},\theta_{\tter},\Dynamics{\tun}{\tter}}$ is
    called a \emph{state reduction across~$\icinterval{\tun}{\tter}$} if we
    have\footnote{%
      Notice that, if only the couple~$\np{\theta_{\tun},\Dynamics{\tun}{\tter}}$ is
      given, we can define the reduction mapping~$\theta_{\tter}$ by~\eqref{gdp:eq:reduction-dynamics},
      and thus obtain a
      triplet~$\np{\theta_{\tun},\theta_{\tter},\Dynamics{\tun}{\tter}}$
      which is a state reduction across~$\icinterval{\tun}{\tter}$.}
    \begin{equation}
      \theta_{\tter}\bp{\np{\history_{\tun}, \history_{\tun{+}1:\tter} }} =
      \Dynamics{\tun}{\tter} \bp{ \theta_{\tun}\np{\history_{\tun}},\history_{\tun{+}1:\tter} }
      \eqsepv \forall \history_{\tter} \in \HISTORY_{\tter}  \eqfinp
      \label{gdp:eq:reduction-dynamics}
    \end{equation}
  \end{subequations}
  The state reduction~$\np{\theta_{\tun},\theta_{\tter},\Dynamics{\tun}{\tter}}$
  is said to be \emph{compatible} with the
  sequence~$\na{\StochKernel{\tbis-1}{\tbis}}_{\tun{+}1 \leq \tbis \leq \tter}$
  of Borel-measurable stochastic kernels~\eqref{gdp:eq:sequence_of_stochastic_kernels} if
  \begin{itemize}
    \item
      there exists a Borel-measurable \emph{reduced stochastic kernel}
      $
      \tildeStochKernel{\tun}{\tun{+}1} : \STATE_{\tun} \to
      \Delta\np{\UNCERTAIN_{\tun{+}1}} $,
      such that the stochastic kernel~$\StochKernel{\tun}{\tun{+}1}$ in~\eqref{gdp:eq:sequence_of_stochastic_kernels}
      can be factored, for all $ \history_{\tun} \in \HISTORY_{\tun} $, as
      $
      \StochKernel{\tun}{\tun{+}1}\kernelargs{\history_{\tun}}{\dd\uncertain_{\tun{+}1}}$
      $=
      \tildeStochKernel{\tun}{\tun{+}1}
      \bkernelargs{\theta_{\tun}\np{\history_{\tun}}}{\dd\uncertain_{\tun{+}1}} $,
    \item
      for all $ \tbis$ in $\ic{\tun{+}2,\tter} $,
      there exists a Borel-measurable \emph{reduced stochastic kernel}
      $
      \tildeStochKernel{\tbis-1}{\tbis} :
      \STATE_{\tun} \times \HISTORY_{\tun{+}1:\tbis-1}
      \to \Delta\np{\UNCERTAIN_{\tbis}} $,
      such that the stochastic kernel~$\StochKernel{\tbis-1}{\tbis}$
      can be factored, for all $\history_{\tbis-1} \in \HISTORY_{\tbis-1}$, as
      $
      \StochKernel{\tbis-1}{\tbis}
      \bkernelargs{\np{\history_{\tun},\history_{\tun{+}1:\tbis-1}}}{\dd\uncertain_{\tbis}}=$
      $\tildeStochKernel{\tbis-1}{\tbis}
      \bkernelargs{\bp{\theta_{\tun}\np{\history_{\tun}},\history_{\tun{+}1:\tbis-1}}}
                  {\dd\uncertain_{\tbis}}$.
  \end{itemize}
\end{definition}

The above definition is similar to the sufficient statistics
idea in stochastic control: the state variable, which summarizes
the history, is sufficient for the controller to design
its control policy (\cite[p.~19]{Whittle:1982-I},
\cite[Definition~10.6]{Bertsekas-Shreve:1996},
\cite{Subramanian-Mahajan:2019}).
However, sufficient statistics in the stochastic control literature
are defined at the original time stages. By contrast,
Definition~\ref{gdp:de:reduction-dynamics}
--- and the coming Definition~\ref{gdp:de:reduction-dynamics_family} ---
consider a notion of sufficient statistics
\emph{only for a subset of time stages}.

\begin{figure}[htb]
  \begin{center}
    \includegraphics[width=\textwidth]{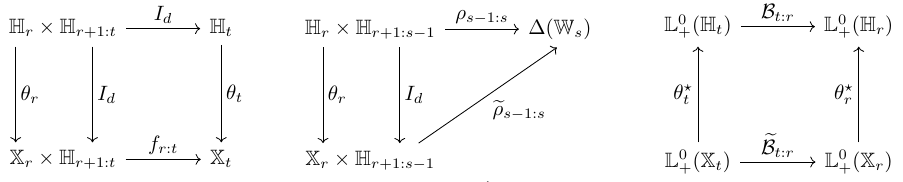}
  \end{center}
  \caption{Commutative diagrams in case of state reduction
    \label{gdp:fig:state_reduction}}
\end{figure}

According to Definition~\ref{gdp:de:reduction-dynamics},
the triplet~$\np{\theta_{\tun},\theta_{\tter},\Dynamics{\tun}{\tter}}$
is a state reduction across~$\icinterval{\tun}{\tter}$ if and only if
the diagram in the left part of Figure~\ref{gdp:fig:state_reduction}
is commutative; it is compatible if and only if the diagram in the
middle part of Figure~\ref{gdp:fig:state_reduction} is commutative.

The following proposition is the key ingredient to formulate
dynamic programming equations with a reduced state.

\begin{proposition}
  \label{gdp:thm:DPB}
  % \mdl{ce qui suit devrait être mis en fin de la section~2 :}
  Under the assumptions of \S\ref{The_Bertsekas-Shreve_setting},
  we define the \emph{Bellman
  operator across~$\icinterval{\tter}{\tun}$},
  $\BelOp{\tter}{\tun} : \espace{L}^{0}_{+}(\HISTORY_{\tter})
  \to \espace{L}^{0}_{+}(\HISTORY_{\tun})$ by
  \begin{equation}
    \BelOp{\tter}{\tun} = \BelOp{\tun{+}1}{\tun}  \circ \cdots \circ  \BelOp{\tter}{\tter-1} \eqsepv
    \label{gdp:eq:BelOp_t:r}
  \end{equation}
  where the one time step operators~$ \BelOp{\tbis}{\tbis-1} $,
  for $ \tbis $ in $ \ic{\tun{+}1,\tter} $ are defined
  in~\eqref{gdp:eq:Bellman_operators_rho}.

  Suppose that there exists a state reduction
  $\np{\theta_{\tun},\theta_{\tter},\Dynamics{\tun}{\tter}}$
  that is compatible with the sequence
  $\na{\StochKernel{\tbis-1}{\tbis}}_{\tbis\in\ic{\tun{+}1,\tter}}$
  of stochastic kernels~\eqref{gdp:eq:sequence_of_stochastic_kernels}.
  Then, there exists a
  \emph{reduced Bellman operator across~$\icinterval{\tter}{\tun}$},
  $ \tildeBelOp{\tter}{\tun} : \espace{L}^{0}_{+}(\STATE_{\tter})
  \to \espace{L}^{0}_{+}(\STATE_{\tun}) $, % \eqfinv
  such that
  \begin{equation}
    \bp{\tildeBelOp{\tter}{\tun} \tilde\varphi_{\tter}} \circ \theta_{\tun} =
    \BelOp{\tter}{\tun} \np{\tilde\varphi_{\tter} \circ \theta_{\tter}}
    \eqsepv \forall \tilde\varphi_{\tter} \in \espace{L}^{0}_{+}(\STATE_{\tter})
    \eqfinp
    \label{gdp:eq:DPB}
  \end{equation}
  For
  any~$\tilde\varphi_{\tter} \in \espace{L}^{0}_{+}(\STATE_{\tter})$
  and for any $\optstate_{\tun}\in\STATE_{\tun}$, we have that
  \begin{align}
    \big(\tildeBelOp{\tter}{\tun}\tilde\varphi_{\tter}\big)
           & =
             \inf_{\control_{\tun}\in\CONTROL_{\tun}} \int_{\UNCERTAIN_{\tun{+}1}}
             \tildeStochKernel{\tun}{\tun{+}1}
             \kernelargs{\optstate_{\tun}}{\dd\uncertain_{\tun{+}1}}
             \nonumber \\
           &
             \qquad
             \inf_{\control_{\tun{+}1}\in\CONTROL_{\tun{+}1}} \int_{\UNCERTAIN_{\tun{+}2}}
             \tildeStochKernel{\tun{+}1}{\tun{+}2}
             \kernelargs{\optstate_{\tun},\control_{\tun},\uncertain_{\tun{+}1}}{\dd\uncertain_{\tun{+}2}}
             \quad \cdots
             \nonumber    \\
           & \qquad\qquad
             \inf_{\control_{\tter-1}\in\CONTROL_{\tter-1}} \int_{\UNCERTAIN_{\tter}}
             \tildeStochKernel{\tter-1}{\tter}
             \kernelargs{\optstate_{\tun},
             \control_{\tun},\uncertain_{\tun{+}1},\ldots,\control_{\tter-2},\uncertain_{\tter-1}}
             {\dd\uncertain_{\tter}}
             \nonumber \\
           & \qquad\qquad\qquad\qquad\qquad\qquad\qquad
             \tilde\varphi_{\tter}\bp{\Dynamics{\tun}{\tter}
             \np{\optstate_{\tun}, \control_{\tun},\uncertain_{\tun{+}1},\ldots,
             \control_{\tter-1},\uncertain_{\tter}}}
             \label{gdp:eq:tildeBelOp-expression} \eqfinp
  \end{align}
\end{proposition}
The formula~\eqref{gdp:eq:tildeBelOp-expression} represents a nested sequence of infima
of integrals (with respect to different stochastic kernels).

The proof of Proposition~\ref{gdp:thm:DPB} is given in
Appendix~\ref{Technical_Details_and_Proofs}.
Proposition~\ref{gdp:thm:DPB} can be interpreted as follows. Denoting by
$\theta_{\tter}^{\star} : \espace{L}^{0}_{+}(\STATE_{\tter})
\to \espace{L}^{0}_{+}(\HISTORY_{\tter})$
the operator defined by
$
\theta_{\tter}^{\star}(\tilde\varphi_{\tter}) =
\tilde\varphi_{\tter} \circ \theta_{\tter}$
for any $ \tilde\varphi_{\tter} \in \espace{L}^{0}_{+}(\STATE_{\tter})
$,
the relation~\eqref{gdp:eq:DPB} rewrites as
$
\theta_{\tun}^{\star} \circ \tildeBelOp{\tter}{\tun} =
\BelOp{\tter}{\tun} \circ \theta_{\tter}^{\star} $,
that is, Proposition~\ref{gdp:thm:DPB} states that the diagram in
the right part of Figure~\ref{gdp:fig:state_reduction} is commutative.

\subsection{State reduction on multiple consecutive time blocks and dynamic programming equations}
\label{State_reduction_on_multiple_consecutive_time_blocks_and_dynamic_programming_equations}

Proposition~\ref{gdp:thm:DPB} can easily be extended to the case of
multiple consecutive time blocks $ \ic{t_{i},t_{i{+}1}} $,
with $N\in\NN^*$, $i \in \ic{0,N{-}1}$ and
$ 0 = t_{0}< \cdots <t_{N} = \horizon $.

\begin{definition}
  Let~$\nseqa{\STATE_{t_{i}}}{i\in \ic{0,N}}$
  be a family of Borel state spaces,
  $\nseqa{\theta_{t_{i}}}{i\in \ic{0,N}}$ be a family of Borel-measurable reduction
  mappings $ \theta_{t_{i}} : \HISTORY_{t_{i}} \to \STATE_{t_{i}} $,
  and $\nseqa{ \Dynamics{t_{i}}{t_{i{+}1}} }{i \in \ic{0,N{-}1}}$ be a family of
  Borel-measurable dynamics
  \[
    \Dynamics{t_{i}}{t_{i{+}1}} : \STATE_{t_{i}} \times \HISTORY_{t_{i}{+}1:t_{i{+}1}}
    \to \STATE_{t_{i{+}1}}
    \eqfinp
  \]
  The family
  $\big( \nseqa{ \STATE_{t_{i}} }{i\in \ic{0,N}}$,
  $\nseqa{\theta_{t_{i}}}{i\in \ic{0,N}}$,
  $\nseqa{ \Dynamics{t_{i}}{t_{i{+}1}} }{i \in \ic{0,N{-}1}}\big)$
  is called a \emph{state reduction across the consecutive time blocks
  $\ic{t_{i},t_{i{+}1}}$, $i \in \ic{0,N{-}1}$} if every triplet
  $\np{ \theta_{t_{i}} , \theta_{t_{i{+}1}} , \Dynamics{t_{i}}{t_{i{+}1}} }$
  is a state reduction across~$\icinterval{t_{i}}{t_{i{+}1}}$,
  for $i$ in $\ic{0,N{-}1}$.

  The state reduction across the consecutive time blocks
  $\ic{t_{i},t_{i{+}1}}$ is said to be \emph{compatible} with the
  family~$\na{\StochKernel{\tbis-1}{\tbis}}_{\tbis\in\ic{1,\horizon}}$
  of stochastic kernels given in~\eqref{gdp:eq:sequence_of_stochastic_kernels}
  if every triplet $ \np{ \theta_{t_{i}} , \theta_{t_{i{+}1}} ,
    \Dynamics{t_{i}}{t_{i{+}1}} } $ is compatible with the
  family~$\na{\StochKernel{\tbis-1}{\tbis}}_{\tbis\in\ic{t_{i}{+}1,t_{i{+}1}}}$,
  for $i$ in $\ic{0,N{-}1}$.
  \label{gdp:de:reduction-dynamics_family}
\end{definition}

\begin{remark}{\bf (Composed state dynamics as a reduction mapping)}
  \label{rem:composed-state-dynamics-as-reduction-mapping}

  There is a practical case where state reductions can readily be
  obtained, namely, when the model is given by controlled state
  dynamics driven by noises.  In that case, we are given a sequence
  $\nseqa{ \STATE_{\tslow}}{\tslow \in\ic{0,T}}$ of Borel state spaces
  and a sequence $ \nseqa{\Dynamics{\tslow}{\tslow{+}1}}{\tslow\in \ic{0,\horizon{-}1}} $
  of Borel-measurable dynamics
  \begin{equation}
    \Dynamics{\tslow}{\tslow{+}1} : \STATE_{\tslow}\times \CONTROL_{\tslow} \times \UNCERTAIN_{\tslow{+}1}
    \to \STATE_{\tslow{+}1}
    \eqfinp
  \end{equation}
  For any time $s\in \ic{0,\horizon{-}1}$, we define the composition
  $\Dynamics{0}{\tslow{+}1}=
  \Dynamics{\tslow}{\tslow{+}1}\circ \Dynamics{\tslow-1}{\tslow}\circ \ldots \circ
  \Dynamics{0}{1}$ with the abuse of notation that the composition is performed on the
  state argument. Setting $\UNCERTAIN_{0}=\STATE_{0}$, we obtain that
  $\Dynamics{0}{\tslow{+}1}:\HISTORY_{\tslow{+}1}\to\STATE_{\tslow{+}1}$ is a Borel-measurable mapping from
  the history space~$\HISTORY_{\tslow{+}1}$ taking values in the state
  space~$\STATE_{\tslow{+}1}$.

  Now, given a natural number $N>0$ and an increasing
  sequence $ 0 = t_{0}< \cdots <t_{N} = \horizon $ of times,
  we define the sequence $\bseqa{\theta_{t_{i}}}{i\in \ic{0,N}}$
  of Borel-measurable reduction mappings
  by $\theta_{t_i}= \Dynamics{0}{t_i}: \HISTORY_{t_i}\to\STATE_{t_i}$ for $i > 0$,
  and by $\theta_0=I_d$ (the identity mapping on $\UNCERTAIN_{0}$) for $i=0$.
  Moreover, given $i$ and $j \in \ic{0,N}$, with $i < j$ we obtain,
  for all $\history_{t_j} \in \HISTORY_{t_j}$,  that
  \begin{equation}
    \theta_{t_j}\np{\history_{t_j}}=
    \theta_{t_j}\bp{\np{\history_{t_i}, \history_{t_i{+}1:t_j} }} =
    \Dynamics{t_i}{t_j} \bp{ \theta_{t_i}\np{\history_{t_i}},\history_{t_i{+}1:t_j}}
    \eqfinv
  \end{equation}
  with
  $\Dynamics{t_i}{t_j}= \Dynamics{t_j-1}{t_j}\circ
  \Dynamics{t_j-2}{t_j-1}\circ \ldots \circ \Dynamics{t_i}{t_i{+}1}$
  which gives the state reduction Equation~\eqref{gdp:eq:reduction-dynamics}.
\end{remark}

\begin{remark}{\bf (Block independent exogenous noises and stochastic kernels)}
  \label{rem:block-independent-noises-and-kernel}

  There is a practical case where \emph{compatible} state
  reductions can readily be obtained. Assume that the
  sequence~$\na{\StochKernel{\tbis-1}{\tbis}}_{\tbis\in\ic{1,\horizon}}$
  of stochastic kernels
  in~\eqref{gdp:eq:sequence_of_stochastic_kernels}
  are mappings whose arguments do not include the control part
  (that is, depend at most on the history uncertainty part
  (see~\eqref{gdp:eq:uncertain_t:r}).
  If we interpret stochastic kernels as (conditional) distributions of noises
  (random process), this means that the system dynamics are
  driven by an exogenous noise process, say $\nseqa{\va{W}_{t}}{t\in\ic{1,T}}$.

  Assume moreover that the stochastic kernels give rise to noises that
  are independent block by block, in the sense that the random vectors
  \( \va{W}_{0} \),
  \( \nseqp{\va{W}_{t}}{t\in\ic{1,t_{1}}} \),
  \( \nseqp{\va{W}_{t}}{t\in\ic{t_1{+}1,t_{2}}} \),
  \ldots,
  \( \nseqp{\va{W}_{t}}{t\in\ic{t_i{+}1,t_{i{+}1}}} \),
  \ldots,
  \( \nseqp{\va{W}_{t}}{t\in\ic{t_{N{-}2}{+}1,t_{N{-}1}}} \),
  \( \nseqp{\va{W}_{t}}{t\in\ic{t_{N{-}1}{+}1,t_{N}}} \)
  are stochastically independent.
  Then, from Definitions~\ref{gdp:de:reduction-dynamics} and
  \ref{gdp:de:reduction-dynamics_family}, we deduce that any
  state reduction across the same time blocks is compatible
  with the stochastic kernels.
\end{remark}

Assuming the existence of a state reduction across the consecutive
time blocks~$\ic{t_{i},t_{i{+}1}}$ compatible with the sequence of stochastic
kernels~\eqref{gdp:eq:sequence_of_stochastic_kernels}, we obtain the existence
of a \emph{sequence of reduced Bellman operators} across the
time blocks~$\icinterval{t_{i}}{t_{i{+}1}}$ as an immediate consequence
of multiple applications of Proposition~\ref{gdp:thm:DPB}, that is,
$
\tildeBelOp{t_{i{+}1}}{t_{i}} :
\espace{L}^{0}_{+}(\STATE_{t_{i{+}1}})
\to \espace{L}^{0}_{+}(\STATE_{t_{i}})
$, $i\in \ic{0,N{-}1} $,
such that, for any function~$\tilde\varphi_{t_{i{+}1}} \in
\espace{L}^{0}_{+}(\STATE_{t_{i{+}1}}) $,
we have that
$\bp{\tildeBelOp{t_{i{+}1}}{t_{i}} \tilde\varphi_{t_{i{+}1}}} \circ \theta_{t_{i}} =
\BelOp{t_{i{+}1}}{t_{i}} \np{\tilde\varphi_{t_{i{+}1}} \circ \theta_{t_{i{+}1}}} $.
We now consider the family of optimization problems defined by
the associated value functions~\eqref{gdp:eq:value_functions}.
Thanks to the state reductions, we can enounce the following
two theorems which establish dynamic programming equations
\emph{across} consecutive time blocks.
The first one,
Theorem~\ref{gdp:thm:DPB_family}, states a dynamic programming equation
for an optimization problem
in Mayer form (that is, just involving a final cost). The second
one, Theorem~\ref{cor:instantaneous-costs-in-addition-to-final-cost},
is more general as it involves both instantaneous costs and
a final cost. As it is well known that the second case can be
reduced to a Mayer form through a state augmentation, the proof
of Theorem~\ref{cor:instantaneous-costs-in-addition-to-final-cost}
easily follows from the proof of Theorem~\ref{gdp:thm:DPB_family}.

\begin{theorem}[Time block decomposition for the Mayer form]
  \label{gdp:thm:DPB_family}
  We assume to be in the setting of~\S\ref{The_Bertsekas-Shreve_setting}.
  Suppose that a state reduction
  $\big( \nseqa{\STATE_{t_{i}}}{i\in \ic{0,N}}$,
  $ \nseqa{\theta_{t_{i}}}{i\in \ic{0,N}}$,
  $\nseqa{ \Dynamics{t_{i}}{t_{i{+}1}}}{i \in \ic{0,N{-}1}}\big)$
  exists across the consecutive time blocks $\nseqa{\ic{t_{i},t_{i{+}1}}}{i\in \ic{0,N{-}1}}$,
  satisfying $0 = t_{0}< \cdots <t_{N} = \horizon $,
  which is compatible with the
  sequence~$\na{\StochKernel{\tbis-1}{\tbis}}_{\tbis\in\ic{1,\horizon}}$
  of stochastic kernels given in~\eqref{gdp:eq:sequence_of_stochastic_kernels}.

  Suppose that there exists a \emph{reduced cost criterion}
  $ \tildecriterion : \STATE_{\horizon} \to \bRR $,
  which is a nonnegative lower semianalytic function and is
  such that the cost function~$\criterion$ in~\eqref{gdp:eq:criterion}
  can be factored as $\criterion = \tildecriterion \circ \theta_{\horizon}$.
  We define the sequence of \emph{reduced value functions}
  $\{\tilde\Value_{t_{i}}\}_{i\in \ic{0,N}}$, where \( \tilde\Value_{t_{i}} :
  \STATE_{t_{i}} \to \bRR \) for $i \in \ic{0,N}$, by
  \begin{equation}
    \tilde\Value_{t_{N}} = \tildecriterion
    \mtext{ and }
    \tilde\Value_{t_{i}} = \tildeBelOp{t_{i{+}1}}{t_{i}} \tilde\Value_{t_{i{+}1}}
    \eqsepv \forall i\in\ic{0,N{-}1}
    \eqfinv
    \label{gdp:eq:reduced_value_functions}
  \end{equation}
  where the reduced Bellman operators~$\nseqa{\tildeBelOp{t_{i{+}1}}{t_{i}}}{i\in \ic{0,N{-}1}}$ across the
  intervals~$\nseqa{\ic{t_{i},t_{i{+}1}}}{i\in \ic{0,N{-}1}}$
  are given in~\eqref{gdp:eq:tildeBelOp-expression}.
  Then, the sequence $\{\Value_{t_{i}}\}_{i\in \ic{0,N}}$
  in~\eqref{gdp:eq:value_functions} satisfies
  $ \Value_{t_{i}} =  \tilde\Value_{t_{i}} \circ \theta_{t_{i}}$,
  for all $i \in \ic{0,N} $.
\end{theorem}
\begin{proof}
  The proof is an immediate consequence of multiple applications
of Theorem~\ref{gdp:pr:DP_withoutstate_third} and
Proposition~\ref{gdp:thm:DPB}.
\end{proof}

Finally, we consider the special case where the criterion
$\criterion:\HIST[\horizon] \to \bRR$ is factored as
\begin{equation}
\criterion\np{\history_\horizon} = \sum_{i=0}^{N{-}1}
\InstantneousCost_{t_i} \bp{ \theta_{t_i}\np{\history_{t_i}},\history_{t_i{+}1:t_{i{+}1}}} +
\FinalCost_{t_N}\bp{\theta_{t_N}\np{\history_{t_N}}}
\eqfinv
\label{gdp:eq:criterion_time_block_additive}
\end{equation}
where the numerical functions $\nseqa{\InstantneousCost_{t_{i}}}{i\in \ic{0,N}}$
are nonnegative lower semianalytic, with
$\InstantneousCost_{t_{i}} : \STATE_{t_{i}}\times\HISTORY_{t_i{+}1:t_{i{+}1}} \to \bRR$
for $i\in \ic{0,N}$. The associated optimization problems,
indexed by $i \in \ic{0,N{-}1}$ and parameterized
by~$\history_{t_{i}} \in \HISTORY_{t_{i}}$, are given by
\begin{subequations}
\begin{align}
\Value_{t_{i}}(\history_{t_{i}})
 & = \inf_{\HistoryFeedback_{t_{i}:\horizon{-}1} \in
   \HISTORYFEEDBACK_{t_{i}:\horizon{-}1}} \int_{\HISTORY_{\horizon}}
   \Bp{ \sum_{j=i}^{N{-}1} \InstantneousCost_{t_j}
        \bp{\theta_{t_j}\np{\history'_{t_j}},\history'_{t_j{+}1:t_{j{+}1}}} +
        \FinalCost_{t_N}\bp{\theta_{t_N}\np{\history'_{t_N}}} }
   \StochKernelFeed{t_{i}}{\horizon}{\HistoryFeedback}
   \kernelargs{\history_{t_{i}}}{\dd \history'_{\horizon}} \eqfinv
\intertext{and, for~$i=N$,}
\Value_{t_{N}}(\history_{t_{N}})
 & = \FinalCost_{t_N}\bp{\theta_{t_N}\np{\history_{t_N}}} \eqfinp
\end{align}
\label{gdp:eq:value_functions_additive}
\end{subequations}
These Bellman equations are a special case
of Equations~\eqref{gdp:eq:value_functions} when
the cost criterion~$\criterion$ is given
by~\eqref{gdp:eq:criterion_time_block_additive}.
It is left to the reader
to prove that the following theorem holds true\footnote{%
The proof uses \cite[Lemma~7.30~(3,4)]{Bertsekas-Shreve:1996}
on the stability of lower semianalytic numerical functions under addition
and under right composition with a Borel-measurable mapping.}.

\begin{theorem}{\bf (Taking care of instantaneous costs in addition to final cost)}
  \label{cor:instantaneous-costs-in-addition-to-final-cost}

  Suppose that the assumptions of Theorem~\ref{gdp:thm:DPB_family} are
  satisfied, but for the cost criterion $\criterion:\HIST[\horizon] \to \bRR$
  defined by Equation~\eqref{gdp:eq:criterion_time_block_additive}.

  We define the sequence of \emph{reduced value functions}
  $\{\tilde\Value_{t_{i}}\}_{i\in \ic{0,N}}$, where \( \tilde\Value_{t_{i}} :
  \STATE_{t_{i}} \to \bRR \) for $i \in \ic{0,N}$, by
  \begin{equation}
    \tilde\Value_{t_{N}} =\FinalCost_{t_N}
    \mtext{ and }
    \tilde\Value_{t_{i}} = \overlineBelOp{t_{i{+}1}}{t_{i}} \tilde\Value_{t_{i{+}1}}
    \eqsepv \forall i\in\ic{0,N{-}1}
    \eqfinv
    \label{gdp:eq:reduced_value_functions_time_block_additive}
  \end{equation}
  where the reduced Bellman operator~$\overlineBelOp{t_{i{+}1}}{t_{i}}$
  across~$\icinterval{t_i}{t_{i{+}1}}$ are given,
  for any $i\in\ic{0,N{-}1}$, for any~$\tilde\varphi_{t_{i{+}1}} \in
  \espace{L}^{0}_{+}(\STATE_{t_{i{+}1}})$
  and for any $\optstate_{t_{i}}\in\STATE_{t_{i}}$, by
  \begin{align}
    \big( \overlineBelOp{t_{i{+}1}}{t_{i}}\tilde\varphi_{t_{i{+}1}}\big) (\optstate_{t_{i}})
    =
    \inf_{\control_{t_{i}}\in\CONTROL_{t_{i}}}
    & \int_{\UNCERTAIN_{t_{i}{+}1}}
      \tildeStochKernel{t_{i}}{t_{i}{+}1}
      \kernelargs{\optstate_{t_{i}}}{\dd\uncertain_{t_{i}{+}1}}
      \nonumber \\
    &
      \inf_{\control_{t_{i}{+}1}\in\CONTROL_{t_{i}{+}1}} \int_{\UNCERTAIN_{t_{i}{+}2}}
      \tildeStochKernel{t_{i}{+}1}{t_{i}{+}2}
      \kernelargs{\optstate_{t_{i}},\control_{t_{i}},\uncertain_{t_{i}{+}1}}{\dd\uncertain_{t_{i}{+}2}}
      \quad \cdots
      \nonumber    \\
    & \quad
      \inf_{\control_{t_{i{+}1}-1}\in\CONTROL_{t_{i{+}1}-1}} \int_{\UNCERTAIN_{t_{i{+}1}}}
      \tildeStochKernel{t_{i{+}1}-1}{t_{i{+}1}}
      \nonumber    \\
    &\quad\quad\quad
      \kernelargs{\optstate_{t_{i}},
      \control_{t_{i}},\uncertain_{t_{i}{+}1},\ldots,\control_{t_{i{+}1}-2},\uncertain_{t_{i{+}1}-1}}
      {\dd\uncertain_{t_{i{+}1}}}
      \nonumber \\
    & \quad\quad
      \Big(
      \InstantneousCost_{t_i}
      \np{\optstate_{t_{i}}, \control_{t_{i}},\uncertain_{t_{i}{+}1},\ldots,\control_{t_{i{+}1}-1},\uncertain_{t_{i{+}1}}}
      \nonumber  \\
    & \quad\quad\quad
      +
      \tilde\varphi_{t_{i{+}1}}\bp{ \Dynamics{t_{i}}{t_{i{+}1}}
      \np{\optstate_{t_{i}}, \control_{t_{i}},\uncertain_{t_{i}{+}1},\ldots,\control_{t_{i{+}1}-1},\uncertain_{t_{i{+}1}}}}
      \Big)
      \label{gdp:eq:overlineBelOp-expression}
      \eqfinp \!\!\!\!\!
  \end{align}
  Then, the sequence $\{\Value_{t_{i}}\}_{i\in \ic{0,N}}$
  in Equations~\eqref{gdp:eq:value_functions_additive} satisfies
  $ \Value_{t_{i}} =  \tilde\Value_{t_{i}} \circ \theta_{t_{i}}$,
  for all $i \in \ic{0,N} $.
\end{theorem}
Here again, Formula~\eqref{gdp:eq:overlineBelOp-expression} 
represents a nested sequence of infima of integrals
(with respect to different stochastic kernels).

Of course, solving Equation~\eqref{gdp:eq:reduced_value_functions} or
Equation~\eqref{gdp:eq:overlineBelOp-expression}
can be as difficult as solving the original Bellman equation.
However, the interest of such time block decomposition will be illustrated on
different applications in
Sect.~\ref{sec:Mixing_dynamic_programming_and_stochastic_programming},
Sect.~\ref{Two_time_scale_optimization_problem}
and Sect.~\ref{gdp:Decision_Hazard_Decision_Dynamic_Programming}.

\subsection{State reduction on multiple consecutive time blocks
            and reduced optimal feedbacks}
\label{State_reduction_on_multiple_consecutive_time_blocks_and_reduced_optimal_feedbacks}

As in the classical dynamic programming framework
\cite[p.~190]{Bertsekas-Shreve:1996},
we recover the property that the search of an optimal
policy among all policies (\emph{history} feedbacks) can be limited to
the search of an optimal \emph{state} feedback.
This is the most important result in practice.

\begin{proposition}
  Under the assumptions of Theorem~\ref{gdp:thm:DPB_family}, the reduced
  value functions $\{\tilde\Value_{t_{i}}\}_{i\in \ic{0,N}}$ defined
  in~\eqref{gdp:eq:reduced_value_functions} are equal to the minimum value
  of the following optimization problems, parameterized by the reduced
  history (state) $\state_{t_i}\in \STATE_{t_i}$
  \begin{subequations}
    \label{eq:reduced-value-function}
    \begin{align}
      \tilde\Value_{t_i}(\state_{t_i})
      & = \inf_{\HistoryFeedback_{t_i:\horizon{-}1} \in
        \HISTORYFEEDBACK^{\state_{t_i}}_{t_i:\horizon{-}1}}
        \int_{\HISTORY_{t_{i{+}1}:\horizon}}
        \barcriterion\np{\state_{t_i},\history'_{t_{i}{+}1:\horizon}}
        \tildevarStochKernelFeed{\tter_i}{\horizon}{\HistoryFeedback}
        \kernelargs{\state_{t_i},\history_{t_{i{+}1}:\tter}}{\dd \history'_{t_i{+}1:\horizon}}
        \eqsepv \forall \state_{\tter_i} \in \STATE_{\tter_i}
        \eqfinv \\
      \textrm{and } \quad \tilde\Value_{\horizon}(\state_{\horizon})
      & = \tildecriterion(\state_{\horizon})
        \eqsepv \forall \state_{\horizon}\in \STATE_{\horizon}
        \eqfinv
    \end{align}
  \end{subequations}
  where the mapping $\barcriterion$ is given by
  $\barcriterion = \tildecriterion \circ \Dynamics{t_{N-1}}{t_{N}}\circ
  \Dynamics{t_{N-2}}{t_{N-1}}\circ \ldots \circ \Dynamics{t_i}{t_{i{+}1}}$
  (with, as already noted, the abuse of notation that the composition
  is performed on the state argument),
  where $\tildevarStochKernelFeed{\tter_i}{\horizon}{\HistoryFeedback}$
  is the reduced stochastic kernel
  (see Definition~\ref{gdp:de:reduction-dynamics})
  associated with the kernel $\varStochKernelFeed{\tter_i}{\horizon}{\HistoryFeedback}$,
  the kernel $\varStochKernelFeed{\tter_i}{\horizon}{\HistoryFeedback}$
  being given in the factorization of the kernel
  $\StochKernelFeed{\tun}{\tter}{\HistoryFeedback}$, namely
  $\StochKernelFeed{\tun}{\tter}{\HistoryFeedback}
  \kernelargs{\history_{\tun}}{\dd\history'_{\tter}} = $
  $\delta_{\history_{\tun}}\np{ \dd\history'_{\tun} }
  \otimes  \varStochKernelFeed{\tun}{\tter}{\HistoryFeedback}
  \kernelargs{\history_{\tun}}{\dd\history'_{\tun{+}1:\tter}}$
  given by~\eqref{eq:kernelfactorization},
  $\delta$ being the Dirac measure,
  and where
  $\HISTORYFEEDBACK^{\state_{t_i}}_{t_i:\horizon{-}1}$ is
  the set of $(t_i{:}\horizon{-}1)$-reduced history feedbacks,
  that is, the set of sequences
  $\bseqa{\HistoryFeedback_{\tbis}}{\tbis\in\ic{t_i,\horizon-1}}$
  of universally measurable mappings
  $\HistoryFeedback_{\tbis} : \STATE_{t_i} \times \HISTORY_{t_i{+}1:\tbis}
  \to \CONTROL_{\tbis} $.
\end{proposition}

\begin{proof}
  Using Theorem~\ref{gdp:thm:DPB_family},
  we have that, for all $i \in \ic{0,N}$,
  $\Value_{t_{i}} = \tilde\Value_{t_{i}} \circ \theta_{t_{i}}$,
  with $\tilde\Value_{t_{i}}$ satisfying the Bellman equation~\eqref{gdp:eq:reduced_value_functions}.
  For establishing that $\tilde\Value_{t_{i}}$ is a value function
  satisfying Equation~\eqref{eq:reduced-value-function}, we now
  prove that, in the definition of $\Value_{t_{i}}$ in Equation~\eqref{gdp:eq:value_functions}, we can replace the
  space~$\HISTORYFEEDBACK_{t_i:\horizon{-}1}$ of history feedbacks
  by the space~$\HISTORYFEEDBACK^{\state_{t_i}}_{t_i:\horizon{-}1}$
  state feedbacks. We proceed as follows.
  Following~\cite[Chapter 8]{Bertsekas-Shreve:1996}, we use the Bellman
  equation~\eqref{gdp:eq:overlineBelOp-expression}
  to obtain $\epsilon$-minimizers
  for each problem~\eqref{gdp:eq:value_functions}.
  As $\epsilon$-minimizers are obtained by recursively solving
  Equations~\eqref{gdp:eq:overlineBelOp-expression}, they are obtained by
  solving (up to $\epsilon$) parametric optimization problems. Thus, we easily
  get, using~\cite[Proposition~8.3, p.~200]{Bertsekas-Shreve:1996}, that an
  $\epsilon$-minimizer at time $t\in \ic{t_i,t_{i{+}1}}$ is a universally
  measurable function of $(\theta(\history_{t_i}),\history_{t_{i{+}1}:t})$.  From
  this last fact, we get that --- in the value function definition
  of~$\Value_{t_{i}}$ given in Equation~\eqref{gdp:eq:value_functions} --- the
  space~$\HISTORYFEEDBACK_{0:\horizon{-}1}$ can be replaced by the space of
  feedbacks given by universally measurable functions of the ordered pair
  $(\theta(\history_{t_i}),\history_{t_{i{+}1}:t})$ without
  changing the value function.

  Finally, when considering Equation~\eqref{gdp:eq:value_functions} --- with
  this restricted space of state feedbacks,
  and considered at time~$t_i$ for $i \in \ic{0,N}$ ---
  we obtain that the cost to be minimized is now parameterized by $\theta(\history_{t_i})$ ---
  as it is the case for the cost to be integrated and
  also for the stochatic kernels induced by the state feedbacks.
  By setting \( \state_i=\theta(\history_{t_i}) \), the obtained  optimization
  problem is the right hand side of~\eqref{eq:reduced-value-function}, that we
  call \(  \tilde{\tilde\Value}_{t_i}(\state_{t_i}) \), and
  we have that, for all $i \in \ic{0,N}$,
  $\Value_{t_{i}} = \tilde{\tilde\Value}_{t_{i}} \circ \theta_{t_{i}}$.
  It remains to prove that \( {\tilde\Value}_{t_{i}}=
  \tilde{\tilde\Value}_{t_{i}}\), for all $i \in \ic{0,N}$.
  By a proof similar to the one of Theorem~\ref{gdp:thm:DPB_family}, we show
  that the sequence $\{\tilde{\tilde\Value}_{t_{i}}\}_{i\in \ic{0,N}}$
  satisfies the Bellman equation~\eqref{gdp:eq:reduced_value_functions}.
  By uniqueness, the sequence $\{\tilde{\tilde\Value}_{t_{i}}\}_{i\in \ic{0,N}}$
  coincides with the sequence $\{{\tilde\Value}_{t_{i}}\}_{i\in \ic{0,N}}$.
\end{proof}

\section{Mixing dynamic programming and stochastic programming}
\label{sec:Mixing_dynamic_programming_and_stochastic_programming}

As a first application of the formalism developed
in~\S\ref{State_reduction_on_multiple_consecutive_time_blocks_and_dynamic_programming_equations},
we show how dynamic programming and stochastic programming can be mixed
(which was the original motivation for the paper, see
Footnote~\ref{ft:Bogota2013}). The proof of the following proposition is a straightforward application of
Theorem~\ref{cor:instantaneous-costs-in-addition-to-final-cost} combined with Remark~\ref{rem:block-independent-noises-and-kernel}.

\begin{proposition}
  \label{prop:Stochprog_Dynprog}
  Suppose that the assumptions of Theorem~\ref{gdp:thm:DPB_family} are satisfied.
  We consider multiple consecutive time blocks $\ic{t_{i},t_{i{+}1}}$,
  with $N\in\NN^*$, $i \in \ic{0,N{-}1}$ and $0 = t_{0}< \cdots <t_{N} = \horizon$,
  and we assume that
  \begin{itemize}
  \item
    a state reduction
    $\big( \nseqa{\STATE_{t_{i}}}{i\in \ic{0,N}}$,
    $\nseqa{\theta_{t_{i}}}{i\in \ic{0,N}}$,
    $\nseqa{ \Dynamics{t_{i}}{t_{i{+}1}}}{i \in \ic{0,N{-}1}}\big)$
    exists across the consecutive time blocks
    $\nseqa{\ic{t_{i},t_{i{+}1}}}{i\in \ic{0,N{-}1}}$,
  \item
    the noises are exogeneous and \emph{time block independent},
    that is, the elements of the
    sequence~$\na{\StochKernel{\tbis-1}{\tbis}}_{\tbis\in\ic{1,\horizon}}$
    in~\eqref{gdp:eq:sequence_of_stochastic_kernels} are,
    for all $i \in \ic{0,N{-}1}$ and $\tun \in [t_i,t_{i{+}1})$, of the form
    \begin{equation}
      \label{block_indep_kernels}
      \StochKernel{\tun}{\tun{+}1}: \UNCERTAIN_{t_i} \times \cdots\times \UNCERTAIN_{\tun} \to \Delta\np{\UNCERTAIN_{\tun{+}1}}
      \eqfinv
    \end{equation}
    which means
    that the distribution of the uncertainty~$\uncertain_{\tun{+}1}$
    is only function of the past uncertainties~$(\uncertain_{t_i},\dots,\uncertain_{\tun})$
    within the time block,

  \item the cost criterion $\criterion:\HIST[\horizon] \to \bRR$
    can be factored as
    \begin{equation}
      \criterion\np{\history_\horizon}
      =
      \sum_{i=0}^{N{-}1}
      \InstantneousCost_{t_i} \bp{ \theta_{t_i}\np{\history_{t_i}},\history_{t_i{+}1:t_{i{+}1}}}
      +
      \FinalCost_{t_N}\bp{\theta_{t_N}\np{\history_{t_N}}}
      \eqfinv
      \label{gdp:eq:criterion_time_block_additive_two}
    \end{equation}
    where the numerical functions $\nseqa{\InstantneousCost_{t_{i}}}{i\in \ic{0,N}}$
    are nonnegative lower semianalytic, with
    \( \InstantneousCost_{t_{i}} : \STATE_{t_{i}} \times
    \HISTORY_{t_i{+}1:t_{i{+}1}}\to\bRR \)
    for $i\in \ic{0,N}$.
  \end{itemize}
  Then, the multistage stochastic optimization problem~\eqref{eq:Vzero}
  can be solved by the following algorithm.
  \begin{description}
  \item[Initialization.]
    Define $\tilde\Value_{t_{N}}=\FinalCost_{t_N} : \STATE_{\horizon} \to \bRR$.
  \item[Backward recursion.] Suppose that the function
    $\tilde\Value_{t_{i{+}1}} : \STATE_{t_{i{+}1}} \to \bRR$
    is known at index $i{+}1\in \ic{1,N}$.
    Then, for each state~$\optstate_{t_{i}}\in\STATE_{t_{i}}$ (for instance
    on a grid approximating the set~$\STATE_{t_{i}}$, or on~$\STATE_{t_{i}}$
    itself when finite and small enough), compute the previous
    Bellman value function~$\tilde\Value_{t_{i}}$ at index~$i$ as
    \begin{align}
      & \tilde\Value_{t_{i}}(\optstate_{t_{i}}) =
        \inf_{\control_{t_{i}}\in\CONTROL_{t_{i}}} \int_{\UNCERTAIN_{t_{i}{+}1}}
        \StochKernel{t_{i}}{t_{i}{+}1}
        \kernelargs{\uncertain_{t_{i}}}{\dd\uncertain_{t_{i}{+}1}}
        \nonumber \\
      & \qquad\qquad
        \inf_{\control_{t_{i}{+}1}\in\CONTROL_{t_{i}{+}1}} \int_{\UNCERTAIN_{t_{i}{+}2}}
        \StochKernel{t_{i}{+}1}{t_{i}{+}2}
        \kernelargs{\uncertain_{t_{i}},\uncertain_{t_{i}{+}1}}{\dd\uncertain_{t_{i}{+}2}}
        \quad \cdots
        \nonumber \\
      & \qquad\qquad
        \inf_{\control_{t_{i{+}1}-1}\in\CONTROL_{t_{i{+}1}-1}} \int_{\UNCERTAIN_{t_{i{+}1}}}
        \StochKernel{t_{i{+}1}-1}{t_{i{+}1}}
        \kernelargs{\uncertain_{t_{i}},\uncertain_{t_{i}{+}1},\ldots,\uncertain_{t_{i{+}1}-1}}
        {\dd\uncertain_{t_{i{+}1}}}
        \nonumber \\
      &
        \quad\quad\quad\quad\quad\quad\quad\quad\quad\quad\quad\quad
        \Big(
        \InstantneousCost_{t_i}
        \np{\optstate_{t_{i}}, \control_{t_{i}},\uncertain_{t_{i}{+}1},\ldots,\control_{t_{i{+}1}-1},\uncertain_{t_{i{+}1}}}
        \nonumber  \\
      &
        \quad\quad\quad\quad\quad\quad\quad\quad\quad\quad\quad\quad
        +
        \tilde\Value_{t_{i{+}1}}
        \bp{
        \Dynamics{t_{i}}{t_{i{+}1}}
        \np{\optstate_{t_{i}}, \control_{t_{i}},\uncertain_{t_{i}{+}1},\ldots,\control_{t_{i{+}1}-1},\uncertain_{t_{i{+}1}}}}
        \Big)
        \label{eq:Stochprog-Dynprog-recursion}
        \eqfinp % \!\!\!\!\!
    \end{align}
  \item[Final step.]
    Compute \( \Value_{0}(\uncertain_{0})
    = \tilde\Value_{t_{0}} \bp{ \theta_{t_{0}}\np{\uncertain_{t_{0}}} } \).
  \end{description}
\end{proposition}

In many practical situations, all the uncertainty
sets~$\UNCERTAIN_{0}$, \ldots, $\UNCERTAIN_{\horizon}$ are finite
and the computation in~\eqref{eq:Stochprog-Dynprog-recursion}
is tractable by using \emph{stochastic programming}
and \emph{scenario tree} techniques, which do not require
stagewise independence of the noises. We are thus able to
take advantage of both the dynamic programming world and
the stochastic programming world:
\begin{itemize}
\item
  use dynamic programming at a selection of time stages (for instance,
  at those of the slow time scale) and across the
  corresponding time blocks (for instance, across consecutive slow
  time stages), when noises are stochastically independent block
  by block; that yields Bellman value functions only for the
  chosen selection of time stages (for instance, at the slow time scale);
\item
  use stochastic programming inside time blocks (for instance,
  at fast time scale, within two consecutive slow time stages);
  the fast time scale final cost function of a block is given
  by the Bellman value function computed at the slow time scale
  which corresponds to the terminal time stage of the block; no stagewise
  independence assumption is required within time blocks
  (for instance, for the short time scale noises).
\end{itemize}

\begin{remark}
  As a special case,
  it is straightforward to check that the triplet
  $\big( \nseqa{\UNCERTAIN_{t_{i}}}{i\in \ic{0,N}}$,
  $\nseqa{\theta_{t_{i}}}{i\in \ic{0,N}}$,
  $\nseqa{ \Dynamics{t_{i}}{t_{i{+}1}} }{i \in \ic{0,N{-}1}}\big)$, with
  \begin{itemize}
  \item the reduction mapping
    $\theta_{t_{i}}$ given by
    $\theta_{t_{i}} ( \history_{t_{i}}) = \uncertain_{t_i}$ for all $i\in \ic{0,N}$,
  \item
    the  dynamics $\Dynamics{t_{i}}{t_{i{+}1}}$ given by
    $\Dynamics{t_{i}}{t_{i{+}1}}(\uncertain_{t_i},\history_{t_{i}{+}1:t_{i{+}1}}) =
    \uncertain_{t_{i{+}1}}$,  for all $i\in \ic{0,N-1}$.
  \end{itemize}
  is a state reduction across the consecutive time blocks
  $ \ic{t_{i},t_{i{+}1}} $, $i \in \ic{0,N{-}1}$ which is compatible
  with the sequence of stochastic kernels given
  by Equation~\eqref{block_indep_kernels}.
  Thus, Proposition~\ref{prop:Stochprog_Dynprog} applies.

  But, in this special case, the optimal controls can be computed
  in parallel with respect to time blocks, as the term
  $\tilde\Value_{t_{i{+}1}}\np{\uncertain_{t_{i{+}1}}}$ is
  a constant in \eqref{eq:Stochprog-Dynprog-recursion}.
  What is interesting in \eqref{eq:Stochprog-Dynprog-recursion}
  is the added fact that the optimal strategy which was, a priori,
  searched as feedbacks depending on the whole history is
  in fact made up of independent strategies, each defined on
  a single time block and made up of feedbacks depending only
  on the block history (the history within the block).
\end{remark}

\subsubsubsection{Numerical illustration}

To numerically illustrate the mixing between dynamic programming
and stochastic programming, we consider a toy optimization problem
over a time span~\( \ic{0,\horizon} \), where $\horizon$
is an even natural number
(for instance $\horizon=24$ for an hourly period problem during a day).
The problem involves a storage, the state~$\state_{t}$
of which is driven by a dynamics
involving a control variable~$\control_{t}$ and a noise variable~$\uncertain_{t{+}1}$.
We assume that the noises during the first half time span,
that is, for $t\in\ic{1,\horizon/2}$,
are independent of the noises during the second half time span,
that is, for $t\in\ic{\horizon/2,\horizon}$.
We also assume that each noise variable~$\uncertain_{t}$ can
only take two possibles values, so that the whole uncertainty
process can be represented by a binary tree.

In this problem, we consider the two consecutive time blocks~$\ic{0,\horizon/2}$
and~$\ic{\horizon/2,\horizon}$, and the state reduction is given
in a straightforward manner by the variable~$\state_{t}$ (as explained
in Remark~\ref{rem:composed-state-dynamics-as-reduction-mapping}).
Thus, we are able to compute Bellman functions by the algorithm given
in Proposition~\ref{prop:Stochprog_Dynprog}.
We illustrate the algorithm for the horizon~$\horizon=24$.
\begin{itemize}
\item
The Bellman function~$\tilde\Value_{24}$ is given
  by the final cost function of the problem.
\item
  The Bellman function~$\tilde\Value_{12}$ is
  approximated by discretization and it is computed on a grid
  involving~$n$ points~$(\state^{1},\dots,\state^{n})$.
  For \( i\in\ic{1,n} \), each value~$\tilde\Value_{12}(\state^{i})$
  is obtained by
  solving a stochastic programming problem on the time span
  $\ic{12,24}$, that is, on a tree involving $2^{12}$ leaves
(as each noise variable~$\uncertain_{t}$ can
only take two possibles values).
\item
  The optimal cost of the optimization problem
  is~$\tilde\Value_{0}(\state_{0})$, obtained again
  by stochastic programming on the time span~$\ic{0,12}$, that is,
  by solving a stochastic optimization problem on a tree involving~$2^{12}$ leaves, the final cost being
  given by the function~$\tilde\Value_{12}$.
\end{itemize}
Gathering the calculations performed by this algorithm, we obtain
that solving the global problem by mixing dynamic programming and
stochastic programming is done by solving $(n+1)$ stochastic
optimization programs on scenario trees, each involving~$2^{12}$
leaves. The total number of leaves to explore --- when solving
the problem by this mixing method --- is
$(n+1) 2^{12} \approx 4(n+1) \, 10^{3}$,
which gives an estimation of the algorithm computational effort.

This mixing method is to be compared with a pure scenario tree method,
that is, when the problem is solved by a stochastic optimization
program on a scenario tree over the whole time horizon on~24 hours,
the total number of leaves to explore being
$2^{24} \approx 1.6 \, 10^{7}$.
Even using a fairly fine state discretization grid, for example
a grid containing~$100$ points, the resolution by mixing dynamic
programming and stochastic programming --- when compared to the pure
stochastic programming approach --- leads to a quite significant gain,
namely a factor $1.6 \times 10^{7}/4(100+1)\times 10^{3} \approx 40$
in our case.

We performed numerical experiments with a single computer equipped
with 12~Intel Core i5-10500 CPU and 16 GB of RAM. We used the LP package
of the solver Gurobi 9.51. Apart from the solver, all our code
has been implemented with the Julia language and the JuMP modeler.
As we failed to obtain a solution for the original problem on a tree
for the horizon~$\horizon=24$, we performed numerical tests for
shorter horizons, hence for smaller numbers of time steps.
For every $\horizon \in\na{12,14,16,18, 20,24}$, we considered
that a state reduction existed at time~$\horizon/2$. The results are
gathered in Table~\ref{fig:mixing_numerical}, and show that the
computational time --- that is, the CPU time needed to create
the LP model by JuMP and to solve it by Gurobi --- needed by the pure
scenario tree method is very rapidly increasing with the number
of time steps, whereas the computational time needed by the mixing
method grows very slowly with the number of time steps, at least
for the different horizons under consideration.
Finally, note that the mixing method can be easily parallelized
since the computation of the~$n$ values
$\ba{\tilde\Value_{\horizon/2}(\state^{i})}_{i\in\ic{1,n}}$
of the Bellman function~$\tilde\Value_{\horizon/2}$ can be
performed in parallel, reducing the CPU time by a factor~$(n+1)/2$,
that is, approximately~$50$ in our case.
\begin{table}[htb]
\begin{center}
  \begin{tabular}{|c|c|c|}
    \hline
    Horizon~$\horizon$ & Mixing method & Pure scenario tree method \\
    \hline
    $12$               &  6.5 s        & 0.5 s                     \\
    \hline
    $14$               &  6.5 s        & 5.0 s                     \\
    \hline
    $16$               &  6.6 s        & 71.2 s                    \\
    \hline
    $18$               &  6.8 s        & 1009.5 s                  \\
    \hline
    $20$               &  6.8 s        & 137,296.0 s               \\
    \hline
    $24$               &  7.2 s        & optimization failed       \\
    \hline
  \end{tabular}
  \caption{Computational time (in seconds) needed for solving the
  problem by both the mixing method and the pure scenario tree method,
  for different horizons
  \label{fig:mixing_numerical}}
\end{center}
\end{table}

\section{Two-time-scale optimization problems}
\label{Two_time_scale_optimization_problem}

As a second application of the formalism developed in \S\ref{State_reduction_on_multiple_consecutive_time_blocks_and_dynamic_programming_equations},
we show how to tackle a class of two-time-scale optimization problems.
Indeed, some decisions problems naturally involve two different
time scales, because of the timing of decisions ---
as for example long term investment decision and short term monitoring of
physical devices.

In~\S\ref{Structure_of_a_two_time_scale_optimization_problem} and
\S\ref{sec:reformulation-extended-timeline}
we detail the structure
and we formulate the two-time-scale optimization problems
that we consider.
In~\S\ref{Two_time_scale_decomposition},
we show how to decompose such problems by time blocks.
In~\S\ref{Crude_Oil_Procurement_Problem}, we
illustrate the approach on a crude oil procurement problem.

\subsection{Structure of a two-time-scale optimization problem}
\label{Structure_of_a_two_time_scale_optimization_problem}

We provide the data for a two-time-scale multistage optimization problem.

\subsubsubsection{Two time scales}

The slow time scale is represented by a finite totally ordered
set~$\np{\SlowSpan,\preceq}$ as follows --- where $\nex\tslow$
denotes the \emph{successor} of $\tslow\in\SlowSpan$ and $\pre\tslow$
its \emph{predecessor}, and where we use the notation~$t \prec t'$
for~$t \preceq t'$ and~$t \neq t'$ ---
\begin{subequations}
  \begin{gather}
    \min\SlowSpan = \bottom\tslow \prec \cdots  \prec \pre\tslow \prec \tslow
    \prec \nex\tslow \prec \cdots \prec \topp\tslow = \max\SlowSpan \eqfinv
    \intertext{and the fast time scale by a finite totally ordered
      set~$\np{\FastSpan,\preceq}$:}
    \min\FastSpan = \bottom\tfast \prec \cdots  \prec \pre\tfast \prec \tfast
    \prec \nex\tfast \prec \cdots \prec \topp\tfast = \max\FastSpan \eqfinp
  \end{gather} \label{eq:chain_def}
\end{subequations}
In a sense to be made more rigorous later (once a unified timeline will have been defined),
each slow time interval~$ \ClosedIntervalOpen{\tslow}{\nex\tslow} $
is made up of $\cardinal{\FastSpan}$ (cardinality of~$\FastSpan$) fast time steps,
hence the denomination ``two-time-scale''.
For instance, $\SlowSpan=\{ Mo$, $ Tu$, $ We$, $ Th$, $ Fr$, $ Sa$, $ Su \}$ may represent
days, whereas $ \FastSpan = \ic{1,24} $ may represent hours within a day.
In some problems, we might even take $ \FastSpan = \ic{0,24} $ to handle
the fact that two decisions (one slow and one fast) are taken at midnight,
hence an additional fast time step~$0$.

\subsubsubsection{Unified timeline}

We define the unified timeline of the decision problem in two steps.
First, we equip the product set $\SlowSpan\times \FastSpan$ with the following lexicographic order:
\begin{align}
  (\bottom\tslow,\bottom\tfast)
  &\prec \cdots\prec
    (\pre\tslow,\topp\tfast) \prec
    (\tslow,\bottom\tfast)\prec(\tslow,\nex{\bottom\tfast})\prec
    \cdots
    \label{eq:timeline}
  \\
  &\cdots \prec (\tslow,\pre{\topp\tfast})\prec(\tslow,\topp\tfast)\prec(\nex\tslow,\bottom\tfast)
    \prec \cdots \prec(\topp\tslow,\topp\tfast)
    \eqfinp
    \nonumber
\end{align}
More formally, we denote by~$\nex{(\tslow,\tfast)}$
the successor of $(\tslow,\tfast)$ in
$\SlowSpan\times\FastSpan \setminus\na{(\topp\tslow,\topp\tfast)}$, with
\begin{subequations}
  \label{eq:lexicographic-order}
  \begin{equation}
    \nex{(\tslow,\tfast)} =
    \begin{cases}
      (\tslow,\nex\tfast) & \text{if  } \tfast \not= \topp\tfast \eqfinv
      \\
      (\nex\tslow,\bottom\tfast) & \text{if  } \tfast = \topp\tfast \eqfinp
    \end{cases}
  \end{equation}
  Similarly, we denote by $\pre{(\tslow,\tfast)}$ the predecessor of
  $(\tslow, \tfast)$ in
  $\SlowSpan\times\FastSpan\setminus\na{(\bottom\tslow,\bottom\tfast)}$,
  with
  \begin{equation}
    \pre{(\tslow,\tfast)} =
    \begin{cases}
      (\tslow,\pre\tfast) & \text{if  } \tfast \not= \bottom\tfast \eqfinv
      \\
      (\pre\tslow,\topp\tfast) & \text{if  } \tfast = \bottom\tfast \eqfinp
    \end{cases}
  \end{equation}
\end{subequations}
In the product set $\SlowSpan{\times}\FastSpan$, the first time
$(\bottom{\tslow},\bottom{\tfast})$ does not coincide with a slow time
(the couple $ \np{Mo,0} $ does not correspond to Monday in our running example in~\S\ref{Crude_Oil_Procurement_Problem}).
Thus, we add to the product set $\SlowSpan{\times}\FastSpan$ an extra time denoted
by $(\pre{\bottom{\tslow}},\topp{\tfast})$,
corresponding to the extra slow time~$\pre{\bottom{\tslow}}$,
which is such that
$\pre{(\bottom{\tslow},\bottom{\tfast})}= (\pre{\bottom{\tslow}},\topp{\tfast})$.
We denote by
$\overline{\SlowSpan}$ the set $\na{\pre{\bottom{\tslow}}}\cup \SlowSpan$
and by $\overline{\SlowSpan{\times}\FastSpan}$ the set $(\pre{\bottom{\tslow}},\topp{\tfast})\cup\np{\SlowSpan{\times}\FastSpan}$,
also called the \emph{extended timeline} when equipped with an order~$\preceq$
as follows (where we use the notation~$(\tslow,\tfast) \prec (\tslow',\tfast')$
for~$(\tslow,\tfast) \preceq (\tslow',\tfast')$
and~$(\tslow,\tfast) \neq (\tslow',\tfast')$)
\begin{align}
  (\pre{\bottom\tslow},\topp\tfast)
  &\prec
    (\bottom\tslow,\bottom\tfast)
    \prec \cdots\prec
    (\pre\tslow,\topp\tfast) \prec
    (\tslow,\bottom\tfast)\prec(\tslow,\nex{\bottom\tfast})\prec
    \cdots
    \nonumber
  \\
  &\cdots \prec (\tslow,\pre{\topp\tfast})\prec(\tslow,\topp\tfast)
    \prec(\nex\tslow,\bottom\tfast)\prec \cdots \prec(\topp\tslow,\topp\tfast)
    \eqfinp
    \label{eq:extended-timeline}
\end{align}
The two-time-scale optimization problem will be formulated on the extended timeline
$\overline{\SlowSpan{\times}\FastSpan}$, which we trivially identify with the
time set $\ic{0,\horizon}$, where $ \horizon =\cardinal{\SlowSpan}\times\cardinal{\FastSpan} $.

\begin{figure}[htb]
  \begin{center}
    \includegraphics[width=\textwidth]{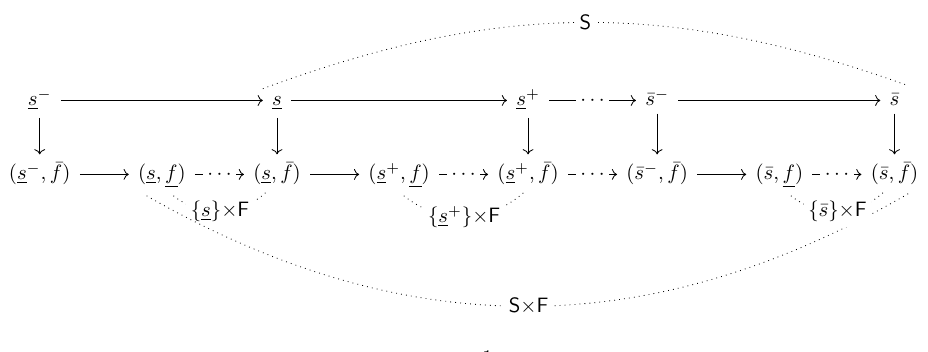}
  \end{center}
  \caption{The product timeline with an extra starting
           point~$(\pre{\bottom{\tslow}},\topp{\tfast})$
           \label{fig:The_product_time_line_with_an_extra_starting_point}}
\end{figure}

\subsubsubsection{Decisions}

We suppose given
\begin{itemize}
\item a family
  $ \nseqa{\CONSO[\tslow]}{\tslow\in\overline{\SlowSpan}\setminus\na{\topp\tslow}} $
  of \emph{slow time scale decision Borel spaces},
  and a family $\nseqa{\PRICE[\tslow]}{\tslow\in \SlowSpan} $
  of \emph{slow time scale uncertainty Borel spaces},
\item a family $
  \nseqa{\BUY[(\tslow,\tfast)]}{(\tslow,\tfast)\in \SlowSpan\times \np{\FastSpan\setminus\na{\topp\tfast}}} $
  of \emph{fast time scale decision Borel spaces},
  and a family\\
  $\nseqa{\COST[(\tslow,\tfast)]}{(\tslow, \tfast)\in\SlowSpan\times \np{\FastSpan\setminus\na{\bottom\tfast}}}$
  of \emph{fast time scale uncertainty Borel spaces}.
\end{itemize}

\subsubsubsection{Dynamics}

We suppose given a family
$\nseqa{\BUFFER[(\tslow,\tfast)]}{(\tslow,\tfast)\in\overline{\SlowSpan\times \FastSpan}}$
of \emph{fast time scale state Borel spaces}. For the sake of simplicity,
we set $\STOCK[\tslow] = \BUFFER[(\tslow,\topp\tfast)]$ for all
$\tslow \in \overline{\SlowSpan}$. Thus, the slow time $\slow \in \overline{\SlowSpan}$ is
identified with the two scale time $(\tslow, \topp{\tfast})$, as illustrated
in Figure~\ref{fig:The_product_time_line_with_an_extra_starting_point}.
We also suppose given a family
$ \nseqa{\hatdyn_{\tslow}^{\slow}}{\tslow\in\overline{\SlowSpan}\backslash\na{\topp\slow}} $
of \emph{slow time scale dynamics Borel-measurable mappings},
that represent the evolution ``driven at the slow time scale'' given,
for $\tslow\in\overline{\SlowSpan}\setminus\na{\topp\slow}$, by\footnote{%
  We stress that the slow time scale dynamics~\eqref{eq:slow-dynamics} yields as output the first
  fast state of the slow period (and not the next slow state).
  Thus, the slow time scale dynamics~\eqref{eq:slow-dynamics} is \emph{not} a
  dynamics from one slow state to the next slow state.}
\begin{subequations}
  \label{eq:fastslow-dynamics}
  \begin{align}
    \hatdyn_{\tslow}^{\slow}
    :
    \STOCK[\tslow]
    {\times} \CONSO[\tslow]
    {\times} \PRICE[\nex\tslow]
    &\to
      \BUFFER[(\nex\tslow,\bottom\tfast)]
      \eqfinv
      \nonumber
    \\
    \bp{\stock[\tslow], \conso[\tslow],\price[\nex\tslow]}
    & \mapsto
      \buffer[(\nex\tslow,\bottom\tfast)]
      = \hatdyn_{\tslow}^{\slow}
      \bp{\stock[\tslow], \conso[\tslow],\price[\nex\tslow]}
      \eqfinp
      \label{eq:slow-dynamics}
      \intertext{%
      We suppose given a family $ \nseqa{\hatdyn^{\fast}_{(\tslow,\tfast)}}
      {(\tslow,\tfast)\in\SlowSpan\times\np{\FastSpan\backslash\na{\topp\tfast}}} $
      of \emph{fast time scale dynamics Borel-measurable mappings},
      that represent the evolution ``driven at the fast time scale''
      given, for all $ \tslow \in \SlowSpan$ and $\tfast \in \FastSpan\setminus\na{\topp\tfast}$, by}
    % \\
    \hatdyn_{(\tslow,\tfast)}^{\fast}
    :%
    \BUFFER[(\tslow,\tfast)]{\times}
    \BUY[(\tslow,\tfast)]{\times}
    \COST[\nex{(\tslow,\tfast)}]
    & \to
      \BUFFER[\nex{(\tslow,\tfast)}]
      \eqfinv
      \nonumber
    \\
    \bp{ \buffer[(\tslow,\tfast)],\buy[(\tslow,\tfast)],
    \cost[\nex{(\tslow,\tfast)}]
    }
    &\mapsto
      \buffer[\nex{(\tslow,\tfast)}] = \hatdyn_{(\tslow,\tfast)}^{\fast}
      \bp{ \buffer[(\tslow,\tfast)],\buy[(\tslow,\tfast)],
      \cost[\nex{(\tslow,\tfast)}]
      }
      \eqfinp
      \label{eq:fast-dynamics}
  \end{align}
\end{subequations}

\subsubsubsection{Cost functions}

\begin{subequations}
  We suppose given a family
  $\nseqa{\instantcost[\tslow]}{\tslow\in\overline{\SlowSpan}\setminus\na{\topp\slow}} $
  of \emph{slow time scale nonnegative lower semianalytic cost functions},
  with
  \begin{align*}
    & \instantcost[\pre\tslow] :
      \STOCK[\pre\tslow]{\times}\CONSO[\pre\tslow]{\times}\PRICE[\tslow]{\times}%
      \underbrace{  \prod_{\tfast\in \FastSpan\backslash\na{\topp\tfast}}
      \bp{
      \BUFFER[(\tslow,\tfast)]{\times}\BUY[(\tslow,\tfast)]{\times}%
      \COST[\nex{(\tslow,\tfast)}]
      } }_{ \textrm{interval } \OpenIntervalOpen{\pre\tslow}{\tslow}
      {=\na{(\tslow,\bottom\tfast),\ldots,(\tslow,\pre{\topp\tfast})}} }
            \to \bRR \eqfinv
      \intertext{for $\tslow\in\SlowSpan$,
      and a \emph{slow time scale nonnegative lower semianalytic final cost function}~$\instantcost[\topp\tslow]$}
    & \instantcost[\topp\tslow] :
      \STOCK[\topp\tslow]
      \to \bRR \eqfinv
  \end{align*}
\end{subequations}
that make up, by summation, an intertemporal cost
\begin{equation}
  \dsum[\tslow\in\SlowSpan] \instantcost[\pre\tslow]
  \Bp{
    \stock[\pre\tslow],\conso[\pre\tslow],\price[\tslow],
    \nseqp{ \buffer[(\tslow,\tfast)],\buy[(\tslow,\tfast)], \cost[\nex{(\tslow,\tfast)}] }
    {\tfast\in \FastSpan\backslash\na{\topp\tfast}}
  }
  + \instantcost[\topp\tslow]\bp{\stock[\topp\tslow]} \eqfinp
  \label{eq:intertemporal-criterion}
\end{equation}

\subsubsubsection{Stochastic kernels}

Finally, we suppose given a family of \emph{constant
  slow time scale Borel-measurable stochastic kernels}
$\nseqa{\rho^{\slow}_{\tslow:\nex\tslow}}{\tslow\in
  \overline{\SlowSpan}\setminus\na{\topp\tslow}}$
\begin{subequations}
  \label{eq:natural_kernel}
  \begin{align}
    \rho^{\slow}_{\tslow:\nex\tslow}
    & \in
      \Delta(\PRICE[\nex\tslow])
      \eqsepv
      \forall \tslow \in \overline{\SlowSpan}\setminus\na{\topp\tslow}
      \eqfinv
      \label{eq:natural_kernel_slow}
      \intertext{and, for each $\tslow \in \SlowSpan$, a family $ \nseqa{\rho^{\fast}_{(\tslow,\tfast):\nex{(\tslow,\tfast)}}}
      {\tfast\in\FastSpan\backslash\na{\topp\tfast}} $ of
      \emph{fast time scale Borel-measurable stochastic kernels}}
    % \\
    \rho^{\fast}_{(\tslow,\tfast):\nex{(\tslow,\tfast)}}
    &:
      \PRICE[\tslow]\times
      \underbrace{  \dprod[\tfast' = \nex{\bottom\tfast}][\tfast]
      \COST[(\tslow,\tfast')]
      }_{ \textrm{in interval } \ClosedIntervalOpen{\pre\tslow}{\tslow} }
      \longrightarrow
      \Delta(\COST[\nex{(\tslow,\tfast)}])
      \eqsepv
      \forall \tslow \in \SlowSpan
      \eqsepv
      \forall \tfast \in \FastSpan \backslash\na{\topp\tfast}
      \eqfinv
      \label{eq:natural_kernel_fast}
  \end{align}
  with the convention that the Cartesian products of spaces in
  Equations~\eqref{eq:natural_kernel_slow} and~\eqref{eq:natural_kernel_fast}
  reduce to nothing when the upper index of the Cartesian product is
  strictly lower that the corresponding lower index.  Note that, for a given
  $\tslow \in \SlowSpan$, each fast time scale stochastic kernel
  $\rho^{\fast}_{(\tslow,\tfast):\nex{(\tslow,\tfast)}}$, only depends on the
  noises of the slow time
  interval~$ \ClosedIntervalOpen{\pre\tslow}{\tslow}
  =\{(\pre\tslow,\topp\tfast),(\tslow,\bottom\tfast),\ldots,(\tslow,\pre{\topp\tfast})\}$.
  The (constant) assumption~\eqref{eq:natural_kernel_slow} and the (single
  block) assumption~\eqref{eq:natural_kernel_fast} correspond to stochastic
  independence between time blocks, and will be useful in the proof of
  Proposition~\ref{gdp:prop:two_time_scale_decomposition}.
\end{subequations}

\subsection{Formulation of a two-time-scale optimization problem on the product timeline}
\label{sec:reformulation-extended-timeline}

To apply Theorem~\ref{gdp:thm:DPB_family}, we introduce sets associated
with the extended timeline~\eqref{eq:extended-timeline} by
\begin{subequations}
  \begin{align}
    \STATE_{(\tslow,\tfast)}
    &=
      \begin{cases}
        \STOCK[\tslow] & \text{if}~\tfast = \topp\tfast\\
        \BUFFER[(\tslow,\tfast)]
        & \text{if}~\tfast \not= \topp\tfast\\
      \end{cases}
    \eqsepv \forall (\tslow,\tfast) \in \overline{\SlowSpan{\times}\FastSpan}
    \eqfinv
    \\
    \CONTROL_{(\tslow,\tfast)}
    &=
      \begin{cases}
        \CONSO[\tslow] & \text{if}~\tfast = \topp\tfast\\
        \BUY[(\tslow,\tfast)]
        & \text{if}~\tfast \not= \topp\tfast\\
      \end{cases}
    \eqsepv \forall  (\tslow,\tfast) \in \overline{\SlowSpan{\times}\FastSpan}\setminus\na{(\topp\tslow,\topp\tfast)}
    \eqfinv
    \\
    \COST[(\tslow,\tfast)][]
    &=
      \begin{cases}
        \PRICE[\tslow]  & \text{if}~\tfast = \bottom\tfast\\
        \COST[(\tslow,\tfast)] & \text{if}~\tfast \not= \bottom\tfast\\
      \end{cases}
    \eqsepv \forall  (\tslow,\tfast) \in \SlowSpan{\times}\FastSpan
    \eqfinv
    \\
    \intertext{with the particular case of the extra initial slow time}
    \COST[(\pre{\bottom\tslow},\topp\tfast)][]
    &= \STOCK[\pre{\bottom\tslow}]
     \label{eq:tts-initial-condition} \eqfinv
  \end{align}
\end{subequations}
and a family of state dynamics
$
\hatdyn_{(\tslow,\tfast)}:
\STATE_{(\tslow,\tfast)}{\times}
\CONTROL_{(\tslow,\tfast)}{\times}
\UNCERTAIN_{\nex{(\tslow,\tfast)}}$
$\to \STATE_{\nex{(\tslow,\tfast)}}$
defined by
\begin{equation}
  \hatdyn_{(\tslow,\tfast)} =
  \begin{cases}
    \hatdyn_{\tslow}^{\slow} &
    \text{if}~\tfast = \topp\tfast\\
    \hatdyn_{(\tslow,\tfast)}^{\fast} &
    \text{if}~\tfast \not= \topp\tfast
  \end{cases}
  \eqsepv \forall
  (\tslow,\tfast) \in
  \overline{\SlowSpan{\times}\FastSpan}\setminus\na{(\topp\tslow,\topp\tfast)}
  \label{eq:extended-dynamics}
  \eqfinp
\end{equation}

\begin{figure}[htb]
  \begin{center}
    \includegraphics[width=\textwidth]{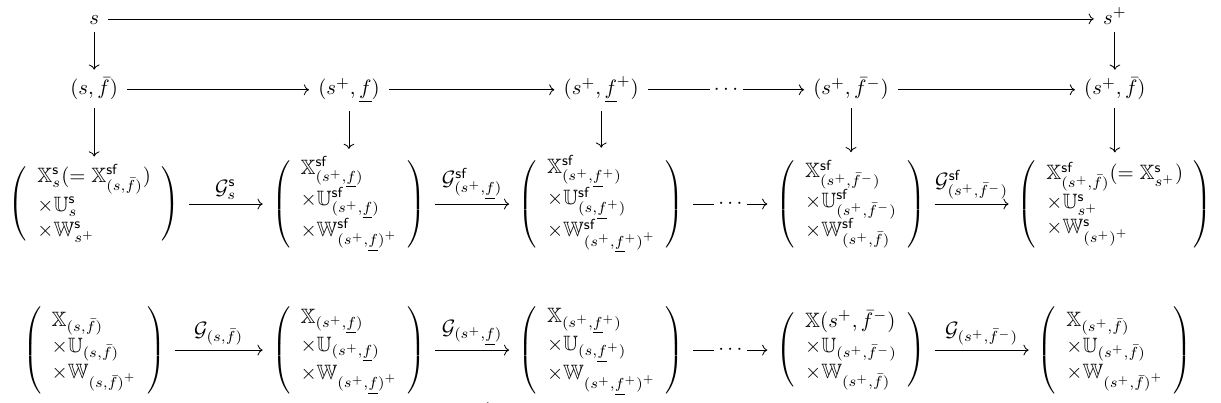}
  \end{center}
  \caption{Original dynamics and their reformulation on the product timeline
    on the slow time  interval $\ClosedIntervalOpen{\tslow}{\nex\tslow}
     {=\na{(\tslow,\topp\tfast),(\nex\tslow,\bottom\tfast),\ldots,(\nex\tslow,\pre{\topp\tfast})}}$}
\end{figure}

From these spaces, we deduce the history spaces and the histories
for all $(\tslow,\tfast) \in \overline{\SlowSpan{\times}\FastSpan}$
\begin{subequations}
  \begin{align}
  \HISTORY_{(\tslow,\tfast)}
    & =
      \UNCERTAIN_{(\pre{\bottom\tslow},\topp\tfast)} \times
      \prod_{(\bottom\tslow,\bottom\tfast) \preceq (\tslow',\tfast') \preceq (\tslow,\tfast)}
      \Bp{\CONTROL_{\pre{(\tslow',\tfast')}}\times\UNCERTAIN_{(\tslow',\tfast')}}
      \eqfinv \\
  \hist[(\tslow,\tfast)]
    & =
      \Bp{\alea[{(\pre{\bottom\tslow},\topp\tfast)}] ,
      \bseqp{
      { \control_{\pre{(\tslow',\tfast')}} ,  \alea[{(\tslow',\tfast')}]}}
      {(\bottom\tslow,\bottom\tfast) \preceq (\tslow',\tfast') \preceq (\tslow,\tfast)}}
      \eqfinv
  \end{align}
\end{subequations}
and, for suitable indices, the partial history sets and the partial histories
\begin{subequations}
  \begin{align}
    \HISTORY_{(\tslow,\tfast):(\tslow',\tfast')}
    &=
      \prod_{(\tslow,\tfast) \preceq (\tslow'',\tfast'') \preceq (\tslow',\tfast')}
      \np{ \CONTROL_{\pre{(\tslow'',\tfast'')}} \times \UNCERTAIN_{(\tslow'',\tfast'')} }
      % \eqsepv \forall (\tslow,\tfast) \in \SlowSpan\times\FastSpan
      \eqfinv
    \\
    \hist[(\tslow,\tfast):(\tslow',\tfast')]
    & =
      \bp{
      \nseqp{
      { \control_{\pre{(\tslow'',\tfast'')}} ,  \alea[{(\tslow'',\tfast'')}]}}
      {(\tslow,\tfast) \preceq (\tslow'',\tfast'') \preceq (\tslow',\tfast')}}
    % \eqsepv \forall (\tslow,\tfast) \in \SlowSpan\times\FastSpan
      \eqfinp
  \end{align}
\end{subequations}

The cost criterion formulated in Equation~\eqref{eq:intertemporal-criterion}
combined with state dynamics leads to a (nonnegative lower semianalytic)
cost criterion $\criterion:\HIST[{(\topp\tslow,\topp\tfast)}] \to \bRR$.

Based on the stochastic kernels~\eqref{eq:natural_kernel_slow}
and~\eqref{eq:natural_kernel_fast}, we introduce stochastic kernels
$ \rho_{(\tslow,\tfast):\nex{(\tslow,\tfast)}}$
associated with the extended timeline~\eqref{eq:extended-timeline}, for each
$(\tslow,\tfast)\in
\overline{\SlowSpan{\times}\FastSpan}\setminus\na{\topp\tslow, \topp\tfast}$, by
$\rho_{(\tslow,\tfast):\nex{(\tslow,\tfast)}}:
\HIST[{(\tslow,\tfast)}]
\longrightarrow \Delta(\ALEA[\nex{(\tslow,\tfast)}])$ with
\begin{align*}
  % &
      \rho_{(\tslow,\tfast):\nex{(\tslow,\tfast)}}
    \bsetp{\dd \alea[\nex{(\tslow,\tfast)}]}{\hist[(\tslow,\tfast)]}
    % \nonumber \\
    % &
      \;=
    \begin{cases}
      \rho^{\slow}_{\tslow:\nex\tslow}
      (\dd \price[\nex\tslow])
      & \text{if }\tfast = \topp\tfast
      \eqfinv
      \\
      \rho^{\fast}_{(\tslow,\tfast):\nex{(\tslow,\tfast)}}
      \bsetp{\dd \cost[(\tslow,\nex{\tfast})]}
      {\price[\tslow],\cost[(\tslow,\nex{\bottom\tfast})],\cdots,\cost[(\tslow,\tfast)]}
      & \text{if }\tfast \not= \topp\tfast
      \eqfinp
    \end{cases}
    % \label{two_scale_kernel}
\end{align*}
Note that, for $\tfast \not= \topp\tfast$, the stochastic kernels
$\rho_{(\tslow,\tfast):\nex{(\tslow,\tfast)}}:
\HIST[{(\tslow,\bottom\tfast):(\tslow,\tfast)}]
\to \Delta(\ALEA[\nex{(\tslow,\tfast)}])$, only depend on the
partial history uncertainty part from $(\tslow,\bottom\tfast)$ to
$(\tslow,\tfast)$, and not on the (past) controls.

The components of the problem are now formulated on the extended timeline
$\overline{\SlowSpan{\times}\FastSpan}$, already identified with the
time set $\ic{0,\horizon}$. Thus, we are in the framework of
\S\ref{The_Bertsekas-Shreve_setting}
and we aim at solving
an optimization problem as formulated in Equation~\eqref{eq:Vzero}.

\subsection{Two-time-scale decomposition}
\label{Two_time_scale_decomposition}

The existence of Bellman equations for a two-time-scale optimization
problem is given by the following proposition.

\begin{proposition}
  Consider a two-time-scale optimization problem as formulated in
  \S\ref{Structure_of_a_two_time_scale_optimization_problem} and \S\ref{sec:reformulation-extended-timeline}.
  The optimization problem~\eqref{eq:Vzero}
  has a solution given by a dynamic programming equation at
  the slow scale. More precisely, let $\nseqa{\tilValue[\tslow]}{\tslow \in \overline{\SlowSpan}}$ be given by
  $\tilValue[\topp\tslow] = \instantcost[\topp\tslow]$
  and, for $\tslow \in \overline{\SlowSpan}\setminus\na{\topp\tslow}$,
  by the backward induction\footnote{Here again, the formula~\eqref{gdp:eq:overlineBelOp-expression-rephrased}
  represents a nested sequence of infima of integrals
  (with respect to different stochastic kernels).}
  \begin{align}
    & \tilValue[\tslow]\np{\stock[\tslow]}
      =
      \inf_{\conso[\tslow]\in\CONSO[\tslow]}
      \int_{\PRICE[\nex\tslow]}
      \rho^{\slow}_{\tslow:\nex\tslow} \np{\dd \price[\nex\tslow]}
      \nonumber \\
    &
      \quad
      \inf_{\buy[(\nex\tslow,\bottom\tfast)]\in\BUY[(\nex\tslow,\bottom\tfast)]}
      \int_{\COST[(\nex\tslow,\nex{\bottom\tfast})]}
      \rho^{\fast}_{(\nex\tslow,{\bottom\tfast}):(\nex\tslow,\nex{\bottom\tfast})}
      \kernelargs{\price[\nex\tslow]}{\dd\cost[(\nex\tslow,\nex{\bottom\tfast})] }
      \; \cdots
      \nonumber    \\
    & \quad\quad
      \inf_{\buy[(\nex\tslow,\pre{\topp\tfast})]\in\BUY[(\nex\tslow,\pre{\topp\tfast})]}
      \int_{\COST[(\nex\tslow,\topp\tfast)]}
          \rho^{\fast}_{(\nex\tslow,\pre{\topp\tfast}):(\nex\tslow,{\topp\tfast})}
          \kernelargs{\price[\nex\tslow],\cost[(\nex\tslow,\nex{\bottom\tfast})],
                      \cdots,\cost[(\nex\tslow,\pre{\topp\tfast})]}
          {\dd\cost[(\nex\tslow,{\topp\tfast})]}
          \nonumber \\
    & \quad\quad\quad
      \Big(
      \instantcost[\tslow]
      \np{\stock[\tslow], \conso[\tslow],\price[\nex\tslow],\ldots,
      \buy[(\nex\tslow,\pre{\topp\tfast})],
      \cost[(\nex\tslow,{\topp\tfast})]}
      \nonumber  \\
    & \quad\quad\quad\;
      +
      \tilValue[\nex\tslow]\bp{\hatdyn_{\tslow:\nex\tslow}
      \np{\stock[\tslow], \conso[\tslow],\price[\nex\tslow],\ldots,
      \buy[(\nex\tslow,\pre{\topp\tfast})],
      \cost[(\nex\tslow,{\topp\tfast})]}}
      \Big)
      \label{gdp:eq:overlineBelOp-expression-rephrased}
      \eqfinv
  \end{align}
  where $\hatdyn_{\tslow:\nex\tslow} $ is the composition
  $\hatdyn_{\tslow:\nex\tslow} =
  \hatdyn_{(\nex{\tslow},\pre{\topp\tfast})}^{\fast}
  \circ \cdots \circ$
  $\hatdyn_{(\nex{\tslow},\bottom\tfast)}^{\fast}$
  $\circ \; \hatdyn_{\tslow}^{\slow}$
  associated with the state dynamics
  defined in~\eqref{eq:fastslow-dynamics}.
  Then, the value of the optimization problem~\eqref{eq:Vzero}
  is given by
  $\tilValue[\pre{\bottom\tslow}]\np{\stock[\pre{\bottom\tslow}]}$,
  where the initial condition $\stock[\pre{\bottom\tslow}]$
  corresponds to~$\uncertain_{0}$ in~\eqref{eq:Vzero},
  as stated by~\eqref{eq:tts-initial-condition}.
  \label{gdp:prop:two_time_scale_decomposition}
\end{proposition}

\begin{proof}
  The proof is an application of
  Theorem~\ref{cor:instantaneous-costs-in-addition-to-final-cost}
  with the help of
  Remarks~\ref{rem:composed-state-dynamics-as-reduction-mapping}
  and~\ref{rem:block-independent-noises-and-kernel}.
  First, we have re-framed in \S\ref{sec:reformulation-extended-timeline}
  the two-time-scale optimization problems
  described in~\S\ref{Structure_of_a_two_time_scale_optimization_problem}
  in the formalism of~\S\ref{The_Bertsekas-Shreve_setting}
  with the help of the extended timeline~\eqref{eq:extended-timeline}.
  Second, as we are given state dynamics~\eqref{eq:extended-dynamics}
  on the extended timeline and thanks to
  Remark~\ref{rem:composed-state-dynamics-as-reduction-mapping},
  we obtain a state reduction
  at times $\nseqa{(\tslow,\topp\tfast)}{\tslow \in \SlowSpan}$
  by composition of the state dynamics. Moreover,
  as the slow time scale stochastic kernels given by Equation~\eqref{eq:natural_kernel_slow} are constant,
  the state reduction across the slow time scale is compatible
  with the stochastic kernels
  (see Remark~\ref{rem:block-independent-noises-and-kernel}).
  We are thus able to apply
  Theorem~\ref{cor:instantaneous-costs-in-addition-to-final-cost}
  and obtain the slow time scale Bellman
  recursion~\eqref{gdp:eq:overlineBelOp-expression-rephrased}
  as a special case of Equation~\eqref{gdp:eq:overlineBelOp-expression}.
\end{proof}

The slow time scale Bellman
equation~\eqref{gdp:eq:overlineBelOp-expression-rephrased}
is as difficult to solve as the Bellman equation on the extended timeline.
However, the interest of~\eqref{gdp:eq:overlineBelOp-expression-rephrased}
lies elsewhere.
Imagine that one is able to obtain, in a relatively easy way, lower
$\underline{\tilValue[\tslow]}$ and upper $\overline{\tilValue[\tslow]}$ approximations of $\tilValue[\tslow]$
in~\eqref{gdp:eq:overlineBelOp-expression-rephrased}.
Then, in order to obtain optimal controls for the
optimization problem~\eqref{eq:Vzero}, one can proceed as follows.
By replacing the last term~$ \tilValue[\nex\tslow] $
of~\eqref{gdp:eq:overlineBelOp-expression-rephrased}
by either $\underline{\tilValue[\nex\tslow]}$
or $\overline{\tilValue[\nex\tslow]}$,
one can now solve a (lower or upper) surrogate of
Equation~\eqref{gdp:eq:overlineBelOp-expression-rephrased}
and thus obtain the optimal controls
on the time interval~{\( \ClosedIntervalOpen{\tslow}{\nex\tslow}
  =\na{(\tslow,\topp\tfast),(\nex\tslow,\bottom\tfast),\ldots,(\nex\tslow,\pre{\topp\tfast})} \).}
For instance, one could use scenario decomposition methods,
like progressive hedging \cite{Rockafellar-Wets:1991}, that do not require statistical
independence of noises within the slow time interval~$
\ClosedIntervalOpen{\tslow}{\nex\tslow} $.
Thus, the two-time-scale stochastic optimization problem as formulated in
\S\ref{Structure_of_a_two_time_scale_optimization_problem}
and \S\ref{sec:reformulation-extended-timeline}
can be approximatively solved, from below and from above, by a mix of
slow time scale dynamic programming and of (for example)
progressive hedging (or any other method, including dynamic programming).

\subsection{Illustration with the crude oil procurement problem}
\label{Crude_Oil_Procurement_Problem}

\begin{figure}[h]
  \begin{center}
    \includegraphics[width = 0.9\textwidth]{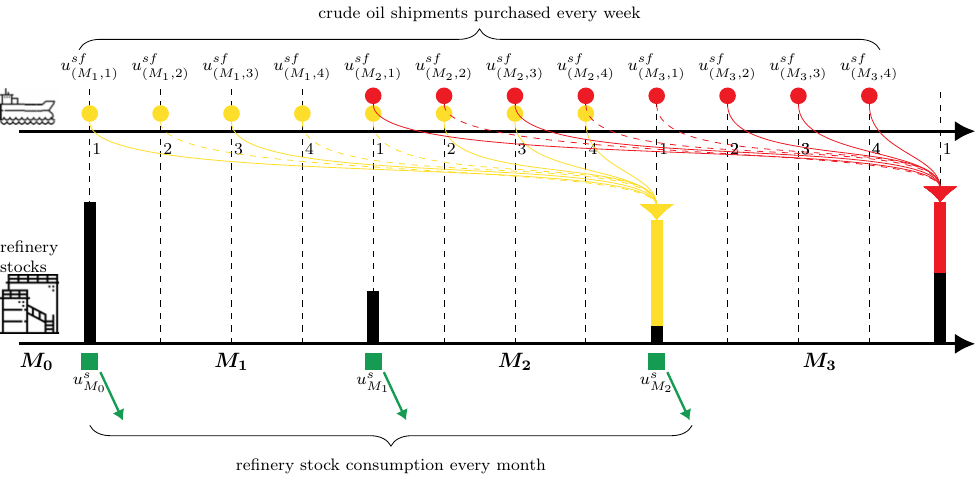}
  \end{center}
  \caption{Procurement of crude oil over 3 months $M_1$, $M_2$ and $M_3$,
    where a circle~\Large$\circ$ \normalsize denotes purchase decisions and
    a square~$\square$ denotes consumption decisions\label{fig:procurement_scheme}}
\end{figure}

This illustration stems from a research work done in partnership
with TotalEnergies, in the context of a PhD thesis \cite{Martin:2021}.
Crude oil procurement is the part of the oil supply chain that sits
between the production of crude oil and its processing in a refinery.
The goal of procurement is to purchase crude oil from various suppliers
around the world and having it delivered in time to the refinery to be
processed. As illustrated in Figure~\ref{fig:procurement_scheme},
every month (on the bottom line) a refinery receives crudes that have
been bought during the 8 previous weeks (on the upper line).

The problem naturally displays two time scales.
On the one hand, deliveries to the refinery are made at the beginning
of each month, and crude consumption is set once a month.
On the other hand, crude oil shipments can be purchased at the frequency
of the week; every week, a selection of shipments is presented to
the decision-maker who must decide which shipments to purchase.
Following the construction of the extended timeline in~\eqref{eq:extended-timeline},
we represent by the sequence
\begin{align}
  % &(M_0,5)
  % \label{eq:total_timeline}\\
  (M_0,5)  &\; \prec (M_1,1)  \prec (M_1,2) \prec (M_1,3) \prec (M_1,4) \prec (M_1,5)
             \nonumber\\
    &\;\; \prec (M_2,1) \prec (M_2,2) \prec (M_2,3) \prec (M_2,4) \prec (M_2,5)
      \label{eq:total_timeline} \\ %\nonumber\\
    &\;\;\; \prec(M_3,1) \prec (M_3,2) \prec (M_3,3) \prec (M_3,4) \prec (M_3,5)
      \nonumber
\end{align}
the timeline associated with Figure~\ref{fig:procurement_scheme}
(notice that we consider that a month is made of 4~weeks).
The initial time stage~$(M_0,5)$ corresponds to the additional time~$
(\pre{\bottom\tslow},\topp\tfast) $ in~\eqref{eq:extended-timeline}.
The times $(M_1,5)$ and $(M_2,5)$ both represent the ``end of the month''
when a consumption decision (slow scale decision $\conso[s]$ on the bottom line
of Figure~\ref{fig:procurement_scheme}) is taken.

We now illustrate how the crude oil procurement problem can be put
in the form of a two-time-scale optimization problem such
as presented in \S\ref{Structure_of_a_two_time_scale_optimization_problem}.
For this purpose, we proceed to the identifications in
Table~\ref{table:total_identification}.

\begin{table}
  \begin{tabular}{|c|c|}
    \hline
    Notations from \S\ref{Structure_of_a_two_time_scale_optimization_problem}
      & Crude oil procurement \\
      & \\
    \hline\hline
    $\SlowSpan$
      & set of months during which we manage the refinery; \\
      & in Figure \ref{fig:procurement_scheme}, $\SlowSpan=\{M_1, M_2, M_3\}$ \\
    \hline
    $\FastSpan$
      & set of weeks in each month; \\
      & in Figure \ref{fig:procurement_scheme}, $\FastSpan=\{1,2,3,4,5\}$ \\
    \hline
    $\CONSO[\tslow]$
      & set of crude oil consumptions during the month $\nex\tslow$ \\
    \hline
    $\PRICE[\nex\tslow]$
      & set of product prices for the month $\nex\tslow$ \\
    \hline
    $\BUY[(\tslow,\tfast)]$
      & set of crude shipments purchased in week $(\tslow,\tfast)$ \\
    \hline
    $\COST[\nex{(\tslow,\tfast)}]$
      & set of crude oil prices in week $(\tslow,\tfast)$ \\
    \hline
    $\hatdyn_{(\tslow,\tfast)}^{\fast}$
      & accumulation of shipments purchased in $(\tslow, \tfast)$ \\
    \hline
    $\hatdyn_{\tslow}^{\slow}$
      & delivery of orders and consumption of crude oil for the month $\nex\tslow$ \\
    \hline
    $\instantcost[\tslow]$
      & operational costs during the month $\tslow$ \\
      & (crude oil purchases during $\tslow$ - earnings from production) \\
    \hline
    $\instantcost[M_4]$
      & end cost associated with the state $\buffer[M_3][s] = \buffer[(M_3,5)][sf]$ \\
      & valuation of the buffers and stocks in the refinery \\
      & before the beginning of the month $M_4$ \\
    \hline
  \end{tabular}
  \caption{Identification of the elements introduced
    in~\S\ref{Structure_of_a_two_time_scale_optimization_problem}
    with elements of the crude oil procurement problem}
  \label{table:total_identification}
\end{table}

We call $\tslow-$buffer (resp. $\pre\tslow-$buffer), the temporary stock
that is created at the beginning of the month $\tslow$ (resp. $\pre\tslow$)
and that will be delivered two months after. For instance, in
Figure~\ref{fig:procurement_scheme}, the yellow disks represent
the $M_1-$buffer and the red disks represent the $M_2-$buffer.
We introduce the state variable
$\buffer[(\tslow,\tfast)] = \Bp{\text{$\pre\tslow-$buffer},
 \text{$\tslow-$buffer},\text{refinery stocks}}$,
together with the accumulation dynamics~$\hatdyn_{(\tslow,\tfast)}^{\fast}$
for the buffers, and the accumulation dynamics~$\hatdyn_{\tslow}^{\slow}$
for the stocks. Regarding the criterion to minimize, it is an intertemporal
cost like in~\eqref{eq:intertemporal-criterion} with slow time scale
cost functions~\( \instantcost[\pre\tslow] \) made of minus the purchases
of crude oil plus the selling of finished products inside a week.
Supposing that the products prices are independent month by month,
we represent this assumption by a family of constant stochastic kernels
$\nseqa{\rho^{\slow}_{\tslow:\nex\tslow}}
{\tslow\in\overline{\SlowSpan}\setminus\na{\topp\tslow}}$.
By contrast, we do not assume that the crude prices are independent
week by week, and the possible dependency is modeled by stochastic
kernels $\rho^{\fast}_{(\tslow,\tfast):\nex{(\tslow,\tfast)}}$
for $\tfast\in\FastSpan\backslash\na{\topp\tfast}$.

Now that all the elements
from~\S\ref{Structure_of_a_two_time_scale_optimization_problem} have been
identified, Proposition \ref{gdp:prop:two_time_scale_decomposition} enables us
to write a dynamic programming equation such
as~\eqref{gdp:eq:overlineBelOp-expression-rephrased} at the scale of the month,
without losing the time-dependency of crude prices inside the month.

\section{Decision-hazard-decision optimization problems}
\label{gdp:Decision_Hazard_Decision_Dynamic_Programming}

In multistage stochastic optimization, the \emph{decision-hazard-decision}
(DHD) framework corresponds to the case where, at the beginning of each time
interval, a decision is taken without knowing the uncertainty that will
materialize at the end of the time interval (decision-hazard), and, at
the end of the time interval, a \emph{recourse} decision is possible
knowing this uncertainty (hazard-decision). The reader is referred
to~\cite[\S1.2.1]{Carpentier-Chancelier-Cohen-DeLara:2015} for the notions
of decision-hazard and hazard-decision in stochastic optimal control.
A discussion about these notions and an application in stochastic thermal
scheduling can also be found in~\cite{Street:2020}.

In~\S\ref{gdp:Decision_Hazard_Decision_Motivation},
we provide motivation for the decision-hazard-decision framework.
In~\S\ref{Decision-Hazard-Decision_Framework_and_dynamic_programming_equation},
we formalize the decision-hazard-decision framework and we provide a
dynamic programming equation.

\subsection{Motivation for the decision-hazard-decision framework}
\label{gdp:Decision_Hazard_Decision_Motivation}

We illustrate our motivation to develop a formalism for the decision-hazard-decision framework
with a single dam management problem.
We suppose given a stochastic process $ \nseqa{\va{\inflow}_{t}}{t\in \ic{1,\horizon}}$
on a probability space.
We can model the dynamics of the water volume in a dam by
\begin{equation}
  \va{\volume}_{t{+}1} = \min
  \ba{ \volume\upper, \va{\volume}_{t} - \va{\turbined}_{t} + \va{\inflow}_{t{+}1} }
  \eqsepv \forall t \in \ic{0,\horizon{-}1}
  \eqfinv
  \label{gdp:eq:dam_model_decision-hazard}
\end{equation}
where $\volume\upper$ is the maximal dam volume,
$\va{\volume}_{t}$ is the volume (stock) of water
at the beginning of period~$[t,t+1[$,
$\va{\inflow}_{t{+}1}$ is the inflow water volume (rain, etc.)  during~$[t,t+1[$,
$\va{\turbined}_{t}$ is the turbined outflow volume during~$[t,t+1[$.
The control variable~$\va{\turbined}_{t}$ is decided at the \emph{beginning} of period~$[t,t+1[$,
chosen such that $ 0\leq \va{\turbined}_{t} \leq \va{\volume}_{t} $,
supposed to depend on the stock~$\va{\volume}_{t}$
but not on the inflow water~$\va{\inflow}_{t{+}1}$ (as $\va{\inflow}_{t{+}1}$ takes place during~$[t,t+1[$,
hence materializes at~$t+1$, hence the time index~$t+1$).
The~$\min$ operation in Equation~\eqref{gdp:eq:dam_model_decision-hazard}
ensures that the dam volume always remains below its maximal capacity~$\volume\upper$,
but induces a nonlinearity in the dynamics.
This nonlinear dynamics is an obstacle to apply stochastic dual dynamic programming
(SDDP).

Alternatively, we can model the dynamics of the water volume in a dam by
\begin{equation}
  \va{\volume}_{t{+}1}= \va{\volume}_{t}-\va{\turbined}_{t}+\va{\inflow}_{t{+}1}-\va{\spilled}_{t{+}1}
  \eqsepv \forall t \in \ic{0,\horizon{-}1}
  \eqfinv
  \label{gdp:eq:dam_model_decision-hazard-decision}
\end{equation}
where the new control variable~$ \va{\spilled}_{t{+}1} $ is the spilled volume,
decided at the \emph{end} of period~$[t,t+1[$,
supposed to depend both on the stock~$\va{\volume}_{t}$ \emph{and}
on the inflow water~$\va{\inflow}_{t{+}1}$,
and chosen such that
$0\leq \va{\volume}_{t} - \va{\turbined}_{t} + \va{\inflow}_{t{+}1}
- \va{\spilled}_{t{+}1} \leq \volume\upper$.
Thus, with the formulation~\eqref{gdp:eq:dam_model_decision-hazard-decision},
we ``pay the price'' to add one control~$ \va{\spilled}_{t{+}1} $,
but we obtain a linear model instead of the
nonlinear model~\eqref{gdp:eq:dam_model_decision-hazard}.
This is especially interesting when using the stochastic dual dynamic
programming (SDDP) algorithm, for which the linearity of the dynamics
is used to obtain the convexity properties required by the algorithm.

\subsection{Decision-hazard-decision framework and dynamic programming equation}
\label{Decision-Hazard-Decision_Framework_and_dynamic_programming_equation}

Let $\nseqa{\CONTROL^\Before_t}{t \in \ic{0,\horizon{-}1}}$ (head or ``before'' controls),
$\nseqa{\CONTROL^\After_t}{t \in \ic{1,\horizon}}$ (tail or ``after'' controls),
$\nseqa{\UNCERTAIN_t}{t \in \ic{1,\horizon}}$ (uncertainties)
and $\nseqa{\STATE_t}{t \in \ic{0,\horizon}}$ (states)
be sequences of Borel spaces.
\begin{subequations}
  Let also be given Borel-measurable dynamics mappings
  \begin{equation}
    {\dynamics_{t} : \STATE_t \times \CONTROL^\Before_t \times\UNCERTAIN_{t{+}1}\times \CONTROL^\After_{t{+}1} }
    \to \STATE_{t{+}1}
    \eqsepv \forall t \in \ic{0,\horizon{-}1}
    \eqfinv
  \end{equation}
  nonnegative lower semianalytic instantaneous cost functions
  \begin{equation}
    {\coutint_{t} : \STATE_t \times \CONTROL^\Before_t \times \UNCERTAIN_{t{+}1} \times \CONTROL^\After_{t{+}1} }
    \to \bRR
    \eqsepv \forall t \in \ic{0,\horizon{-}1}
    \eqfinv
  \end{equation}
  and a nonnegative lower semianalytic final cost function
  \begin{equation}
    \coutfin : \STATE_{\horizon} \to \bRR
    \eqfinp
  \end{equation}
  \label{eq:DHD_dynamics_cost}
\end{subequations}
Finally, we suppose given a Borel probability space~$\epro$,
a random variable \( \va\State_{0} : \Omega \to \STATE_{0} \), and a
stochastic process $ \nseqa{\va{\Uncertain}_{t}}{t\in \ic{1,\horizon}}$ (noise process),
where \( \va{\Uncertain}_{t} : \Omega \to \UNCERTAIN_{t} \)
for $t\in \ic{1,\horizon}$.

Thus equipped, we consider the following multistage stochastic optimization problem
\begin{subequations}
  \begin{align}
    \inf_{ \nseqa{\np{\va{\Control}^\Before_{t}, \va{\Control}^\After_{t{+}1}}}{t \in \ic{0,\horizon{-}1}} }
    & \besp{ \sum_{t=0}^{\horizon{-}1}
      \coutint_t\np{ \va\State_{t}, \va{\Control}^\Before_{t}, \va{\Uncertain}_{t{+}1}, \va{\Control}^\After_{t{+}1} }
      + \coutfin\np{ \va\State_{\horizon} } }
      \eqfinv
    \\
    &
      \sigma(\va{\Control}^\Before_{t}) \subset \sigma\np{\va{\State}_{0},
      \va{\Uncertain}_{1}, \ldots, \va{\Uncertain}_{t}}
      \eqsepv \forall t \in \ic{0,\horizon{-}1}
      \eqfinv
    \\
    &
      \sigma(\va{\Control}^\After_{t}) \subset \sigma\np{\va{\State}_{0},
      \va{\Uncertain}_{1}, \ldots, \va{\Uncertain}_{t}}
      \eqsepv \forall t \in \ic{1,\horizon}
      \eqfinv
    \\
    & \va\State_{t{+}1} =
      \dynamics_{t}\np{ \va\State_{t}, \va{\Control}^\Before_{t},\va{\Uncertain}_{t{+}1}, \va{\Control}^\After_{t{+}1}
      }
      \eqsepv \forall t \in \ic{0,\horizon{-}1}
      \eqfinp
  \end{align}
  \label{eq:DHD_multistage_stochastic_optimization_problem}
\end{subequations}
Thus, in the above setting,
during the time interval between two time steps, the decision-maker
makes two decisions.
At the end of the time interval~$[t-1,t[$,
a random variable~\( \va{\Uncertain}_{t} \) is revealed,
and then, at the beginning of the time interval~$[t,t+1[$,
the decision-maker makes a \emph{head decision}~$\va{\Control}_{t}^\Before$.
What is new  --- in comparison with the classical decision-hazard framework ---
is that, at the end of the time interval~$[t,t+1[$,
when a next random variable~$\va{\Uncertain}_{t{+}1}$ is revealed,
the decision-maker has the possibility to make a
\emph{tail decision}~$\va{\Control}_{t{+}1}^\After$.
This latter decision~$\va{\Control}_{t{+}1}^\After$ can be thought
as a \emph{recourse} variable for a two stage stochastic optimization problem
that would take place inside the time interval~$[t,t+1[$. Note that,
because of the term~$\np{\va{\Control}^\Before_{t},\va{\Control}^\After_{t{+}1}}$
in the cost function~$\coutint_{t}$ and in the dynamics~$\dynamics_{t}$,
considering the pair~$\np{\va{\Control}^\Before_{t},\va{\Control}^\After_{t}}$
as the control variable at time~$t$ would not satisfy the assumptions of
Theorem~\ref{cor:instantaneous-costs-in-addition-to-final-cost}.

\begin{proposition}
  If the a random variable \( \va\State_{0} \) is deterministic with value~\( \state_{0} \),
  and if the noise process $ \nseqa{\va{\Uncertain}_{t}}{t\in \ic{1,\horizon}}$ is white,
  that is, is made of independent random variables,
  then the value~\( \Value_{0}(\state_{0}) \)
  of the multistage stochastic optimization problem~\eqref{eq:DHD_multistage_stochastic_optimization_problem}
  is given by the dynamic programming backward induction
  \begin{subequations}
    \label{gdp:eq:Bellman_operators_DHD}
    \begin{align}
      \Value_{\horizon}\np{\state}
      &=
        \coutfin\np{\state}
        \eqfinv
        \intertext{and, for \( t\in\ic{0,\horizon{-}1} \),}
      \Value_{t}\np{\state}
      &=
        \inf_{\control^\Before_{t}\in\CONTROL^\Before_{t}}
        \bgesp{
        \inf_{\control_{t{+}1}^\After\in\CONTROL_{t{+}1}^\After}
        \Ba{
        \coutint_t\np{\state,\control^\Before_{t},\va{\Uncertain}_{t{+}1},\control^\After_{t{+}1}}
        +
        \Value_{t{+}1}\bp{ \dynamics_{t}\np{\state,\control^\Before_{t},\va{\Uncertain}_{t{+}1},\control^\After_{t{+}1}} }
        } }
        \eqfinp
    \end{align}
  \end{subequations}
    \label{pr:DHD_multistage_stochastic_optimization_problem}
\end{proposition}

\begin{proof}
  As the statement is made with random variables, whereas
  the theory has been developed with
  stochastic kernels, we make the link as follows:
  for each time~$\tDHD\in\ic{1,\horizonDHD}$,
  the stochastic kernel~\( \StochKernel{t-1}{t} \) in~\eqref{gdp:eq:sequence_of_stochastic_kernels}
  is the probability distribution of the random variable~\( \va{\Uncertain}_{t} \).
  This done, the proof is an application of
  Theorem~\ref{cor:instantaneous-costs-in-addition-to-final-cost},
  as follows (we just sketch the procedure, as the detailed proof can be found
  in the preprint~\cite{Carpentier-Chancelier-DeLara-Rigaut:hal-01757113}).
  We rename the uncertainty sets $\nseqa{\UNCERTAIN_t}{t \in \ic{1,\horizon}}$ as
  $\nseqa{\UNCERTAIN_t^\After}{t \in \ic{1,\horizon}}$ and
  for each time~$\tDHD\in\ic{1,\horizonDHD}$ 
  we introduce a spurious uncertainty variable~$\uncertain^\Before_\tDHD$
  taking values in a singleton set~$\UNCERTAIN^\Before_{\tDHD} =
  \{\bar\uncertain^\Before_{\tDHD}\}$, so that we obtain
  the following sequence of events
  \begin{multline*}
    \underbrace{\uncertain^\Before_{0}}_{=\state_0}
    \rightsquigarrow \control^\Before_{0}
    \rightsquigarrow \uncertain^\After_{1}
    \rightsquigarrow \control^\After_{1}
    \rightsquigarrow \underbrace{\uncertain^\Before_{1}}_{\text{spurious}}
    \rightsquigarrow \control^\Before_{1}
    \rightsquigarrow \uncertain^\After_{2}
    \rightsquigarrow \control^\After_{2}
    \rightsquigarrow \underbrace{\uncertain^\Before_{2}}_{\text{spurious}}
    \rightsquigarrow \control^\Before_{2}
    \rightsquigarrow  \quad \dots \\ \quad
    \rightsquigarrow \uncertain^\After_{\horizonDHD{-}1}
    \rightsquigarrow \control^\After_{\horizonDHD{-}1}
    \rightsquigarrow \underbrace{\uncertain^\Before_{\horizonDHD{-}1}}_{\text{spurious}}
    \rightsquigarrow \control^\Before_{\horizonDHD{-}1}
    \rightsquigarrow \uncertain^\After_{\horizonDHD}
    \rightsquigarrow \control^\After_{\horizonDHD}
    \rightsquigarrow \underbrace{\uncertain^\Before_{\horizonDHD}}_{\text{spurious}}
    \eqfinp
  \end{multline*}
  {Proceeding this way, we have doubled the timeline as time $t$ has been ``duplicated'' in
  the ordered pair $(t,\After)$ and $(t,\Before)$.}
  With this, we embed the decision-hazard-decision structure as a particular
  case of multiple consecutive time blocks (of size~2) as
  in~\S\ref{State_reduction_on_multiple_consecutive_time_blocks_and_dynamic_programming_equations}.
\end{proof}

Equation~\eqref{gdp:eq:Bellman_operators_DHD} can
be solved using the stochastic dual dynamic programming (SDDP)
algorithm provided that lower semicontinuity and convexity of
the value functions are preserved. This is ensured first by
assuming linearity in the dynamics -- a feature that may be obtained
by modeling the problem in the decision-hazard-decision framework
as illustrated in~\S\ref{gdp:Decision_Hazard_Decision_Motivation} --
and second by assuming lower semicontinuity for the cost
functions as well as compactness for the existence
of optimal controls.
% (see~\cite[Theorem~3.4]{LeFranc-Carpentier-Chancelier-DeLara:2022}).

\section{Conclusion and perspectives}

As said in the introduction, decomposition methods are appealing
to tackle multistage stochastic optimization problems, as they
are naturally large scale. The most common approaches are time
decomposition (and state-based resolution methods, like stochastic
dynamic programming, in stochastic optimal control), and scenario
decomposition (like progressive hedging in stochastic programming).

This paper is part of a general research program that consists
in \emph{mixing} different decomposition bricks. Space decomposition
methods have been investigated in  \cite{Barty_RAIRO_2010}
and \cite{Carpentier-Chancelier-DeLara-Pacaud:2020}.
Here, we have tackled the issue of using time block decomposition
in such a way that stochastic dynamic programming is used
at the slow time scale with an appropriate white noise assumption,
whereas stochastic programming methods such as progressive hedging
can be used at the fast time scale where such an independence assumption
does not hold. This approach paves the way of mixing time decomposition
with scenario decomposition.
For this purpose, we have revisited the notion of state,
and have provided a way to perform time decomposition
but only across specified time blocks.

\bigskip

\textbf{Acknowledgements.} We thank Roger Wets for fruitful
discussions about the possibility of mixing stochastic dynamic
programming with progressive hedging.

\appendix

\section{Technical details and proofs}
% (Sect.~\ref{gdp:Stochastic_Dynamic_Programming_and_State_Reduction_by_Time_Blocks})}
\label{Technical_Details_and_Proofs}

We suppose to be in the framework of \S\ref{The_Bertsekas-Shreve_setting}.
We introduce the notations
\begin{equation}
   \UNCERTAIN_{\tun:\tter}
  = \prod_{\tbis=\tun}^{\tter} \UNCERTAIN_{\tbis} \eqsepv
    0 \leq \tun \leq \tter \leq \horizon
    \eqsepv
  \CONTROL_{\tun:\tter}
    = \prod_{\tbis=\tun}^{\tter} \CONTROL_{\tbis} \eqsepv
    0 \leq \tun \leq \tter \leq \horizon{-}1 % \\
  \eqfinp
\end{equation}
\begin{subequations}
  Let $ 0 \leq \tun \leq \tbis \leq \tter \leq \horizon $.
  From a history~$\history_{\tter} \in \HISTORY_{\tter}$,
  we extract
  the $\interval{\tun}{\tbis}$-\emph{history uncertainty part}
  \begin{equation}
    \nc{\history_{\tter}}_{\tun:\tbis}^{\UNCERTAIN}=
    \np{\uncertain_{\tun},\ldots, \uncertain_{\tbis} } = \uncertain_{\tun:\tbis}
    \in \UNCERTAIN_{\tun:\tbis}
    \eqsepv 0 \leq \tun \leq \tbis \leq \tter
    \eqfinv
    \label{gdp:eq:uncertain_t:r}
  \end{equation}
  the $\interval{\tun}{\tbis}$-\emph{history control part}
  (notice that the indices are special)
  \begin{align}
    \nc{\history_{\tter}}_{\tun:\tbis}^{\CONTROL}
%    &
      =
      \np{\control_{\tun-1},\ldots, \control_{\tbis-1} }
    %   \nonumber \\
    % &
      = \control_{\tun-1:\tbis-1}
      \in \CONTROL_{\tun-1:\tbis-1}
      \eqsepv 1 \leq \tun \leq \tbis \leq \tter
      \eqfinp
      \label{gdp:eq:control_t:r}
  \end{align}
\end{subequations}

\subsubsubsection{Flows}

Let $\tun$ and $\tter$ be given such that
$ 0 \leq \tun < \tter \leq \horizon $.
For a $\interval{\tun}{\tter-1}$-history feedback
$\HistoryFeedback
=\nseqa{\HistoryFeedback_{\tbis}}{\tbis\in\ic{\tun,\tter}-1}
\in\HISTORYFEEDBACK_{\tun:\tter-1} $,
we define the \emph{flow} $\Flow{\tun}{\tter}{\HistoryFeedback}$ by
\begin{subequations}
  \begin{align}
    \Flow{\tun}{\tter}{\HistoryFeedback}
    : \HISTORY_{\tun} \times \UNCERTAIN_{\tun{+}1:\tter}
    & \to  \HISTORY_{\tter}
      \label{gdp:eq:flow-def} \\
    \np{ \history_{\tun}, \uncertain_{\tun{+}1:\tter} }
    & \mapsto
      \Big(\history_{\tun},
      \HistoryFeedback_{\tun}(\history_{\tun}),\uncertain_{\tun{+}1},
                       \HistoryFeedback_{\tun{+}1}\bp{\history_{\tun},
                       \HistoryFeedback_{\tun}(\history_{\tun}),
                       \uncertain_{\tun{+}1}},\uncertain_{\tun{+}2},
                                      \cdots, \HistoryFeedback_{\tter-1}(\history_{\tter-1}),
                                      \uncertain_{\tter}\Big)
                                      \eqfinp \nonumber
  \end{align}
  {Otherwise stated, the flow is given by}
  \begin{align}
    \Flow{\tun}{\tter}{\HistoryFeedback}
    \np{ \history_{\tun}, \uncertain_{\tun{+}1:\tter} }
    &=
      \np{\history_{\tun}, \control_{\tun},\uncertain_{\tun{+}1},
      \control_{\tun{+}1}, \uncertain_{\tun{+}2}, \ldots,
      \control_{\tter-1},\uncertain_{\tter}} \eqfinv
    \\
    \mtext{with } \history_{\tbis}
    &=\np{\history_{\tun}, \control_{\tun},\uncertain_{\tun{+}1},
      \ldots, \control_{\tbis-1}, \uncertain_{\tbis} }
      \eqsepv \tun < \tbis \leq \tter  \eqfinv
    \\
    \mtext{and } \control_{\tbis}
    &= \HistoryFeedback_{\tbis}
      \np{\history_{\tbis}}
      \eqsepv \tun \leq \tbis \leq \tter-1
      \eqfinp
  \end{align}
  When $ 0 \leq \tun = \tter \leq \horizon $, we put
  $ \Flow{\tun}{\tun}{\HistoryFeedback} : \HISTORY_{\tun}
  \to  \HISTORY_{\tun} $, $ \history_{\tun}  \mapsto
  \history_{\tun}$.
  \label{gdp:eq:flow}
\end{subequations}
With this convention, the expression $\Flow{\tun}{\tter}{\HistoryFeedback}$
makes sense when $ 0 \leq \tun \leq \tter \leq \horizon $.
The mapping $\Flow{\tun}{\tter}{\HistoryFeedback}$
gives the history at time~$\tter$ as a
function of the initial history~$\history_{\tun}$ at time~$\tun$
and of the history feedbacks
$\nseqa{\HistoryFeedback_{\tbis}}{\tbis
 \in\ic{\tun,\tter}-1}\in\HISTORYFEEDBACK_{\tun:\tter-1} $.
\begin{subequations}
  \label{gdp:eq:flowconnseqas}
  An immediate consequence of this definition are the \emph{flow properties}:
  \begin{align}
    \Flow{\tun}{\tter+1}{\HistoryFeedback}
    \np{ \history_{\tun}, \uncertain_{\tun{+}1:\tter+1}}
    &=
      \Bp{ \Flow{\tun}{\tter}{\HistoryFeedback}
      \np{ \history_{\tun}, \uncertain_{\tun{+}1:\tter}},
      \HistoryFeedback_{\tter}\bp{\Flow{\tun}{\tter}{\HistoryFeedback}
      \np{ \history_{\tun}, \uncertain_{\tun{+}1:\tter}}},\uncertain_{\tter+1}}
      \eqsepv
                         0 \leq \tun \leq \tter \leq \horizon{-}1
                         \eqfinv
                         \label{gdp:eq:flowconnseqa1}
    \\
    \Flow{\tun}{\tter}{\HistoryFeedback}
    \np{\history_{\tun},\uncertain_{\tun{+}1:\tter}}
    &=
      \Flow{\tun{+}1}{\tter}{\HistoryFeedback}
      \Bp{\bp{\history_{\tun},\HistoryFeedback_{\tun}(\history_{\tun}),
      \uncertain_{\tun{+}1}},\uncertain_{\tun{+}2:\tter}}
      \eqsepv 0 \leq \tun < \tter \leq \horizon
      \eqfinp
      \label{gdp:eq:flowconnseqa2}
  \end{align}
\end{subequations}

We recall that $\espace{L}^{0}_{+}(\HISTORY_{t}) $ denotes
    the space of lower semianalytic nonnegative numerical functions over~$\HISTORY_{t}$.
\begin{definition}
  \label{gdp:de:stochastic_kernels_rho}
  Let $\tun$ and $\tter$ be given such that $0 \leq \tun \leq \tter \leq \horizon$.
  \begin{itemize}
  \item When $ 0 \leq \tun < \tter \leq \horizon $, for
    a $\interval{\tun}{\tter-1}$-history feedback
    $\HistoryFeedback
    =\nseqa{\HistoryFeedback_{\tbis}}{\tbis\in \ic{\tun,\tter-1}}
    \in\HISTORYFEEDBACK_{\tun:\tter-1} $,
    and for a family $ \nseqa{ \StochKernel{\tbis-1}{\tbis}}%
    {\tbis\in\ic{\tun{+}1,\tter}} $ of Borel-measurable stochastic kernels
    $       \StochKernel{\tbis-1}{\tbis} :  \HISTORY_{\tbis-1} \to \Delta\np{\UNCERTAIN_{\tbis}}
    \eqsepv \tbis \in \ic{\tun{+}1,\tter} $,
    we define a Borel-measurable stochastic kernel
    $ \StochKernelFeed{\tun}{\tter}{\HistoryFeedback} :
      \HISTORY_{\tun} \to \Delta\np{ \HISTORY_{\tter} }
    $
    such that, for any numerical function
    $ \varphi \in \espace{L}^{0}_{+}(\HISTORY_{t}) $, we have that
    \begin{align}
      \int_{\HISTORY_{\tter}}
      & \varphi\np{\history'_{\tun},\history'_{\tun{+}1:\tter}}
        \StochKernelFeed{\tun}{\tter}{\HistoryFeedback}
        \kernelargs{\history_{\tun}}{\dd\history'_{\tter}}
        \nonumber \\
      &=
        \int_{\UNCERTAIN_{\tun{+}1:\tter}}
        \varphi \bp{ \Flow{\tun}{\tter}{\HistoryFeedback} \np{\history_{\tun}, \uncertain_{\tun{+}1:\tter} }}
        \prod_{\tbis=\tun{+}1}^{\tter}
        \StochKernel{\tbis-1}{\tbis}
        \bkernelargs{ \Flow{\tun}{\tbis-1}{\HistoryFeedback}
        \np{\history_{\tun},\uncertain_{\tun{+}1:\tbis-1}}}
        {\dd\uncertain_{\tbis}}
        \eqfinp \label{gdp:eq:stochastic_kernels_rho_b}
    \end{align}
  \item
    When $ 0 \leq \tun = \tter \leq \horizon $, we define
    $
    \StochKernelFeed{\tun}{\tun}{\HistoryFeedback} :
    \HISTORY_{\tun} \to \Delta\np{ \HISTORY_{\tun} } $
    by $  \StochKernelFeed{\tun}{\tun}{\HistoryFeedback} \kernelargs{\history_{\tun}}{ \dd\history'_{\tun} } =
    \delta_{\history_{\tun}}\np{ \dd\history'_{\tun} } $
    where~$\delta$ represents the Dirac measure.
  \end{itemize}
\end{definition}
The stochastic kernels~$ \StochKernelFeed{\tun}{\tter}{\HistoryFeedback} $
on~$\HISTORY_{\tter} $, given by~\eqref{gdp:eq:stochastic_kernels_rho_b},
are of the form
\begin{equation}
\label{eq:kernelfactorization}
\StochKernelFeed{\tun}{\tter}{\HistoryFeedback}
\kernelargs{\history_{\tun}}{\dd\history'_{\tter}} =
\StochKernelFeed{\tun}{\tter}{\HistoryFeedback}
\kernelargs{\history_{\tun}}{\dd\history'_{\tun}\dd\history'_{\tun{+}1:\tter}} =
\delta_{\history_{\tun}}\np{ \dd\history'_{\tun} }
\otimes  \varStochKernelFeed{\tun}{\tter}{\HistoryFeedback}
\kernelargs{\history_{\tun}}{\dd\history'_{\tun{+}1:\tter}} \eqfinv
\end{equation}
where, for each $ \history_{\tun} \in \HISTORY_{\tun} $,
the probability distribution~$ \varStochKernelFeed{\tun}{\tter}{\HistoryFeedback}
\kernelargs{\history_{\tun}}{\dd\history'_{\tun{+}1:\tter}} $ only charges
the histories visited by the flow from~$\tun{+}1$ to~$\tter$. The construction
of the stochastic kernels~$ \StochKernelFeed{\tun}{\tter}{\HistoryFeedback}$
is developed in \cite[p. 190]{Bertsekas-Shreve:1996} for relaxed
history feedbacks and obtained by using
\cite[Proposition~7.45, p.~175]{Bertsekas-Shreve:1996}.

\begin{proposition}
  \label{gdp:pr:stochastic_kernels_rho}
  The family
  $\nseqa{ \StochKernelFeed{\tbis}{\tter}{\HistoryFeedback}}
  {\tbis\in\ic{\tun,\tter}} $ of stochastic kernels of
  Definition~\ref{gdp:de:stochastic_kernels_rho}
  has the flow property: % for $ \tbis < \tter $,
  \begin{equation}
    \StochKernelFeed{\tbis}{\tter}{\HistoryFeedback}
    \kernelargs{\history_{\tbis}}{\dd\history'_{\tter}}
    =
    \int_{\UNCERTAIN_{\tbis{+}1}}
    \StochKernel{\tbis}{\tbis{+}1}\kernelargs{\history_{\tbis}}{\dd\uncertain_{\tbis{+}1} }
    \StochKernel{\tbis{+}1}{\tter}^{\HistoryFeedback}
    \Bkernelargs{\bp{\history_{\tbis},\HistoryFeedback_{\tbis}\np{\history_{\tbis}},
                     \uncertain_{\tbis{+}1}}}
    {\dd\history'_{\tter}}
    \eqsepv \forall \tbis < \tter
    \eqfinp
    \label{gdp:eq:stochastic_kernels_rho_property}
  \end{equation}
\end{proposition}

\begin{proof}
  Let $ \tbis < \tter $.
  For any $\varphi \in \espace{L}^{0}_{+}(\HISTORY_{t})$, we have that
  \begin{subequations}
    \begin{align}
      \int_{\HISTORY_{\tter}}
      \varphi
      & \np{\history'_{\tbis},\history'_{\tbis{+}1:\tter}}
        \StochKernelFeed{\tbis}{\tter}{\HistoryFeedback}\bkernelargs{\history_{\tbis}}{ \dd\history'_{\tter}}
        \label{gdp:eq:StochasticKernel_Flow_Ini}
      \\
      =& \int_{\UNCERTAIN_{\tbis{+}1:\tter}}
         \varphi \bp{ \Flow{\tbis}{\tter}{\HistoryFeedback}
         \np{\history_{\tbis}, \uncertain_{\tbis{+}1:\tter} } }
\int_{\UNCERTAIN_{\tbis{+}1:\tter}}
          \prod_{{\tbis'}=\tbis{+}1}^{\tter}
        \StochKernel{{\tbis'}-1}{{\tbis'}}
        \bkernelargs{ \Flow{\tbis}{{\tbis'}-1}{\HistoryFeedback}\np{\history_{\tbis},
        \uncertain_{\tbis{+}1:{\tbis'}-1} }}{\dd\uncertain_{{\tbis'}} }
        \nonumber
\intertext{by Definition~\eqref{gdp:eq:stochastic_kernels_rho_b}}
      =&
         \int_{\UNCERTAIN_{\tbis{+}1:\tter}}
         \varphi \bp{ \Flow{\tbis}{\tter}{\HistoryFeedback}
         \np{\history_{\tbis}, \uncertain_{\tbis{+}1:\tter} } }
         \StochKernel{\tbis}{\tbis{+}1}\bkernelargs{ \history_{\tbis}}{\dd\uncertain_{\tbis{+}1} }
          \int_{\UNCERTAIN_{\tbis{+}1:\tter}}
          \prod_{{\tbis'}=\tbis+2}^{\tter}
        \StochKernel{{\tbis'}-1}{{\tbis'}}\bkernelargs{ \Flow{\tbis}{{\tbis'}-1}{\HistoryFeedback}\np{\history_{\tbis},
        \uncertain_{\tbis{+}1:{\tbis'}-1} }}{\dd\uncertain_{{\tbis'}} }
        \nonumber \\
      =&
         \int_{\UNCERTAIN_{\tbis{+}1:\tter}}
         \varphi \bp{\Flow{\tbis{+}1}{\tter}{\HistoryFeedback}
         \bp{\np{\history_{\tbis},\HistoryFeedback_{\tbis}(\history_{\tbis}),
                 \uncertain_{\tbis{+}1}},\uncertain_{\tbis+2:\tter}}}
        \StochKernel{\tbis}{\tbis{+}1}\bkernelargs{ \history_{\tbis}}{\dd\uncertain_{\tbis{+}1} }
        \nonumber \\
      & \;\;
        \prod_{{\tbis'}=\tbis+2}^{\tter}
        \StochKernel{{\tbis'}-1}{{\tbis'}}\bkernelargs{ \Flow{\tbis{+}1}{{\tbis'}-1}{\HistoryFeedback}
        \bp{\np{\history_{\tbis},\HistoryFeedback_{\tbis}(\history_{\tbis}),
                \uncertain_{\tbis{+}1}},\uncertain_{\tbis+2:{\tbis'}-1}}
        }{\dd\uncertain_{{\tbis'}} }
        \tag{by the flow property~\eqref{gdp:eq:flowconnseqa2}}
      \\
      =&
         \int_{\UNCERTAIN_{\tbis{+}1}}
         \StochKernel{\tbis}{\tbis{+}1}\bkernelargs{ \history_{\tbis}}{\dd\uncertain_{\tbis{+}1} }
        \int_{\UNCERTAIN_{\tbis+2:\tter}}
        \varphi \bp{\Flow{\tbis{+}1}{\tter}{\HistoryFeedback}
        \bp{\np{\history_{\tbis},\HistoryFeedback_{\tbis}(\history_{\tbis}),
                \uncertain_{\tbis{+}1}},\uncertain_{\tbis+2:\tter}}}
        \nonumber \\
      & \;\;  \prod_{{\tbis'}=\tbis+2}^{\tter}
        \StochKernel{{\tbis'}-1}{{\tbis'}}\bkernelargs{ \Flow{\tbis{+}1}{{\tbis'}-1}{\HistoryFeedback}
        \bp{\np{\history_{\tbis},\HistoryFeedback_{\tbis}(\history_{\tbis}),
                \uncertain_{\tbis{+}1}},\uncertain_{\tbis+2:{\tbis'}-1}}
        }{\dd\uncertain_{{\tbis'}} }
        \nonumber %\\
      \intertext{by Fubini Theorem} %\cite[p.137]{Loeve:1977}}
      % \\
      =& \int_{\UNCERTAIN_{\tbis{+}1}} \StochKernel{\tbis}{\tbis{+}1}\bkernelargs{\history_{\tbis}}{\dd\uncertain_{\tbis{+}1} }
         \int_{\HISTORY_{\tter}} \varphi\bp{
         \np{\history'_{\tbis},\HistoryFeedback_{\tbis}(\history'_{\tbis}),
             \uncertain'_{\tbis{+}1}}, \history'_{\tbis+2:\tter}}
        \StochKernel{\tbis{+}1}{\tter}^{\HistoryFeedback} \bkernelargs{
        \np{\history_{\tbis},\HistoryFeedback_{\tbis}(\history_{\tbis}),
            \uncertain_{\tbis{+}1}}} {\dd\history'_{\tter}}
        \nonumber
      \intertext{by Definition~\eqref{gdp:eq:stochastic_kernels_rho_b}}
      =&
         \int_{\HISTORY_{\tter}} \varphi\bp{
         \np{\history'_{\tbis},\HistoryFeedback_{\tbis}(\history'_{\tbis}),
             \uncertain'_{\tbis{+}1}},\history'_{\tbis+2:\tter}}
         \nonumber \\
      &
        \;\;\int_{\UNCERTAIN_{\tbis{+}1}}
        \StochKernel{\tbis}{\tbis{+}1}
        \bkernelargs{\history_{\tbis}}{\dd\uncertain_{\tbis{+}1}}
        \StochKernel{\tbis{+}1}{\tter}^{\HistoryFeedback}
        \bkernelargs{\np{\history_{\tbis},\HistoryFeedback_{\tbis}(\history_{\tbis}),
                         \uncertain_{\tbis{+}1}}}{\dd\history'_{\tter}}
        \label{gdp:eq:StochasticKernel_Flow_Last}
    \end{align}
    by Fubini Theorem.
  \end{subequations}
  As the two expressions~\eqref{gdp:eq:StochasticKernel_Flow_Ini}
  and~\eqref{gdp:eq:StochasticKernel_Flow_Last} are equal for any
  $\varphi \in \espace{L}^{0}_{+}(\HISTORY_{t})$,
  we deduce the flow property
  \eqref{gdp:eq:stochastic_kernels_rho_property}.
\end{proof}

\subsubsubsection{Proof of Theorem~\ref{gdp:pr:DP_withoutstate_third}}

We only give a sketch of the proof, as it is a variation on
different results of~\cite{Bertsekas-Shreve:1996},
the framework of which we follow.

\begin{proof}
  We are in the setting of~\cite[Chap.~7,~Chap.~8]{Bertsekas-Shreve:1996}.
  We take the history space $\HISTORY_{t}$ for state space,
  and the state dynamics
  \begin{equation}
    f\bp{\history_{t},\control_{t},\uncertain_{t{+}1}}
    =\bp{\history_{t},\control_{t},\uncertain_{t{+}1}}=
    \history_{t{+}1}
    \in \HISTORY_{t{+}1} = \HISTORY_{t} \times \CONTROL_{t} \times \UNCERTAIN_{t{+}1}
    \eqfinp
    \label{gdp:eq:history=state_dynamics}
  \end{equation}
  Then, the family
  $ \nseqa{ \StochKernel{\tbis-1}{\tbis}}{\tbis\in\ic{1,\horizon}}$ of Borel-measurable
  stochastic kernels~\eqref{gdp:eq:sequence_of_stochastic_kernels}
  gives a family of disturbance kernels (vocabulary of~\cite[p~189]{Bertsekas-Shreve:1996})
  that do not depend on the current control.
  The criterion to be minimized~\eqref{gdp:eq:criterion}
  is a function of the history at time $\horizon$,
  thus of the state at time $\horizon$. We consider the finite
  horizon model with final cost corresponding to the optimization
  problem  defined by the associated value function~\eqref{eq:Vzero}:
  \begin{align*}
    \Value_{0}(\uncertain_{0})
    & = \inf_{\HistoryFeedback_{0:\horizon{-}1} \in \HISTORYFEEDBACK_{0:\horizon{-}1}}
      \int_{\HISTORY_{\horizon}} \criterion\np{\history'_{\horizon}}
      \StochKernelFeed{0}{\horizon}{\HistoryFeedback}
      \kernelargs{\uncertain_{0}}{\dd \history'_{\horizon}}
      \nonumber
    \\
    & = \inf_{\HistoryFeedback_{0:\horizon{-}1} \in \HISTORYFEEDBACK_{0:\horizon{-}1}}
      \int_{\UNCERTAIN_{1:\horizon}}
      \criterion
      \bp{ \Flow{0}{\horizon}{\HistoryFeedback} \np{\uncertain_{0:\horizon} }}
      \prod_{\tbis=1}^{\horizon}
      \StochKernel{\tbis-1}{\tbis}
      \bkernelargs{ \Flow{0}{\tbis-1}{\HistoryFeedback}\np{\uncertain_{0:\tbis-1}}}
      {\dd\uncertain_{\tbis} }
      \eqfinv
      \tag{by~\eqref{gdp:eq:stochastic_kernels_rho_b}}
  \end{align*}
  where the flows $\Flow{0}{\tbis}{\HistoryFeedback}$ for
  $\tbis \in \ic{0,\horizon{-}1}$ are defined by
  Equation~\eqref{gdp:eq:flow-def}, and
where we are minimizing over the so-called state-feedbacks.
  Then, the proof of Theorem~\ref{gdp:pr:DP_withoutstate_third}
  follows from the results developed in Chap.~7,~8 and~10
  of~\cite{Bertsekas-Shreve:1996} in a Borel setting.

  The Bellman operators in~\eqref{gdp:eq:Bellman_operators_rho}
  satisfy~\eqref{gdp:eq:Bellman_operators_rho_lower_semianalytic_preservation}
  because, by Lemma 7.30(4) and Propositions 7.47 and 7.48 in
  \cite{Bertsekas-Shreve:1996},
  we have that\footnote{More precisely, the
  property~\eqref{gdp:eq:Bellman_operators_rho_lower_semianalytic_preservation}
  is a consequence of the properties that i) the Bellman
  operator~\eqref{gdp:eq:Bellman_operators_rho} corresponds to the operator~$T$
  (with $g=0$ and $\alpha=1$) in \cite[Definition~8.5,
  p.~195]{Bertsekas-Shreve:1996}
  ii) that $T(J)$ is lower semianalytic whenever $J$ also is,
  as explained right above \cite[Lemma~8.2, p.~196]{Bertsekas-Shreve:1996}.}
  \( \varphi\in  \espace{L}^{0}_{+}(\HISTORY_{t{+}1}) \)
  $\implies$ $\BelOp{t+1}{t}\varphi \in \espace{L}^{0}_{+}(\HISTORY_{t})$,
  for $t$ in $\ic{0,\horizon{-}1}$.

  Since we are considering a finite horizon model with a final cost,
  we detail the steps needed to use the results
  of~\cite[Chap.~8]{Bertsekas-Shreve:1996}.
  The final cost at time~$\horizon$ can be turned into an instantaneous cost
  at time $\horizon{-}1$~by inserting
  the state dynamics~\eqref{gdp:eq:history=state_dynamics}
  in the final cost.
  Getting rid of the disturbance in the expected cost
  by using the disturbance kernel is standard practice.
  Then, we can turn this non-homogeneous finite horizon model
  into a finite horizon model with
  homogeneous dynamics and costs
  by following the steps of~\cite[Chap.~10]{Bertsekas-Shreve:1996}.
  Using~\cite[Proposition~8.2, p.~198]{Bertsekas-Shreve:1996},
  we obtain that the family of optimization problems
  defined by the associated value functions~\eqref{gdp:eq:value_functions},
  when minimizing over the relaxed state feedbacks,
  satisfies the Bellman equation~\eqref{gdp:eq:Bellman_equation};
  we conclude
with~\cite[Proposition~8.4, p.~203]{Bertsekas-Shreve:1996}
  which covers the minimization over state feedbacks.

The Bellman equation~\eqref{gdp:eq:Bellman_equation} is a consequence of
  \cite[Proposition~8.2, p.~198]{Bertsekas-Shreve:1996}.

To finish, Theorem~\ref{gdp:pr:DP_withoutstate_third} is valid
under the general Borel assumptions of~\cite[Chap.~8]{Bertsekas-Shreve:1996}
and with the specific $(F^{+})$ assumption needed
for~\cite[Proposition~8.4, p.~203]{Bertsekas-Shreve:1996}; this last assumption
is fulfilled here since we have assumed that the
cost criterion~\eqref{gdp:eq:criterion} is nonnegative.
\end{proof}

%\section{Technical details and proofs (Sect.~\ref{gdp:State_Reduction_by_Time_Blocks})}
%\label{Technical_Details_and_Proofs_gdp:State_Reduction_by_Time_Blocks}

\subsubsubsection{Proof of Proposition~\ref{gdp:thm:DPB}}\quad

We suppose to be in the framework of \S\ref{The_Bertsekas-Shreve_setting}.

\begin{proof}
  We are in the setting of~\cite[Chap.~7,~Chap.~8]{Bertsekas-Shreve:1996}.
  Let~$\tilde\varphi_{\tter}: \STATE_{\tter} \to \bRR$ be a
  given lower semianalytic nonnegative numerical function,
  and let~$\varphi_{\tter}: \HISTORY_{\tter} \to \bRR$ be
  \begin{equation}
    \varphi_{\tter} = \tilde\varphi_{\tter} \circ \theta_{\tter} \eqfinp
    \label{gdp:eq:varphi}
  \end{equation}
  Let~$\varphi_{\tun} : \HISTORY_{\tun} \to \bRR$
  be the lower semianalytic nonnegative numerical function obtained
  by applying the Bellman operator $\BelOp{\tter}{\tun}$
  across~$\icinterval{\tter}{\tun}$ (see~\eqref{gdp:eq:BelOp_t:r})
  to the lower semianalytic nonnegative numerical function~$\varphi_{\tter}$:
  \begin{equation}
    \varphi_{\tun} = \BelOp{\tter}{\tun} \varphi_{\tter} =
    \BelOp{\tun{+}1}{\tun} \circ \cdots \circ
    \BelOp{\tter}{\tter-1} \varphi_{\tter} \eqfinp
    \label{gdp:eq:varphi_tun}
  \end{equation}
  By \cite[Lemma 7.30(3), p.~178]{Bertsekas-Shreve:1996} --- on the
  stability of lower semianalytic functions under right composition
  with a Borel-measurable mapping --- we get that the nonnegative
  numerical function~$\varphi_{\tun}$ is lower semianalytic.
  We show that there exists a lower semianalytic nonnegative numerical
  function $\tilde\varphi_{\tun} : \STATE_{\tun} \to \bRR$ such that
  \begin{equation}
    \varphi_{\tun} =  \tilde\varphi_{\tun} \circ \theta_{\tun}  \eqfinp
    \label{gdp:eq:proof_DPB}
  \end{equation}

  First, we show by backward induction that,
  for all $\tbis \in \ic{\tun,\tter}$,
  there exists a lower semianalytic nonnegative numerical
  function~$\overline\varphi_{\tbis}$ such that
  $\varphi_{\tbis}(h_{\tbis})= \overline\varphi_{\tbis}
  ( \theta_{\tun} \np{\history_{\tun}},\history_{\tun{+}1:\tbis})$.
  Second, we prove that the function~$\tilde\varphi_{\tun}=\overline\varphi_{\tun}$
  satisfies~\eqref{gdp:eq:proof_DPB} and is lower semianalytic.
  \begin{itemize}
  \item
    For $\tbis=\tter$, we have, by~\eqref{gdp:eq:varphi} and
    by~\eqref{gdp:eq:reduction-dynamics}, that
    $
    \varphi_{\tter}(\history_{\tter}) =
    \tilde\varphi_{\tter}\bp{\theta_{\tter} \np{\history_{\tter}}} =
    \tilde\varphi_{\tter}\bp{\Dynamics{\tun}{\tter}\np{\theta_{\tun}\np{\history_{\tun}},
        \history_{\tun{+}1:\tter}}} $,
    so that the nonnegative numerical function~$\overline\varphi_{\tter}$
    is given by
    $\tilde\varphi_{\tter} \circ \Dynamics{\tun}{\tter}$.
    By \cite[Lemma 7.30(3), p.~178]{Bertsekas-Shreve:1996},
    $\overline\varphi_{\tter}$ is a lower semianalytic numerical function.
  \item
    Assume that, at~$\tbis{+}1$, the result holds true, that is,
    $ \varphi_{\tbis{+}1}(h_{\tbis{+}1})=
    \overline\varphi_{\tbis{+}1}(\theta_{\tun}\np{\history_{\tun}},
    \history_{\tun{+}1:\tbis{+}1}) $,
    where the numerical function~$\overline\varphi_{\tbis{+}1}$
    is nonnegative lower semianalytic.
    Then, by~\eqref{gdp:eq:varphi_tun},
    \begin{subequations}
      \begin{align*}
        \varphi_{\tbis}(\history_{\tbis})
        & =
          \bp{\BelOp{\tbis{+}1}{\tbis} \varphi_{\tbis{+}1}}(\history_{\tbis})
        \\
        &
          =
          \inf_{\control_{\tbis}\in\CONTROL_{\tbis}} \int_{\UNCERTAIN_{\tbis{+}1}}
          \varphi_{\tbis{+}1}\bp{\np{\history_{\tbis},\control_{\tbis},\uncertain_{\tbis{+}1}}}
          \StochKernel{\tbis}{\tbis{+}1}\nkernelargs{\history_{\tbis}}{\dd\uncertain_{\tbis{+}1} }
          \tag{by definition~\eqref{gdp:eq:Bellman_operators_rho} of the Bellman operator }
        \\
        & =
          \inf_{\control_{\tbis}\in\CONTROL_{\tbis}} \int_{\UNCERTAIN_{\tbis{+}1}}
          \overline\varphi_{\tbis{+}1}\bp{\np{\theta_{\tun} \np{\history_\tun},
          \np{\history_{\tun{+}1:\tbis},\control_{\tbis},\uncertain_{\tbis{+}1}}}}
            \inf_{\control_{\tbis}\in\CONTROL_{\tbis}} \int_{\UNCERTAIN_{\tbis{+}1}}
            \StochKernel{\tbis}{\tbis{+}1}\nkernelargs{\history_{\tbis}}
            {\dd\uncertain_{\tbis{+}1}}
          \tag{by the induction assumption}
        \\ %
        & =
          \inf_{\control_{\tbis}\in\CONTROL_{\tbis}} \int_{\UNCERTAIN_{\tbis{+}1}}
          \overline\varphi_{\tbis{+}1}\bp{\np{\theta_{\tun}\np{\history_\tun},
          \np{\history_{\tun{+}1:\tbis},\control_{\tbis},\uncertain_{\tbis{+}1}}}}
            \inf_{\control_{\tbis}\in\CONTROL_{\tbis}} \int_{\UNCERTAIN_{\tbis{+}1}}
            \tildeStochKernel{\tbis}{\tbis{+}1}
          \bkernelargs{\np{\theta_{\tun}\np{\history_{\tun}},\history_{\tun{+}1:\tbis}}}
          {\dd\uncertain_{\tbis{+}1}}
          \tag{by compatibility of the stochastic kernel}
        \\ %
        & =
          \overline\varphi_{\tbis}\bp{\theta_{\tun}\np{\history_{\tun}},\history_{\tun{+}1:\tbis}}
          \eqfinv
      \end{align*}
    \end{subequations}
    \[
       \textrm{where } \quad
       \overline\varphi_{\tbis}\bp{\optstate_{\tun},\history_{\tun{+}1:\tbis}}
        = \inf_{\control_{\tbis}\in\CONTROL_{\tbis}} \int_{\UNCERTAIN_{\tbis{+}1}}
        \overline\varphi_{\tbis{+}1}\bp{\np{\optstate_{\tun},
        \np{\history_{\tun{+}1:\tbis},\control_{\tbis},\uncertain_{\tbis{+}1}}}}
        \tildeStochKernel{\tbis}{\tbis{+}1}
        \bkernelargs{\np{\optstate_{\tun},\history_{\tun{+}1:\tbis}}}
        {\dd\uncertain_{\tbis{+}1}}
        \eqfinp
    \]
    By \cite[p.~196]{Bertsekas-Shreve:1996} (right before Lemma~8.2),
    we get that the numerical function~$\overline\varphi_{\tbis}$
    is nonnegative lower semianalytic.
    Thus, we have shown that the result holds true at time~$\tbis$.
  \end{itemize}
  The induction implies that, at time~$\tun$, the expression
  of~$\varphi_{\tun}(\history_{\tun})$ is
  $\varphi_{\tun}(\history_{\tun}) =
  \overline\varphi_{\tun}\bp{\theta_{\tun}\np{\history_{\tun}}}$,
  since the term ~$\history_{\tun{+}1:\tun}$ vanishes.
  Choosing~$\tilde\varphi_{\tun}=\overline\varphi_{\tun}$
  gives the expected result.
\end{proof}

\ifx\undefined\noopsort\newcommand{\noopsort}[1]{}\fi\ifx\undefined\allcaps\def\allcaps#1{#1}\fi

\ifpreprint
\newpage
\preprintstart

\begin{center}
  \small
  {\bfseries Forewords for additional material} \\[0.2cm]
  \begin{minipage}{0.8\textwidth}
    We provide here additional material to the paper~\citeAM{tts:carpentier:JOCA}.
    In Appendix~\ref{gdp:Background_on_Stochastic_Dynamic_Programming},
    we survey several frameworks and approaches to solve, by dynamic programming,
    a stochastic optimal control problem formulated in discrete time.
    In Appendix~\ref{Link_Classical_Framework}, we make the link between
    the setting of two-time-scale optimization problems
    (as developed in  Sect.~\ref{Two_time_scale_optimization_problem} with
    stochastic kernels)
    and the framework of stochastic optimal control
    (with random variables). 
    In Appendix~\ref{Technical_details_and_proofs_DHD},
we give a detailed proof of
Proposition~\ref{pr:DHD_multistage_stochastic_optimization_problem}
concerning the decision-hazard-decision approach. 
In Appendix~\ref{Two-time-scale_dynamic_programming}, 
we present a framework for two-time-scale multistage optimization problems
which is more general 
than in Sect.~\ref{Two_time_scale_optimization_problem}. 
  \end{minipage}
\end{center}

\section{A brief survey of frameworks for stochastic dynamic programming in discrete time}
\label{gdp:Background_on_Stochastic_Dynamic_Programming}

We sketch mathematical frameworks for stochastic
dynamic programming in discrete time to be found in the literature.
In what follows, $\tinitial \in \NN$ and $\horizon \in \NN^*$
are two natural numbers such that $\tinitial < \horizon$.
We use the notation $ \ic{r,s}=\na{r,r+1,\ldots,s-1,s} $ for any two
natural numbers~$r,s$ such that $r \leq s$.

\subsubsubsection{Witsenhausen approach}

The most general stochastic dynamic programming principle is sketched by
Witsenhausen at the end of~\citeAM{Witsenhausen:1975b}.
However, we do not detail it as its formalism is
too far from the following ones, though we will touch the subject when we
discuss Y\"{u}ksel's approach below.
We present here what Witsenhausen calls
an optimal stochastic control problem in \emph{standard form}
(see~\citeAM{Witsenhausen:1973}).
The ingredients are the following:
\begin{enumerate}
\item
  time $t\in\ic{\tinitial,\horizon}$
  is discrete and runs among a finite set of consecutive natural numbers;
\item
  $(\espacea{\OptState}_{\tinitial},\tribu{\OptState}_{\tinitial})$ (nature),
  $(\espacea{\OptState}_{\tinitial{+}1},\tribu{\OptState}_{\tinitial{+}1})$,
  \ldots, $(\espacea{\OptState}_{\horizon},\tribu{\OptState}_{\horizon})$
  (state spaces) are measurable spaces;
\item
  $(\espacea{\Control}_{\tinitial},\tribu{\Control}_{\tinitial})$,\ldots,
  $(\espacea{\Control}_{\horizon{-}1},\tribu{\Control}_{\horizon{-}1})$
  are measurable spaces (control spaces);
\item
  $\tribu{\Information}_t$ is a subfield of $\tribu{\OptState}_t$,
  for $t\in\ic{{\tinitial},\horizon{-}1}$ (information);
  \label{gdp:it:standard_form_information}
\item
  $\dynamics_t :(\espacea{\OptState}_{t} \times \espacea{\Control}_t, \tribu{\OptState}_{t} \otimes \tribu{\Control}_t)
  \rightarrow (\espacea{\OptState}_{t{+}1},\tribu{\OptState}_{t{+}1})$
  is measurable, for $t\in\ic{{\tinitial},\horizon{-}1}$ (dynamics);
\item
  $\pi_{\tinitial}$ is a probability
  on~$(\espacea{\OptState}_{\tinitial},\tribu{\OptState}_{\tinitial})$;
\item
  $\criterion :
  (\espacea{\OptState}_{\horizon},\tribu{\OptState}_{\horizon}) \rightarrow
  \RR$ is a measurable function (criterion).
\end{enumerate}
With these ingredients, Witsenhausen formulates a stochastic optimization
problem, whose solutions are to be searched among adapted feedbacks, namely
$ \policy_{t} : (\espacea{\OptState}_{t},\tribu{\OptState}_{t})
\rightarrow (\espacea{\Control}_{t},\tribu{\Control}_{t}) $ with the property
that $ \policy_{t}^{-1}(\tribu{\Control}_{t}) \subset
\tribu{\Information}_{t} $ for all $t\in\ic{{\tinitial},\horizon{-}1}$.
Then, he establishes a dynamic programming equation,
where the Bellman functions are function of the
(unconditional) distribution of the original state~$\optstate_t \in \espacea{\OptState}_t$,
and where the minimization is done over adapted feedbacks.
The main objective of Witsenhausen is to establish a
dynamic programming equation for nonclassical information patterns.

\subsubsubsection{Evstigneev approach}

The ingredients of the approach developed in~\citeAM{Evstigneev:1976}
are the following:
\begin{enumerate}
\item
  time $t\in\ic{\tinitial,\horizon}$
  is discrete and runs among a finite set of consecutive natural numbers;
\item
  $(\espacea{\Control}_{\tinitial},\tribu{\Control}_{\tinitial})$,\ldots,
  $(\espacea{\Control}_{\horizon{-}1},\tribu{\Control}_{\horizon{-}1})$ are measurable
  spaces (control spaces);
\item
  $ \np{\Omega,\tribu{\NatureField}} $ is a measurable space (nature);
\item
  $ \bseqa{ \tribu{\NatureField}_{t} }{ t\in\ic{{\tinitial},\horizon{-}1} } $
  is a filtration of~$\tribu{\NatureField}$ (information);
\item
  $ \PP $ is a probability on $ \np{\Omega,\tribu{\NatureField}} $;
\item
  $ \criterion : \np{
    \produit{\Omega}{\prod_{t\in\ic{{\tinitial},\horizon{-}1}} \espacea{\Control}_{t}} ,
    \oproduit{\tribu{\NatureField}}{\bigotimes_{t\in\ic{{\tinitial},\horizon{-}1}} \tribu{\Control}_{t}}
  } \to \RR $
  is a measurable function (criterion).
\end{enumerate}
With these ingredients, Evstigneev formulates a stochastic optimization
problem, whose solutions are to be searched among adapted processes,
namely random processes with values in
$\prod_{t\in\ic{{\tinitial},\horizon{-}1}} \espacea{\Control}_{t} $
and adapted to the filtration
$\bseqa{ \tribu{\NatureField}_{t} }{{t\in\ic{{\tinitial},\horizon{-}1}}}$.
Then, he establishes a dynamic programming equation,
where the Bellman function at time~$t$ is an $\tribu{\NatureField}_{t}$-integrand
depending on controls up to time~$t$ (random variables)
and where the minimization is done over $\tribu{\NatureField}_{t}$-measurable
random variables at time~$t$.
The main objective of Evstigneev is to establish an existence theorem
for an optimal adapted process (under proper technical assumptions,
especially on the objective function~$\criterion$, that we do not detail here).
Notice that there is no notion of state variable.

\subsubsubsection{Puterman approach}

The ingredients of the approach developed in~\citeAM[Sect.~2.1]{Puterman:1994}
are the following:
\begin{enumerate}
\item
  time $t\in\ic{\tinitial,\horizon}$
  is discrete and runs among a finite set of consecutive natural numbers;
\item
  $(\espacea{\OptState}_{\tinitial},\tribu{\OptState}_{\tinitial})$, \ldots,
  $(\espacea{\OptState}_{\horizon},\tribu{\OptState}_{\horizon})$ are measurable
  spaces (state spaces);
\item
  $(\espacea{\Control}_{\tinitial},\tribu{\Control}_{\tinitial})$,\ldots,
  $(\espacea{\Control}_{\horizon{-}1},\tribu{\Control}_{\horizon{-}1})$ are measurable
  spaces (control spaces);
\item
  $ \nu_{t:t+1} :  \espacea{\OptState}_{t} \times
  \espacea{\Control}_{t} \to \Delta\np{\espacea{\OptState}_{t{+}1}} $
  is a stochastic kernel, for $t\in\ic{{\tinitial},\horizon{-}1}$ (transitions);
\item
  $ \coutint_{\tter} : \STATE_{\tter} \times \CONTROL_{\tter}
  \to \RR $, for~$t\in\ic{{\tinitial},\horizon{-}1}$, and
  $ \coutfin : \STATE_{\horizon} \to \RR $, are measurable functions
  (instantaneous and final costs).
  % Voir p 190 Def 8.2
\end{enumerate}
With these ingredients, Puterman formulates a stochastic optimization
problem with a time additive cost function over given state
and control spaces, whose solutions are to be searched among history
feedbacks, namely sequences of mappings
$ \STATE_{\tinitial} \times \prod_{\tbis=\tinitial}^{t-1}
\np{ \CONTROL_{\tbis} \times \STATE_{\tbis{+}1} } \to \CONTROL_t $.
Then, he establishes a dynamic programming equation,
where the Bellman functions are function of the history
$ \history_t \in \STATE_{\tinitial} \times \prod_{\tbis=\tinitial}^{t-1}
\np{ \CONTROL_{\tbis} \times \STATE_{\tbis{+}1} } $.
He identifies cases where no loss of optimality results from reducing the search
to Markovian feedbacks $ \STATE_t \to \CONTROL_t $. In such cases,
the Bellman functions are function of the state~$\optstate_t \in \espacea{\OptState}_t$,
and the minimization in the dynamic programming equation is done over controls
$ \control_t \in \espacea{\Control}_t $.
The main objective of Puterman is to explore infinite horizon criteria,
average reward criteria, the continuous time case, and to present many examples.

\subsubsubsection{Hern\'andez-Lerma and Lasserre approach}

The ingredients of the approach developed
in~\citeAM[\S2.2, \S3.2, \S3.3]{Hernandez-Lerma-Lasserre:1996}
are the following:
\begin{enumerate}
\item
  time $t\in\ic{\tinitial,\horizon}$
  is discrete and runs among a finite set of consecutive natural numbers;
\item
  $(\espacea{\OptState}_{\tinitial},\tribu{\OptState}_{\tinitial})$, \ldots,
  $(\espacea{\OptState}_{\horizon},\tribu{\OptState}_{\horizon})$ are Borel %measurable
  spaces (state spaces);
\item
  $(\espacea{\Control}_{\tinitial},\tribu{\Control}_{\tinitial})$,\ldots,
  $(\espacea{\Control}_{\horizon{-}1},\tribu{\Control}_{\horizon{-}1})$
  are Borel spaces (control spaces); there are also feasible state-dependent
  control constraints that we do not present here;
\item
  $ \nu_{t:t+1} :  \espacea{\OptState}_{t} \times
  \espacea{\Control}_{t} \to \Delta\np{\espacea{\OptState}_{t{+}1}} $,
  for $t\in\ic{{\tinitial},\horizon{-}1}$,
  are Borel-measurable stochastic kernels (transitions);
\item
  $ \coutint_{\tter} : \STATE_{\tter} \times \CONTROL_{\tter}
  \to \RR $, for~$t\in\ic{{\tinitial},\horizon{-}1}$, and
  $ \coutfin : \STATE_{\horizon} \to \RR $ are Borel-measurable functions
  (instantaneous and final costs).
\end{enumerate}
With these ingredients, Hern\'andez-Lerma and Lasserre
formulate a stochastic optimization problem with a time additive cost
function over given state and control spaces.
They introduce the ``canonical construction'' %p.15
where the history at time~$t$ consists
in the states and the controls prior to~$t$. Then, they study optimization
problems whose solutions (policies) are to be searched among history feedbacks
(or randomized history feedbacks), namely sequences of mappings
$ \STATE_{\tinitial} \times \prod_{\tbis=\tinitial}^{t-1}
\np{ \CONTROL_{\tbis} \times \STATE_{\tbis{+}1} } \to \CONTROL_t $.
They identify cases where no loss of optimality results from reducing the search
to (relaxed) Markovian feedbacks $ \STATE_t \to \CONTROL_t $.
Then, they establish a dynamic programming equation,
where the Bellman functions are function of the state~$\optstate_t \in \espacea{\OptState}_t$,
and where the minimization is done over controls
$ \control_t \in \espacea{\Control}_t $.
For finite horizon problems, the mathematical challenge is
to set up a mathematical framework --- the Borel assumptions
plus additional topological ones presented
in~\citeAM[\S3.3]{Hernandez-Lerma-Lasserre:1996} ---
for which optimal policies exists.
The main objective of~\citeAM{Hernandez-Lerma-Lasserre:1996} is to offer
a unified and comprehensive treatment of discrete-time Markov control
processes, with emphasis on the case of Borel state and control spaces,
and possibly unbounded costs and noncompact control constraint sets.

\subsubsubsection{Bertsekas and Shreve approach}

The ingredients of the approach developed
in~\citeAM{Bertsekas-Shreve:1996} (more precisely in~\citeAM[Definition~10.1]{Bertsekas-Shreve:1996})
are the following:
\begin{enumerate}
\item
  time $t\in\ic{\tinitial,\horizon}$
  is discrete and runs among a finite set of consecutive natural numbers;
\item
  $(\espacea{\OptState}_{\tinitial},\tribu{\OptState}_{\tinitial})$, \ldots,
  $(\espacea{\OptState}_{\horizon},\tribu{\OptState}_{\horizon})$ are Borel %measurable
  spaces (state spaces);
\item
  $(\espacea{\Control}_{\tinitial},\tribu{\Control}_{\tinitial})$,\ldots,
  $(\espacea{\Control}_{\horizon{-}1},\tribu{\Control}_{\horizon{-}1})$ are Borel spaces (control spaces); there are also feasible state-dependent control
  constraints that we do not present here;
\item
  $(\espacea{\Uncertain}_{\tinitial},\tribu{\Uncertain}_{\tinitial})$,\ldots,
  $(\espacea{\Uncertain}_{\horizon},\tribu{\Uncertain}_{\horizon})$ are Borel spaces (noise);
\item
  $\dynamics_t :(\espacea{\OptState}_{t} \times \espacea{\Control}_t \times \espacea{\Uncertain}_t,
  \tribu{\OptState}_{t} \otimes \tribu{\Control}_t \otimes \tribu{\Uncertain}_t)
  \rightarrow (\espacea{\OptState}_{t{+}1},\tribu{\OptState}_{t{+}1})$,
  for $t\in\ic{{\tinitial},\horizon{-}1}$,
  are Borel-measurable mappings (dynamics);
\item
  $ \StochKernel{t}{t+1} :  \espacea{\OptState}_{t} \times
  \espacea{\Control}_{t} \to \Delta\np{\espacea{\Uncertain}_{t{+}1}} $,
  for $t\in\ic{{\tinitial},\horizon{-}1}$,
  are Borel-measurable stochastic kernels (noise distributions);
\item
  $ \coutint_{\tter} : \STATE_{\tter} \times \CONTROL_{\tter}
  \to \RR $, for~$t\in\ic{{\tinitial},\horizon{-}1}$, and
  $ \coutfin : \STATE_{\horizon} \to \RR $ are lower semianalytic
  functions (instantaneous and final costs).
\end{enumerate}
With these ingredients, Bertsekas and Shreve formulate a stochastic
optimization problem with a time additive cost function over given
state spaces, control spaces and uncertainty spaces.
They introduce the notion of history at time~$t$ which consists in the
states and the controls prior to~$t$ and study optimization problems
whose solutions (policies) are to be searched among history feedbacks
(or relaxed history feedbacks), namely sequences of mappings from history
space~$ \STATE_{\tinitial} \times \prod_{\tbis=\tinitial}^{t-1}
\np{ \CONTROL_{\tbis} \times \STATE_{\tbis{+}1} } \to \CONTROL_t $.
They identify cases where no loss of optimality results from reducing
the search to (relaxed) Markovian feedbacks $ \STATE_t \to \CONTROL_t $.
Then, they establish a dynamic programming equation,
where the Bellman functions are function of the state~$\optstate_t \in \espacea{\OptState}_t$,
and where the minimization is done over controls
$ \control_t \in \espacea{\Control}_t $.
For finite horizon problems, the mathematical challenge is
to set up a mathematical framework (the Borel assumptions)
for which optimal policies exists.
The main objective of Bertsekas and Shreve is to state conditions
under which the dynamic programming equation is mathematically sound
in the context of Borel spaces. The interested reader will find all
the subtleties in~\citeAM[Chapter 7]{Bertsekas-Shreve:1996}.

\subsubsubsection{Y\"{u}ksel approach}

As said at the beginning, the most general stochastic dynamic
programming principle is sketched by
Witsenhausen at the end of~\citeAM{Witsenhausen:1975b}.
This approach builds upon the so-called Witsenhausen intrinsic model
\citeAM{Witsenhausen:1975} which does not consider state, but information under
the form of $\sigma$-fields (see \citeAM{Witsenhausen:1988} for the functional
form).
In~\citeAM{Witsenhausen:1973}, Witsenhausen provides conditions to express
stochastic control optimization problems --- with information constraints, but
without state --- in standard form with a state (the first approach that we have
considered above).

Although Witsenhausen established a dynamic programming equation
in~\citeAM{Witsenhausen:1973}, Y\"{u}ksel notes in~\citeAM{Yuksel:2020}
that ``Witsenhausen's construction [\ldots] does not address the well-posedness of
such a dynamic program'' and that ``the existence problem was not considered''.
In the spirit of~\citeAM{Witsenhausen:1973}, Y\"{u}ksel entails in~\citeAM{Yuksel:2020}
``a general approach establishing that any sequential team optimization may admit
a formulation appropriate for a dynamic programming analysis''.
One of the contributions of~\citeAM{Yuksel:2020} is to propose a construction
of standard Borel controlled state and action spaces and to establish
a universal dynamic program for stochastic control optimization problems ---
with information constraints, but without state --- thus addressing some
of the issues raised and left open by Witsenhausen.
The ingredients are the following:
\begin{enumerate}
\item
  time $t\in\ic{\tinitial,\horizon}$
  is discrete and runs among a finite set of consecutive natural numbers;
\item
  $ \np{\Omega,\tribu{\NatureField}} $ is a measurable space (nature);
\item
  $(\espacea{\Control}_{\tinitial},\tribu{\Control}_{\tinitial})$,\ldots, $(\espacea{\Control}_{\horizon{-}1},\tribu{\Control}_{\horizon{-}1})$ are measurable
  spaces (control spaces);
\item
  $(\espacea{\Observation}_{\tinitial},\tribu{\Observation}_{\tinitial})$,
  \ldots, $(\espacea{\Observation}_{\horizon{-}1},\tribu{\Observation}_{\horizon{-}1})$ are measurable
  spaces (``observation'' spaces);
\item
  $\big\{ \eta_{t} : (
  \produit{\Omega}{\prod_{s\in\ic{{\tinitial},t}} \espacea{\Control}_{s}}$
  , $\oproduit{\tribu{\NatureField}}{\bigotimes_{s\in\ic{{\tinitial},t}} \tribu{\Control}_{s}}
  ) \to$
  \\ $(\espacea{\Control}_{t},\tribu{\Control}_{t}) \big\}_{t\in\ic{{\tinitial},\horizon{-}1}}$
  are measurable mappings (``measurement constraints'');
\item
  $ \PP $ is a probability on $ \np{\Omega,\tribu{\NatureField}} $;
\item
  $ \criterion : \np{
    \produit{\Omega}{\prod_{t\in\ic{{\tinitial},\horizon{-}1}} \espacea{\Control}_{t}} ,
    \oproduit{\tribu{\NatureField}}{\bigotimes_{t\in\ic{{\tinitial},\horizon{-}1}} \tribu{\Control}_{t}}
  } \to \RR_+ $
  is a measurable function (criterion).
\end{enumerate}
With these ingredients, Y\"{u}ksel
formulates a stochastic team optimization problem
whose solutions (policies) are to be searched among sequences of measurable
mappings (``design constraints'')
$ \espacea{\Observation}_{t{-}1} \to \CONTROL_t $,
and their ``randomized'' versions (so-called strategic measures).
He establishes a dynamic programming equation,
where the Bellman functions are function of probability distributions
and where the minimization is done over proper design mappings.
One objective of Y\"{u}ksel is
to set up a mathematical framework under which
the dynamic programming equation is mathematically sound \citeAM[Theorem~3.6]{Yuksel:2020}.

\subsubsubsection{Our approach}

The ingredients that we use (in Sect.~\ref{gdp:Stochastic_Dynamic_Programming_and_State_Reduction_by_Time_Blocks}
and in Sect.~\ref{gdp:State_Reduction_by_Time_Blocks}) are the following:
\begin{enumerate}
\item
  time $t\in\ic{\tinitial,\horizon}$
  is discrete and runs among a finite set of consecutive natural numbers;
\item
  $(\espacea{\Control}_{\tinitial},\tribu{\Control}_{\tinitial})$,\ldots,
  $(\espacea{\Control}_{\horizon{-}1},\tribu{\Control}_{\horizon{-}1})$
  are Borel   spaces (control spaces);
\item
  $(\espacea{\Uncertain}_{\tinitial},\tribu{\Uncertain}_{\tinitial})$,\ldots,
  $(\espacea{\Uncertain}_{\horizon},\tribu{\Uncertain}_{\horizon})$
  are Borel spaces (noise);
\item
  $ \StochKernel{t}{t+1} :
  \UNCERTAIN_{\tinitial} \times \prod_{\tbis=\tinitial}^{t-1}
  \np{ \CONTROL_{\tbis} \times \UNCERTAIN_{\tbis{+}1} }
  \to \Delta\np{\UNCERTAIN_{t{+}1}} $, for $t\in\ic{{\tinitial},\horizon{-}1}$,
  are Borel-measurable stochastic kernels (noise distributions);
\item
  $ \criterion: \np{
    \UNCERTAIN_{\tinitial}{\times}\prod_{\tbis=\tinitial}^{\horizon{-}1}
    \np{ \CONTROL_{\tbis}{\times}\UNCERTAIN_{\tbis{+}1} },
    \tribu{\Uncertain}_{\tinitial}{\otimes}\bigotimes_{\tbis=\tinitial}^{\horizon{-}1}
    \np{ \tribu{\Control}_{\tbis}{\otimes}\tribu{\Uncertain}_{\tbis{+}1}} }$
  $\to \bRR $
  is a nonnegative lower semianalytic function (criterion);
\item
  $ t_{0}< \cdots <t_{N} = \horizon $
  are the indices of multiple consecutive time blocks
  $ \ic{t_{0},t_{1}} $, \ldots, $ \ic{t_{N{-}1},t_{N}} $, with $N\geq 1$ a natural number;
\item
  $ \bseqa{ \np{ \STATE_{t_{j}},\tribu{\OptState}_{t_{j}}} }{j\in\ic{0,N}}  $
  are Borel spaces (time block state spaces);
\item
  $ \Bseqa{ \theta_{t_{j}} :  \UNCERTAIN_{\tinitial} \times \prod\limits_{\tbis=\tinitial}^{t_{j}-1}
    \np{ \CONTROL_{\tbis} \times \UNCERTAIN_{\tbis{+}1} } \to \STATE_{t_{j}} }{j\in\ic{1,N}} $
  and $ \theta_{t_{0}} :  \UNCERTAIN_{\tinitial} \to \STATE_{\tinitial} $
  are Borel-measurable mappings (time block reduction of history towards state);
\item
  $ \Bseqa{ \Dynamics{t_{j}}{t_{j{+}1}} :
    \STATE_{t_{j}}{\times}\prod\limits_{\tbis=t_{j}}^{t_{j{+}1}-1}
    \np{ \CONTROL_{\tbis}{\times}\UNCERTAIN_{\tbis{+}1} }
    \to \STATE_{t_{j{+}1}} }{j\in\ic{0,N{-}1}} $
  are Borel-measurable mappings (time block dynamics).
\end{enumerate}
The framework developed in the paper~\citeAM{tts:carpentier:JOCA} is intermediate between the ones of
Evstigneev in~\citeAM{Evstigneev:1976} and of Y\"{u}ksel in~\citeAM{Yuksel:2020}
--- notable by the absence of a state space ---
and the ones of
Witsenhausen~\citeAM{Witsenhausen:1973},
Hern\'andez-Lerma and Lasserre~\citeAM{Hernandez-Lerma-Lasserre:1996},
Bertsekas and Shreve~\citeAM{Bertsekas-Shreve:1996}
and Puterman~\citeAM{Puterman:1994} ---
where the state spaces are given for all times.

This said, our preoccupation could be adapted to any of the above frameworks.
Indeed, our objective is to establish a dynamic programming equation with a state,
not at any time~$t\in\ic{{\tinitial},\horizon}$, but at some specified instants
$ \tinitial <t_{1} < \cdots <t_{N} = \horizon $.
In~\S\ref{State_reduction_on_multiple_consecutive_time_blocks_and_dynamic_programming_equations},
the state spaces are introduced as image sets (codomains) of what we call
\emph{(time block) history reduction mappings} (where history at time~$t$ consists
of all uncertainties and controls prior to time~$t$).

\section{Supplement to Sect.~\ref{Two_time_scale_optimization_problem}}
\label{Link_Classical_Framework}

We make the link between the setting of two-time-scale optimization problems (as
developed in Sect.~\ref{Two_time_scale_optimization_problem} with stochastic
kernels) and the framework of stochastic optimal control (with random
variables).

The property that the stochastic kernels~\eqref{eq:natural_kernel} do not depend
on any decision variable makes it possible to build a probability
$\rho_{(\bottom\tslow,\bottom\tfast): (\topp\tslow,\topp\tfast)}$ on the product
space $\UNCERTAIN_{(\bottom\tslow,\bottom\tfast): (\topp\tslow,\topp\tfast)}$
by
\begin{multline}
  \rho_{(\bottom\tslow,\bottom\tfast): (\topp\tslow,\topp\tfast)} =
  \bigotimes_{\tslow\in \overline{\SlowSpan}}
  \Big(
  \rho^{\slow}_{\tslow:\nex\tslow} \np{\dd \price[\nex\tslow]} \otimes
  \rho^{\fast}_{(\nex\tslow,{\bottom\tfast}):(\nex\tslow,\nex{\bottom\tfast})}
  \kernelargs{\price[\nex\tslow]}{\dd\cost[(\nex\tslow,\nex{\bottom\tfast})] }
  \otimes \cdots \\
  \otimes
  \rho^{\fast}_{(\nex\tslow,\pre{\topp\tfast}):(\nex\tslow,{\topp\tfast})}
  \kernelargs{\price[\nex\tslow],\cost[(\nex\tslow,\nex{\bottom\tfast})],
              \cdots,\cost[(\nex\tslow,\pre{\topp\tfast})]}
             {\dd\cost[(\nex\tslow,{\topp\tfast})]}
  \Big) \eqfinp
\label{eq:probability}
\end{multline}

Then, with the notations given in
\S\ref{Structure_of_a_two_time_scale_optimization_problem} and using the
probability definied in Equation~\eqref{eq:probability},
Problem~\eqref{eq:Vzero}, may be rewritten as
\begin{subequations}
  \label{eq:pb-two-time-scale}
  \begin{align}
    \tilValue[\pre{\bottom\tslow}]
    \np{\stock[\pre{\bottom\tslow}]}
    & = \inf_{\gamma}
      \int_{\UNCERTAIN_{(\bottom\tslow,\bottom\tfast): (\topp\tslow,\topp\tfast)}}
      \Big( \instantcost[\topp\tslow]\bp{\stock[\topp\tslow]}
      \nonumber \\
    & \hspace{2cm}
      + \dsum[\tslow\in\SlowSpan] \instantcost[\tslow]
      \bp{\stock[\pre\tslow],\conso[\pre\tslow],\price[\tslow],
      \nseqa{\buffer[(\tslow,\tfast)],\buy[(\tslow,\tfast)],
      \cost[\nex{(\tslow,\tfast)}]}
      {\tfast\in \FastSpan\backslash\na{\topp\tfast}}}
      \Big) \nonumber \\
    & \hspace{4.5cm}
      \rho_{(\bottom\tslow,\bottom\tfast):(\topp\tslow,\topp\tfast)}
      \bp{\dd \price[\bottom\tslow], \dd \cost[(\bottom\tslow,\nex{\bottom\tfast})]
      \cdots \dd \cost[(\topp\tslow,\pre{\topp\tfast})],
      \dd \cost[(\topp\tslow,\topp\tfast)]}
      \label{eq:cost-two-time-scale} \\
    & s.t. \quad \vabuffer[ \nex{(\tslow,\tfast)}] =
      \hatdyn_{(\tslow,\tfast)}^{\fast}(\vabuffer[(\tslow,\tfast)],
      \vabuy[(\tslow,\tfast)], \vacost[\nex{(\tslow,\tfast)}])
      \eqsepv
      \forall \tslow \in \SlowSpan\eqsepv
      \forall \tfast \in \FastSpan\setminus\na{\topp\tfast}
      \eqfinv \label{eq:dynamics-1} \\
    & \hphantom{s.t. \quad}
      \vabuffer[(\nex\tslow,\bottom\tfast)] =
      \hatdyn_{\tslow}^{\slow}(\vastock[\tslow], \vaconso[\tslow], \vaprice[\nex\tslow])
      \eqsepv \forall \tslow \in \overline{\SlowSpan}\setminus\na{\topp\slow}
      \eqfinv \label{eq:dynamics-2} \\
    & \hphantom{s.t. \quad}
      \vaconso[\tslow] = \gamma_{\tslow}
      \bp{\nseqa{\control_{(\tslow',\tfast')},
      \uncertain_{\nex{(\tslow',\tfast')}}}
      {(\tslow',\tfast') \prec(\tslow,\topp\tfast)}}
      \eqsepv \forall \tslow \in \overline{\SlowSpan}\setminus\na{\topp\slow}
      \eqsepv \label{eq:controls-1} \\
    & \hphantom{s.t. \quad}
      \vabuy[(\tslow,\tfast)]= \gamma_{(\tslow,\tfast)}
      \bp{\nseqa{\control_{(\tslow',\tfast')},\uncertain_{\nex{(\tslow',\tfast')}}}
      {(\tslow',\tfast') \prec(\tslow,\tfast)}} \eqsepv
      \forall \tslow \in \SlowSpan\eqsepv \forall \tfast
      \in \FastSpan\setminus\na{\topp\tfast} \eqfinp
      \label{eq:controls-2}
  \end{align}
\end{subequations}

The integral cost given in the right hand side
of Equation~\eqref{eq:cost-two-time-scale} can be reformulated as
a mathematical expectation, denoted by~$\EE$, with respect to the
probability~$\rho_{(\bottom\tslow,\bottom\tfast): (\topp\tslow,\topp\tfast)}$
by introducing random variables for the exogeneous noises as projection mappings from
$\UNCERTAIN_{(\bottom\tslow,\bottom\tfast): (\topp\tslow,\topp\tfast)}$ to $\UNCERTAIN_{(\tslow,\tfast)}$ for all
$(\tslow,\tfast)\in \SlowSpan{\times}\FastSpan$
\begin{equation}
  \va{W}_{(\tslow,\tfast)} : \UNCERTAIN_{(\bottom\tslow,\bottom\tfast): (\topp\tslow,\topp\tfast)} \to \UNCERTAIN_{(\tslow,\tfast)} \eqsepv
  \forall (\tslow,\tfast)\in \SlowSpan{\times}\FastSpan \eqfinv
\end{equation}
and obtaining random variables for the states and the control through
the dynamics equations~\eqref{eq:dynamics-1}--\eqref{eq:dynamics-2}
and the feedback equations~\eqref{eq:controls-1}--\eqref{eq:controls-2}.

This leads to a reformulation of Problem~\eqref{eq:pb-two-time-scale}
as a classical stochastic optimal control problem
\begin{subequations}
\begin{align}
\inf \; \EE
  & \Big[
    \dsum[\tslow\in\SlowSpan] \instantcost[\tslow]
    \bp{
    \va{\OptState}^{\slow}_{\pre\tslow}, \va{\Control}^{\slow}_{\pre\tslow}, \va{\Uncertain}_{\tslow},
    \nseqa{ \va{\OptState}^{\fast}_{(\tslow,\tfast)},\va{\Control}^{\fast}_{(\tslow,\tfast)},
    \va{\Uncertain}^{\tfast}_{\nex{(\tslow,\tfast)}}}
    {\tfast\in \FastSpan\backslash\na{\topp\tfast}}} +
    \instantcost[\topp\tslow]\bp{\va{\OptState}^{\slow}_{\topp\tslow}}
    \Big] \\
  & s.t. \quad
    \va{\OptState}^{\fast}_{ \nex{(\tslow,\tfast)}} =
    \hatdyn_{(\tslow,\tfast)}^{\fast}(\va{\OptState}^{\fast}_{(\tslow,\tfast)},
    \va{\Control}^{\fast}_{(\tslow,\tfast)},
    \va{\Uncertain}^{\fast}_{\nex{(\tslow,\tfast)}})
    \eqsepv \forall \tslow \in \SlowSpan\eqsepv
    \forall \tfast \in \FastSpan\setminus\na{\topp\tfast}
    \eqfinv \label{eq:dynamics-1-va}\\
  & \hphantom{s.t. \quad} \va{\OptState}^{\fast}_{(\nex\tslow,\bottom\tfast)}=
    \hatdyn_{\tslow}^{\slow}
    \bp{\va{\OptState}^{\slow}_{\tslow}, \va{\Control}^{\slow}_{\tslow},
    \va{\Uncertain}^{\slow}_{\nex\tslow}}
    \eqsepv \forall \tslow \in \overline{\SlowSpan}\backslash\na{\topp\slow}
    \eqfinv \label{eq:dynamics-2-va} \\
  & \hphantom{s.t. \quad}
    \va{\Control}^{\slow}_{\tslow} \in \CONSO[\tslow]
    \eqsepv \forall \tslow \in \overline{\SlowSpan}\backslash\na{\topp\slow} \eqsepv \\
  & \hphantom{s.t. \quad}
    \sigma(\va{\Control}^{\slow}_{\tslow})
    \subset \sigma(\nseqa{\va{\Uncertain}^{\slow}_{\tslow}}{\tslow' \preceq \tslow },
    \nseqa{\va{\Uncertain}^{\fast}_{(\tslow',\tfast')}}
    {(\tslow',\tfast') \preceq (\tslow,\topp\tfast)})
    \eqsepv \forall \tslow \in \overline{\SlowSpan}\backslash\na{\topp\slow}
    \label{eq:pb-two-time-scale-va-mes-slow} \eqsepv \\
  & \hphantom{s.t. \quad}
    \va{\Control}^{\fast}_{(\tslow,\tfast)} \in \BUY[(\tslow,\tfast)]
    \eqsepv \forall \tslow \in \SlowSpan\eqsepv
    \forall \tfast \in \FastSpan\setminus\na{\topp\tfast} \eqsepv \\
  & \hphantom{s.t. \quad}
    \sigma\bp{\va{\Control}^{\fast}_{(\tslow,\tfast)}}
    \subset
    \sigma\bp{\nseqa{\va{\Uncertain}^{\slow}_{\tslow'}}{\tslow' \preceq \tslow },
    \nseqa{\va{\Uncertain}^{\fast}_{(\tslow',\tfast')}}
    {(\tslow',\tfast') \preceq (\tslow,\tfast)}} \eqsepv \forall \tslow \in \SlowSpan\eqsepv
    \forall \tfast \in \FastSpan\setminus\na{\topp\tfast}
    \label{eq:pb-two-time-scale-va-mes-fast} \eqfinv
\end{align}
\label{eq:pb-two-time-scale-va}
\end{subequations}
where the two feedback constraints in~\eqref{eq:controls-1}
and~\eqref{eq:controls-2} are reformulated as measurability
constraints~\eqref{eq:pb-two-time-scale-va-mes-slow} and~\eqref{eq:pb-two-time-scale-va-mes-fast},
using the $\sigma$-fields generated by random variables
(of course, a formal equivalence would require to be more
specific about spaces to use Doob functional Lemma).

\section{Detailed proof of
  Proposition~\ref{pr:DHD_multistage_stochastic_optimization_problem}
%   Technical details and proofs
  (Sect.\ref{gdp:Decision_Hazard_Decision_Dynamic_Programming})}
\label{Technical_details_and_proofs_DHD}

As indicated in the sketchy proof of
Proposition~\ref{pr:DHD_multistage_stochastic_optimization_problem},
we give here a detailed proof of this latter.
For this purpose, we first flesh out the ingredients necessary to formulate a multistage stochastic
optimization problem with decision-hazard-decision information structure.
Second, we provide in Theorem~\ref{cor:instantaneous-costs-in-addition-to-final-cost_DHD} a
Bellman equation for such a multistage stochastic optimization problem (see
Equation~\eqref{gdp:eq:BelOpnatural-expression}). Third, we give the proof of
Theorem~\ref{cor:instantaneous-costs-in-addition-to-final-cost_DHD} which is
derived through an embeding of the decision-hazard-decision structure as a
particular case of multiple consecutive time blocks followed by an application
of Theorem~\ref{cor:instantaneous-costs-in-addition-to-final-cost}.

\subsubsubsection{History}

Let $ \horizonDHD \in \NN^* $.  For each time~$\tDHD\in\ic{0,\horizonDHD{-}1}$,
the \emph{head decision}~$ \control_{\tDHD}^\Before$ takes values in a Borel space~$\CONTROL^\Before_{\tDHD}$.
For each time~$\tDHD \in \ic{1,\horizonDHD}$, the \emph{tail decision}~$\control_{\tDHD}^\After$ takes values in
a Borel space~$\CONTROL^\After_{\tDHD}$.  For each
time~$\tDHD\in\ic{1,\horizonDHD}$, the uncertainty~$\uncertain^\After_{\tDHD}$
takes its values in a Borel space~$\UNCERTAIN^\After_{\tDHD}$.
For time~$\tDHD=0$, the
uncertainty~$\uncertain^\Before_{0}$ takes its values in a Borel space~$\UNCERTAIN^\Before_{0}$.

At the beginning of the time interval~$[\tDHD,\tDHD+1[$,
the decision-maker makes a \emph{head decision}~$\control_{\tDHD}^\Before$.
What is new --- in comparison with the classical decision-hazard framework ---
is that, at the end of the time interval~$[\tDHD,\tDHD+1[$,
when an uncertainty variable~$\uncertain^\After_{\tDHD+1}$ is revealed,
the decision-maker has the possibility to make a
\emph{tail decision}~$\control_{\tDHD+1}^\After$.
This latter decision~$\control_{\tDHD+1}^\After$ can be thought
as a \emph{recourse} variable for a two stage stochastic optimization problem
that would take place inside the time interval~$[\tDHD,\tDHD+1[$.
We call~$\uncertain^\Before_0$ the uncertainty happening right before the first decision.
The interplay between uncertainties and decisions is thus as follows
(compare the chronology with the one in~\eqref{gdp:eq:interplay_Noise_Control}):
\begin{multline*}
  \uncertain^\Before_{0} \rightsquigarrow \control^\Before_{0}
  \rightsquigarrow \uncertain^\After_{1} \rightsquigarrow \control^\After_{1}
  \rightsquigarrow \control^\Before_{1} \rightsquigarrow  \uncertain^\After_{2}
  \rightsquigarrow  \;\; \dots  \;\; \rightsquigarrow
  \uncertain^\After_{\horizonDHD{-}1} \rightsquigarrow \control^\After_{\horizonDHD{-}1} \rightsquigarrow \control^\Before_{\horizonDHD{-}1}
  \rightsquigarrow \uncertain^\After_{\horizonDHD} \rightsquigarrow \control^\After_{\horizonDHD} \eqfinp
\end{multline*}

\subsubsubsection{History Spaces}

For $\tDHD\in\ic{0,\horizonDHD}$, we define the \emph{head history space}
as the product Borel space
\begin{subequations}
  \begin{align}
    \HISTORY^\Before_{\tDHD}
    &=
      \UNCERTAIN_{0}^\Before \times \prod_{\tDHDIt=0}^{\tDHD-1}
      \bp{ \CONTROL^\Before_{\tDHDIt} \times \UNCERTAIN_{\tDHDIt+1}^\After
      \times \CONTROL^\After_{\tDHDIt+1} } \eqfinp
      \label{gdp:eq:head_history_space}
      \intertext{We also define,
      for $\tDHD\in\ic{1, \horizonDHD}$,
      the \emph{tail history space} as the product Borel space }
      \HISTORY^\After_{\tDHD}
    &=
      \HISTORY^\Before_{\tDHD-1}
      \times \CONTROL^\Before_{\tDHD-1}
      \times \UNCERTAIN_{\tDHD}^\After
      \eqfinp
      \label{gdp:eq:tail_history_space}
  \end{align}
\end{subequations}

\subsubsubsection{Stochastic kernels}

We introduce a family
$\na{\StochKernel{\tDHD-1}{\tDHD}^{\After}}_{\tDHD\in\ic{1,\horizonDHD}}$
of probability distributions (constant Borel-measurable stochastic kernels), with
\begin{equation}
  \StochKernel{\tDHD-1}{\tDHD}^{\After}
  \in  \Delta\np{\UNCERTAIN_{\tDHD}^\After}
  \eqfinp
  \label{gdp:eq:rhodhd}
\end{equation}

\subsubsubsection{History feedbacks}

\begin{subequations}
  For $\tDHD\in \ic{0,\horizonDHD{-}1}$,
  a \emph{head history feedback} at time~$\tDHD$
  is a universally measurable mapping
  \begin{equation}
    \HistoryFeedback^\Before_\tDHD:
    \HISTORY^\Before_{\tDHD} \to \CONTROL^\Before_\tDHD \eqfinp
  \end{equation}
  We call~$\HISTORYFEEDBACK^\Before_{\tDHD}$ the
  \emph{set of head history feedbacks at time~$\tDHD$},
  and we define
  $\HISTORYFEEDBACK^\Before_{\tDHD:\horizonDHD} =
  \HISTORYFEEDBACK^\Before_{\tDHD} \times \cdots
  \times \HISTORYFEEDBACK^\Before_{\horizonDHD}$.
  We also define, for all $\tDHD \in\ic{1,\horizonDHD}$,
  a \emph{tail history feedback} at time~$\tDHD$
  as a universally measurable mapping
  \begin{equation}
    \HistoryFeedback^\After_\tDHD:
    \HISTORY^\After_{\tDHD} \to \CONTROL^\After_\tDHD \eqfinp
  \end{equation}
  We call~$\HISTORYFEEDBACK^\After_{\tDHD}$
  the \emph{set of tail history feedbacks at time~$\tDHD$},
  and we define
  $\HISTORYFEEDBACK^\After_{\tDHD:\horizonDHD} =
  \HISTORYFEEDBACK^\After_{\tDHD} \times \cdots
  \times \HISTORYFEEDBACK^\After_{\horizonDHD}$.
  \label{eq:DHD_history_feedbacks}
\end{subequations}

\subsubsubsection{Value functions}

Let $\nseqa{\STATE_t^\Before}{t \in \ic{0,\horizon}}$ (states)
be sequences of Borel spaces with $\STATE_0^\Before=\UNCERTAIN_0^\Before$.
\begin{subequations}
  Let also be given Borel-measurable dynamics mappings
  \begin{equation}
    {\dynamics_{t} : \STATE_t^\Before \times \CONTROL^\Before_t \times\UNCERTAIN_{t{+}1}^\After\times \CONTROL^\After_{t{+}1} }
    \to \STATE_{t{+}1}^\Before
    \eqsepv \forall t \in \ic{0,\horizon{-}1}
    \eqfinv
    \label{eq:DHD_dynamics}
  \end{equation}
  nonnegative lower semianalytic instantaneous cost functions
  \begin{equation}
    {\coutint_{t} : \STATE_t^\Before \times \CONTROL^\Before_t \times \UNCERTAIN_{t{+}1}^\After \times \CONTROL^\After_{t{+}1} }
    \to \bRR
    \eqsepv \forall t \in \ic{0,\horizon{-}1}
    \eqfinv
  \end{equation}
  and a nonnegative lower semianalytic final cost function
  \begin{equation}
    \coutfin : \STATE_{\horizon}^\Before \to \bRR
    \eqfinp
  \end{equation}
  \label{eq:DHD_dynamics_cost_bis}
\end{subequations}

Second, we recursively define a family of Borel-measurable reduction
mappings $ \theta_{t}^\Before: \HISTORY_{t}^\Before \to \STATE_{t}^\Before$ by the following forward recursion
for $t\in \ic{0,\horizon{-}1}$
\begin{equation}
  \theta_{t+1}^\Before\np{\history_{t+1}^\Before}=
  \dynamics_{t} \bp{ \theta_{t}^\Before\np{\history_{t}^\Before},
    {\control_t^\Before,\uncertain_{t+1}^\After,\control_{t+1}^\After}}
  \eqfinv
  \text{ and } \theta_{0}^\Before \np{\history_{0}^\Before} = \history_{0}^\Before
  \label{eq:theta_natural}
  \eqfinv
\end{equation}
and consider a nonnegative and lower semianalytic numerical function
$\criterion :  \HISTORY^{\Before}_{\horizonDHD} \to \bRR$ defined by
\begin{equation}
  \criterion\np{ \history_{\horizonDHD}^\Before }
  =
  \sum_{t=0}^{\horizon{-}1}
  \coutint_{t} \bp{ \theta_{t}^\Before\np{\history_{t}}, {\control_t^\Before,\uncertain_{t+1}^\After,\control_{t+1}^\After}}
  +
  \coutfin \bp{\theta_{\horizon}^\Before\np{\history_{\horizon}^\Before}}
  \eqfinp
  \label{gdp:eq:criterion_time_block_additive_DHD}
\end{equation}

For~$\tDHD \in\ic{0,\horizonDHD}$, we define \emph{value functions} by
\begin{equation}
  \Value_{\tDHD}(\history^\Before_{\tDHD}) =
  \inf_{\HistoryFeedback^\Before \in \HISTORYFEEDBACK^\Before_{\tDHD:\horizonDHD{-}1},
    \HistoryFeedback^\After \in \HISTORYFEEDBACK^\Before_{\tDHD+1:\horizonDHD}}
  \int_{\HISTORY^\Before_{\horizonDHD}}
  \criterion\np{\history'_{\horizonDHD}}
  {\StochKernelFeed{\tDHD}{\horizonDHD}{\HistoryFeedback^\Before,
      \HistoryFeedback^\After}}
  \nkernelargs{\history^\Before_\tDHD}{\dd\history'_{\horizonDHD}}
  \eqsepv \forall \history^\Before_{\tDHD} \in \HISTORY^\Before_{\tDHD}
  \eqfinv
  \label{gdp:eq:value_function_dhd}
\end{equation}
where $ {\StochKernelFeed{\tDHD}{\horizonDHD}{\HistoryFeedback^\Before,
    \HistoryFeedback^\After}}$ has to
be understood as $ \StochKernelFeed{\tDHD}{\horizonDHD}{\HistoryFeedback} $
(see Definition~\ref{gdp:de:stochastic_kernels_rho}), with
\begin{subequations}
  \label{gdp:eq:feedback_dhd}
  \begin{align}
    \HistoryFeedback_{\tDHD}(\history_{\tDHD}^\Before)
    &=
      \HistoryFeedback^\Before_{\tDHD}(\history_{\tDHD}^\Before)
      \eqsepv
      \forall \history_{\tDHD}^\Before \in \HISTORY_{\tDHD}^\Before
      \eqfinv
    \\
    \HistoryFeedback_{\tDHDIt}(\history_{\tDHDIt}^\After)
    &=
      \Bp{\HistoryFeedback^\After_{\tDHDIt}\np{\history_{\tDHDIt}^\After},
      \HistoryFeedback^\Before_{\tDHDIt}\bp{\history_{\tDHDIt}^\After,
      \HistoryFeedback^\After_{\tDHDIt}\np{\history_{\tDHDIt}^\After} } }
      \eqsepv
      \forall \tDHDIt\in\ic{\tDHD+1,\horizonDHD{-}1}
      \eqsepv
      \forall \history_{\tDHDIt}^\After \in \HISTORY_{\tDHDIt}^\After \eqfinv
    \\
    \HistoryFeedback_{\horizonDHD}(\history_{\horizonDHD}^\After)
    &=
      \HistoryFeedback^\After_{\horizonDHD}(\history_{\horizonDHD}^\After)
      \eqsepv
      \forall \history_{\horizonDHD}^\After \in \HISTORY_{\horizonDHD}^\After
      \eqfinp
  \end{align}
\end{subequations}
In fact, in the special case~\eqref{gdp:eq:rhodhd} we are developing here, $ {\StochKernelFeed{\tDHD}{\horizonDHD}{\HistoryFeedback^\Before,
    \HistoryFeedback^\After}}$ does not depend on the history feedbacks but reduces to the product
\begin{equation}
   {\StochKernelFeed{\tDHD}{\horizonDHD}{\HistoryFeedback^\Before,
       \HistoryFeedback^\After}} =
   \StochKernel{\tDHD-1}{\tDHD}^{\After} \otimes \ldots \otimes  \StochKernel{\horizon-1}{\horizon}^{\After}
   \eqfinp
\end{equation}

\def\doubletime{\tau}

\begin{theorem}
  \label{cor:instantaneous-costs-in-addition-to-final-cost_DHD}
  We assume to be in the setting of~\S\ref{The_Bertsekas-Shreve_setting}.
  % Suppose that the assumptions of Theorem~\ref{gdp:thm:DPB_family} are
  % satisfied, but for the cost criterion $\criterion:\HISTORY_{\horizon}^\Before \to \bRR$
  % defined by Equation~\eqref{gdp:eq:criterion_time_block_additive_DHD}. %\eqref{gdp:eq:criterion_time_block_additive}.
  We define the sequence of \emph{reduced value functions}
  $\na{\Value_{t}^\Before}_{t\in \ic{0,\horizon}}$, where \( \Value_{{t}}^\Before :
  \STATE_{t}^\Before \to \bRR \) for $t \in \ic{0,\horizon}$, by
  \begin{equation}
    \Value_{\horizon}^\Before = \coutfin_{\horizon}
    \mtext{ and }
    \Value_{t}^\Before = \naturalBelOp{t{+}1}{t} \Value_{t{+}1}^\Before
    \eqsepv \forall t\in\ic{0,\horizon{-}1}
    \eqfinv
    \label{gdp:eq:reduced_value_functions_time_block_additive_dhd}
  \end{equation}
  where the reduced Bellman operator~$\naturalBelOp{{t{+}1}}{{t}}$ are given,
  for any $t\in\ic{0,\horizon{-}1}$, for any~$\varphi_{t{+}1}^\Before \in
  \espace{L}^{0}_{+}(\STATE_{{t{+}1}}^\Before)$
  and for any $\optstate_{t}^\Before\in\STATE_{t}^\Before$, by
  \begin{align}
    \big( \naturalBelOp{{t{+}1}}{{t}}\varphi_{{t{+}1}}^\Before\big)
    (\optstate_{{t}}^\Before)
    =
    \inf_{\control_{t}^\Before\in\CONTROL_{t}^\Before}
    & \int_{\UNCERTAIN_{t}^\After}
      \StochKernel{t-1}{t}^\After \np{\dd\uncertain_{t}^\After}
      \nonumber \\
    &
      \inf_{\control_{t+1}^\After\in\CONTROL_{t+1}^\After}
      \big(
      \coutint_{t}
      \np{\optstate_{{t}}^\Before, \control_{{t}}^\Before,\uncertain_{{t}{+}1}^\After,\control_{t{+}1}^\After}
      \nonumber  \\
    & \quad\quad\quad
      +
      \varphi_{t{+}1}^\Before\bp{ \dynamics_{t}
      \np{\optstate_{t}^\Before, \control_{t}^\Before,\uncertain_{{t}{+}1}^\After,\control_{t+1}^\After}
      }
      \big)
      \label{gdp:eq:BelOpnatural-expression}
      \eqfinp \!\!\!\!\!
  \end{align}
  Then, the value function  $\Value_{0}$ given in Equation~\eqref{gdp:eq:value_function_dhd}
  coincides with the value function $\Value_{{0}}^\Before$.
\end{theorem}

\begin{proof} We embed the decision-hazard-decision structure as a particular
  case of multiple consecutive time blocks (of size~2) as in~\S\ref{State_reduction_on_multiple_consecutive_time_blocks_and_dynamic_programming_equations}
  in order to use Theorem~\ref{cor:instantaneous-costs-in-addition-to-final-cost} to obtain
  the reduced Bellman equation~\eqref{gdp:eq:overlineBelOp-expression} which boils down to~Equation~\eqref{gdp:eq:BelOpnatural-expression}
  in the  decision-hazard-decision case.
  For each time~$\tDHD\in\ic{1,\horizonDHD}$,
  we introduce a spurious uncertainty variable~$\uncertain^\Before_\tDHD$
  taking values in a singleton set~$\UNCERTAIN^\Before_{\tDHD} =
  \{\bar\uncertain^\Before_{\tDHD}\}$, so that we obtain
  the following chronology
  \begin{multline*}
    \underbrace{\uncertain^\Before_{0}}_{\state^\Before_{0}}
    \rightsquigarrow \control^\Before_{0}
    \rightsquigarrow \uncertain^\After_{1}
    \rightsquigarrow \control^\After_{1}
    \rightsquigarrow \underbrace{\uncertain^\Before_{1}}_{\text{spurious}}
    \rightsquigarrow \control^\Before_{1}
    \rightsquigarrow \uncertain^\After_{2}
    \rightsquigarrow \control^\After_{2}
    \rightsquigarrow \underbrace{\uncertain^\Before_{2}}_{\text{spurious}}
    \rightsquigarrow \control^\Before_{2}
    \rightsquigarrow  \quad \dots \\ \quad
    \rightsquigarrow \uncertain^\After_{\horizonDHD{-}1}
    \rightsquigarrow \control^\After_{\horizonDHD{-}1}
    \rightsquigarrow \underbrace{\uncertain^\Before_{\horizonDHD{-}1}}_{\text{spurious}}
    \rightsquigarrow \control^\Before_{\horizonDHD{-}1}
    \rightsquigarrow \uncertain^\After_{\horizonDHD}
    \rightsquigarrow \control^\After_{\horizonDHD}
    \rightsquigarrow \underbrace{\uncertain^\Before_{\horizonDHD}}_{\text{spurious}}
    \eqfinp
  \end{multline*}
  Proceeding this way, we have doubled the timeline as time $t$ has been
  ``duplicated'' in the ordered pair $(t,\After)$ and $(t,\Before)$.  We
  introduce new notations to explicitely deal with the new duplicated
  timeline. An element of the duplicated timeline is denoted by $\doubletime$.
  For $\doubletime\in \ic{0,2\horizonDHD}$, we introduce the
  sets
  \begin{subequations}
    \label{eq:double_time_line_controls_noises}
  \begin{align}
    \CONTROL_{\doubletime}
    &=
      \begin{cases}
        \CONTROL_{s}^{\Before} & \text{if}~\doubletime = 2s 
        \\
        \CONTROL_{s}^{\After} & \text{if}~\doubletime = 2s+1 
      \end{cases}
    \eqsepv \forall \doubletime \in \ic{0,2\horizonDHD -1}
      \eqfinv
      \label{eq:double_time_line_controls}
    \\
    \UNCERTAIN_{\doubletime}
    &=
      \begin{cases}
        \UNCERTAIN_{s}^{\Before} & \text{if}~\doubletime = 2s
        \\
        \UNCERTAIN_{s}^{\After} & \text{if}~\doubletime = 2s +1
      \end{cases}
      \eqsepv \forall   \doubletime\in \ic{0,2\horizonDHD}
      \label{eq:double_time_line_noises}
    \eqfinv
  \end{align}
\end{subequations}

For $\doubletime\in\ic{0, 2 \horizon}$, we define the history space~$\HISTORY_{\doubletime}$
as the product Borel space
\begin{equation}
    \HISTORY_{\doubletime}
    = \UNCERTAIN_{0} \times \prod_{\tbis=1}^{\doubletime}
    \np{ \CONTROL_{\tbis{-}1} \times \UNCERTAIN_{\tbis} }
    \eqfinv
\end{equation}

Then, given the times $\doubletime_i=2i$ for $i\in \ic{0,\horizon}$, we consider the consecutive time blocks $ \ic{\doubletime_{i},\doubletime_{i{+}1}} $, 
for $i \in \ic{0,\horizon{-}1}$ whose union covers the doubled timeline.
For $i \in \ic{0,\horizon}$, we define the state spaces 
\(  \STATE_{\doubletime_i}= \STATE_{i}^\Before \).
Using Equation~\eqref{eq:DHD_dynamics}, we define a family of Borel-measurable
dynamics $\nseqa{ \Dynamics{\doubletime_{i}}{\doubletime_{i{+}1}}}{i \in \ic{0,\horizon{-}1}}$ by
\begin{align}
  \Dynamics{\doubletime_{i}}{\doubletime_{i{+}1}} :\hspace{1cm} \STATE_{\doubletime_{i}} \times \HISTORY_{\doubletime_{i}{+}1:\doubletime_{i{+}1}}
  & \to \STATE_{\doubletime_{i{+}1}}
    \nonumber
  \\
  \bp{\state^\Before_i,\underbrace{\control^\Before_i ,\uncertain_{i{+}1}^\After,
  \control^\After_{i{+}1}, \uncertain_{i{+}1}^\Before }_{\history_{\doubletime_{i}{+}1:\doubletime_{i{+}1}}}
  }
  &
    \mapsto
    \dynamics_{i}\np{\state_i^\Before, \control^\Before_i , \uncertain_{i{+}1}^\After, \control^\After_{i{+}1}}
    \eqfinp
    \label{eq:dynamics_dhd_proof}
\end{align}
We recursively define a family of Borel-measurable reduction
mappings $ \theta_{\doubletime_{i}} : \HISTORY_{t_{i}} \to \STATE_{t_{i}} $ by the following forward recursion
for $i\in \ic{0,\horizon{-}1}$
\begin{equation}
  \theta_{\doubletime_{i+1}}\np{\history_{\doubletime_{i+1}}}=
  \theta_{\doubletime_{i+1}}\bp{\np{\history_{\doubletime_i}, \history_{\doubletime_i{+}1:\doubletime_{i+1}} }} =
  \Dynamics{\doubletime_i}{\doubletime_{i+1}} \bp{ \theta_{\doubletime_i}\np{\history_{\doubletime_i}},\history_{\doubletime_i{+}1:\doubletime_{i+1}}}
  \eqfinv
  \text{ and } \theta_{\doubletime_{0}}\np{\history_{0}} = \history_{0}
  \eqfinv
\end{equation}
to obtain 
the family
$\big( \nseqa{ \STATE_{\doubletime_{i}} }{i\in \ic{0,\horizon}}$,
$\nseqa{\theta_{\doubletime_{i}}}{i\in \ic{0,\horizon}}$,
$\nseqa{ \Dynamics{\doubletime_{i}}{\doubletime_{i{+}1}} }{i \in \ic{0,\horizon{-}1}}\big)$
which gives a state reduction across the consecutive time blocks
$\ic{\doubletime_{i},\doubletime_{i{+}1}}$, $i \in \ic{0,\horizon{-}1}$.
It is worth noting that, for $i\in \ic{0,\horizon}$, we have that $\theta_{\doubletime_{i}}=\theta_{i}^\Before$ where the mapping $\theta_{i}^\Before$ is defined in Equation~\eqref{eq:theta_natural}.

Now, for each $i \in \ic{1,\horizon}$ we consider the family of stochastic kernels
$\na{\StochKernel{i-1}{i}^\After}_{i\in\ic{1,\horizon}}$ given by Equation~\eqref{gdp:eq:rhodhd} which are probability distributions
and the  family of spurious stochastic kernels $\na{\StochKernel{i-1}{i}^\Before}_{i\in\ic{1,\horizon}}$ which are
Dirac measures at fixed points $\na{\bar\uncertain_{i}^\Before}_{i\in\ic{1,\horizon}}$. The state reduction across the consecutive time blocks
$\ic{\doubletime_{i},\doubletime_{i{+}1}}$ is indeed compatible with the
family~$\na{\StochKernel{\doubletime-1}{\doubletime}}_{\doubletime\in\ic{0,2\horizon}}$ 
of stochastic kernels given by
\begin{align}
  \StochKernel{\doubletime-1}{\doubletime}
  &=
    \begin{cases}
      \StochKernel{\tbis-1}{\tbis}^{\Before}
      & \text{if}~\doubletime = 2s
      \\
      \StochKernel{\tbis-1}{\tbis}^{\After}
      & \text{if}~\doubletime = 2s +1
    \end{cases}
    \eqsepv \forall \doubletime \in \ic{0,2\horizonDHD}
    \eqfinp
    \label{eq:def_rho_dhd}
\end{align}
Finally we introduce the familly of cost functions $\na{\InstantneousCost_{\doubletime_i}}_{i\in \ic{0,\horizon}}$, given by
\begin{subequations}
  \label{eq:couts_dhd_proof}
\begin{align}
  \InstantneousCost_{\doubletime_i}
  &: \STATE_{\doubletime_{i}} \times \HISTORY_{\doubletime_{i}{+}1:\doubletime_{i{+}1}}
  \ni 
  \bp{\state^\Before_i,\control^\Before_i ,\uncertain_{i{+}1}^\After,
  \control^\After_{i{+}1}, \uncertain_{i{+}1}^\Before }
    \mapsto
    \coutint_{i}\np{\state_i^\Before, \control^\Before_i , \uncertain_{i{+}1}^\After, \control^\After_{i{+}1}}
     \eqsepv \forall i \in \ic{0,\horizon{-}1}
     \eqfinv
     \label{eq:couti_dhd_proof}
  \\
  \InstantneousCost_{\doubletime_{\horizon}}
  &: \STATE_{\doubletime_{\horizon}}
    \ni \state^\Before_{\horizon} 
    \mapsto
    \coutfin \np{\state_i^\Before}
    \eqfinp
    \label{eq:coutf_dhd_proof}
\end{align}
\end{subequations}

Now, using Theorem~\ref{cor:instantaneous-costs-in-addition-to-final-cost} we obtain
that the value function  $\Value_{0}$ given in Equation~\eqref{gdp:eq:value_function_dhd}
coincides with the value function $\tilde\Value_{{0}}$ where the sequence of reduced value functions 
  $\{\tilde\Value_{\doubletime_{i}}\}_{i\in \ic{0,\horizon}}$, with \( \tilde\Value_{\doubletime_{i}} :
  \STATE_{\doubletime_{i}} \to \bRR \) for $i \in \ic{0,N}$, is defined by
  \begin{equation}
    \tilde\Value_{\doubletime_{\horizon}} =\FinalCost_{\doubletime_{\horizon}}
    \mtext{ and }
    \tilde\Value_{\doubletime_{i}} = \overlineBelOp{\doubletime_{i{+}1}}{\doubletime_{i}} \tilde\Value_{\doubletime_{i{+}1}}
    \eqsepv \forall i\in\ic{0,\horizon{-}1}
    \eqfinv
    \label{gdp:eq:reduced_value_functions_time_block_additive_bis}
  \end{equation}
  and where the reduced Bellman operator~$\overlineBelOp{\doubletime_{i{+}1}}{\doubletime_{i}}$
  across~$\icinterval{\doubletime_i}{\doubletime_{i{+}1}}$ are given,
  for any $i\in\ic{0,N{-}1}$, for any~$\tilde\varphi_{\doubletime_{i{+}1}} \in
  \espace{L}^{0}_{+}(\STATE_{\doubletime_{i{+}1}})$
  and for any $\optstate_{\doubletime_{i}}\in\STATE_{\doubletime_{i}}$, by
  Equation~\eqref{gdp:eq:overlineBelOp-expression} that we reproduce here:
  \begin{align}
    \big( \overlineBelOp{\doubletime_{i{+}1}}{\doubletime_{i}}\tilde\varphi_{\doubletime_{i{+}1}}\big) (\optstate_{\doubletime_{i}})
    =
    \inf_{\control_{\doubletime_{i}}\in\CONTROL_{\doubletime_{i}}}
    & \int_{\UNCERTAIN_{\doubletime_{i}{+}1}}
      \StochKernel{\doubletime_{i}}{\doubletime_{i}{+}1}
      \kernelargs{\optstate_{\doubletime_{i}}}{\dd\uncertain_{\doubletime_{i}{+}1}}
      \nonumber \\
    &
      \inf_{\control_{\doubletime_{i}{+}1}\in\CONTROL_{\doubletime_{i}{+}1}} \int_{\UNCERTAIN_{\doubletime_{i}{+}2}}
      \StochKernel{\doubletime_{i}{+}1}{\doubletime_{i}{+}2}
      \kernelargs{\optstate_{\doubletime_{i}},\control_{\doubletime_{i}},\uncertain_{\doubletime_{i}{+}1}}{\dd\uncertain_{\doubletime_{i}{+}2}}
      \quad \cdots
      \nonumber    \\
    & \quad
      \inf_{\control_{\doubletime_{i{+}1}-1}\in\CONTROL_{\doubletime_{i{+}1}-1}} \int_{\UNCERTAIN_{\doubletime_{i{+}1}}}
      \StochKernel{\doubletime_{i{+}1}-1}{\doubletime_{i{+}1}}
      \nonumber    \\
    &\quad\quad\quad
      \kernelargs{\optstate_{\doubletime_{i}},
      \control_{\doubletime_{i}},\uncertain_{\doubletime_{i}{+}1},\ldots,\control_{\doubletime_{i{+}1}-2},\uncertain_{\doubletime_{i{+}1}-1}}
      {\dd\uncertain_{\doubletime_{i{+}1}}}
      \nonumber \\
    & \quad\quad
      \Big(
      \InstantneousCost_{\doubletime_i}
      \np{\optstate_{\doubletime_{i}}, \control_{\doubletime_{i}},\uncertain_{\doubletime_{i}{+}1},\ldots,\control_{\doubletime_{i{+}1}-1},\uncertain_{\doubletime_{i{+}1}}}
      \nonumber  \\
    & \quad\quad\quad
      +
      \tilde\varphi_{\doubletime_{i{+}1}}\bp{\Dynamics{\doubletime_{i}}{\doubletime_{i{+}1}}
      \np{\optstate_{\doubletime_{i}}, \control_{\doubletime_{i}},\uncertain_{\doubletime_{i}{+}1},\ldots,\control_{\doubletime_{i{+}1}-1},\uncertain_{\doubletime_{i{+}1}}}}
      \Big)
      \label{gdp:eq:overlineBelOp-expression_DHD}
      \eqfinp \!\!\!\!\!
  \end{align}
  
  It remains to show that Equation~\eqref{gdp:eq:overlineBelOp-expression_DHD}
  gives Equation~\eqref{gdp:eq:BelOpnatural-expression} that is, for $i\in \ic{0,\horizon{-}1}$, 
  $\naturalBelOp{i+1}{i}= \overlineBelOp{\doubletime_{i{+}1}}{\doubletime_{i}}$
  and for $i \in \ic{0,\horizon}$ $\Value_{i}^\Before= \tilde\Value_{\doubletime_i}$.
  Note that, using the definition of $\tau_{i+1}$, we have that $\tau_{i+1}-1=2(i+1)-1= 2i+1=\tau_i+1$ and thus
  we only have two minimization problems in Equation~\eqref{gdp:eq:overlineBelOp-expression_DHD}.
  We obtain successively

  \begin{align}
    \big( \overlineBelOp{\doubletime_{i{+}1}}{\doubletime_{i}}\tilde\varphi_{\doubletime_{i{+}1}}\big) (\optstate_{\doubletime_{i}})
    =
    \inf_{\control_{\doubletime_{i}}\in\CONTROL_{\doubletime_{i}}}
    & \int_{\UNCERTAIN_{\doubletime_{i}{+}1}}
      \StochKernel{\doubletime_{i}}{\doubletime_{i}{+}1}
      \kernelargs{\optstate_{\doubletime_{i}}}{\dd\uncertain_{\doubletime_{i}{+}1}}
      \nonumber \\
    & \quad
      \inf_{\control_{\doubletime_{i}+1}\in\CONTROL_{\doubletime_{i}+1}} \int_{\UNCERTAIN_{\doubletime_{i{+}1}}}
      \StochKernel{\doubletime_{i{+}1}-1}{\doubletime_{i{+}1}}
      \kernelargs{\optstate_{\doubletime_{i}},
      \control_{\doubletime_{i}},\uncertain_{\doubletime_{i}{+}1},\ldots,\control_{\doubletime_{i{+}1}-2},\uncertain_{\doubletime_{i{+}1}-1}}
      {\dd\uncertain_{\doubletime_{i{+}1}}}
      \nonumber \\
    & \quad\quad\quad
      \Big(
      \InstantneousCost_{\doubletime_i}
      \np{\optstate_{\doubletime_{i}}, \control_{\doubletime_{i}},\uncertain_{\doubletime_{i}{+}1},\control_{\doubletime_{i}+1},\uncertain_{\doubletime_{i{+}1}}}
      \nonumber  \\
    & \quad\quad\quad\quad
      +
      \tilde\varphi_{\doubletime_{i{+}1}}\bp{\Dynamics{\doubletime_{i}}{\doubletime_{i{+}1}}
      \np{\optstate_{\doubletime_{i}}, \control_{\doubletime_{i}},\uncertain_{\doubletime_{i}{+}1},\control_{\doubletime_{i}+1},\uncertain_{\doubletime_{i{+}1}}}}
      \Big)
      \tag{by~\eqref{gdp:eq:overlineBelOp-expression_DHD}}   \nonumber
    \\
    =\inf_{\control_{\doubletime_{i}}\in\CONTROL_{\doubletime_{i}}}
    & \int_{\UNCERTAIN_{i}^\After}
      \StochKernel{i-1}{i}^\After \np{\dd\uncertain_{i}^\After}
      \tag{by~\eqref{eq:def_rho_dhd} second case} \\
    & \quad
      \inf_{\control_{\doubletime_{i}+1}\in\CONTROL_{\doubletime_{i}+1}}
      \int_{\UNCERTAIN_{i{+}1}^\Before}
      \StochKernel{i}{{i{+}1}}^\Before \np{\dd\uncertain_{i{+}1}^\Before}
      \tag{by~\eqref{eq:def_rho_dhd} first case} \\
    & \quad\quad\quad
      \Big(
      \coutint_i  % \InstantneousCost_{\doubletime_i}
      \np{\optstate_{\doubletime_{i}}, \control_{\doubletime_{i}},\uncertain_{i}^\After,\control_{\doubletime_{i}+1}} %,\uncertain_{i{+}1}^\Before}
      % \nonumber  \\
      % & \quad\quad\quad
      +
      \tilde\varphi_{\doubletime_{i{+}1}}\bp{
      \dynamics_{i} %   \Dynamics{\doubletime_{i}}{\doubletime_{i{+}1}}
      \np{\optstate_{\doubletime_{i}}, \control_{\doubletime_{i}},\uncertain_{i}^\After,\control_{\doubletime_{i}+1}}}% ,\uncertain_{i{+}1}^\Before}}
      \Big)
      \tag{by~\eqref{eq:couts_dhd_proof} and \eqref{eq:dynamics_dhd_proof}}
    \\
    = \inf_{\control_{i}^\Before\in\CONTROL_{i}^\Before}
    & \int_{\UNCERTAIN_{i}^\After}
      \StochKernel{i-1}{i}^\After \np{\dd\uncertain_{i}^\After}
      \tag{by~\eqref{eq:double_time_line_controls}} \nonumber \\
    & \quad
      \inf_{\control_{i+1}^\After\in\CONTROL_{{i}+1}^\After}
      \Big(
      \coutint_i  % \InstantneousCost_{\doubletime_i}
      \np{\optstate_{\doubletime_{i}}, \control_{{i}}^\Before,\uncertain_{i}^\After,\control_{{i}+1}^\After}
      % ,\uncertain_{i{+}1}^\Before}
      +
      \tilde\varphi_{\doubletime_{i{+}1}}\bp{
      \dynamics_{i} %   \Dynamics{\doubletime_{i}}{\doubletime_{i{+}1}}
      \np{\optstate_{\doubletime_{i}}, \control_{{i}}^\Before,\uncertain_{i}^\After,\control_{{i}+1}^\After}}
    % ,\uncertain_{i{+}1}^\Before}}
      \Big)
      \nonumber
      \eqfinv
  \end{align}
  from which the end of the proof follows.
\end{proof}

\section{Two-time-scale dynamic programming}
\label{Two-time-scale_dynamic_programming}

We present a framework for two-time-scale multistage optimization problems
which is more general 
than in Sect.~\ref{Two_time_scale_optimization_problem},
as we do not require dynamics and states at the fast time-scale.
In~\S\ref{Data_for_two-time-scale_multistage_optimization_problem},
we detail the data needed to formulate the two-time-scale optimization problems
that we consider.
In~\S\ref{Two-time-scale_dynamic_programming_bis}, 
we show how to decompose such problems by slow time blocks.

\subsection{Data for two-time-scale multistage optimization problem}
\label{Data_for_two-time-scale_multistage_optimization_problem}

\subsubsubsection{Two time scales and unified extended timeline}

We consider the same setting as
in~\S\ref{Structure_of_a_two_time_scale_optimization_problem}
in what regards the two time-scales $\SlowSpan$ and $\FastSpan$
\begin{subequations}
  \begin{align}
    \min\SlowSpan &= \bottom\tslow \prec \cdots  \prec \pre\tslow \prec \tslow
    \prec \nex\tslow \prec \cdots \prec \topp\tslow = \max\SlowSpan \eqfinv
    \\ 
    \min\FastSpan & = \bottom\tfast \prec \cdots  \prec \pre\tfast \prec \tfast
                    \prec \nex\tfast \prec \cdots \prec \topp\tfast = \max\FastSpan
                    \eqfinv
  \end{align}
but with an extra \emph{final} slow time~$\nex{\topp{\tslow}}$
(with $\overline{\SlowSpan}=\SlowSpan \cup \na{\nex{\topp{\tslow}}} $) and 
the unified extended timeline $\overline{\SlowSpan{\times}\FastSpan}
=\np{\SlowSpan{\times}\FastSpan}
\cup (\nex{\topp{\tslow}},\bottom{\tfast}) $
\begin{align}
  (\bottom\tslow,\bottom\tfast)
  &\prec (\bottom\tslow,\nex{\bottom\tfast}) \prec 
    \cdots\prec
    (\pre\tslow,\topp\tfast) \prec
    (\tslow,\bottom\tfast)\prec(\tslow,\nex{\bottom\tfast})\prec
    \cdots
    \nonumber
  \\
  &\cdots \prec (\tslow,\pre{\topp\tfast})\prec(\tslow,\topp\tfast)
    \prec(\nex\tslow,\bottom\tfast)\prec \cdots \prec(\topp\tslow,\topp\tfast)
    \prec (\nex{\topp{\tslow}},\bottom{\tfast})
    \eqfinp
\label{eq:extended-timeline_bis}
\end{align}
Note that, at the difference of \S\ref{Structure_of_a_two_time_scale_optimization_problem},
as the state dynamics is only defined at the slow time scale, it is more convenient
to add an extra final time $\nex{\topp{\tslow}}$ rather than an initial extra time
$(\pre{\bottom\tslow},\topp\tfast)$.

The extended timeline
$\overline{\SlowSpan{\times}\FastSpan}$ can be identified with the
time set $\ic{0,\horizon}$, so that we are in the framework of
\S\ref{The_Bertsekas-Shreve_setting}. 
                    
In conformity with the above unified extended timeline ---
and, as in~\S\ref{State_reduction_on_multiple_consecutive_time_blocks_and_dynamic_programming_equations},
we will consider state reduction on multiple consecutive time blocks, but in the
special case where each block is made of all the fast time steps between two
consecutive slow time steps --- we define
    \begin{align}
      \ClosedIntervalClosed{(\tslow,\bottom{\tfast})}{(\tslow,\topp\tfast)}
      &=
        \na{\tslow} \times \FastSpan
 = \na{ (\tslow,\bottom\tfast), (\tslow,\nex{\bottom\tfast}), \ldots, (\tslow,\topp\tfast)}
       \eqsepv
      \forall \tslow \in \SlowSpan
        \eqfinv
        \label{eq:double_time_def}
      \\
     \ClosedIntervalClosed{(\tslow,\bottom{\tfast})}{(\nex\tslow,\bottom\tfast)}
      &=
  \na{ (\tslow,\bottom\tfast), (\tslow,\nex{\bottom\tfast}), \ldots,
        (\tslow,\topp\tfast), (\nex\tslow,\bottom\tfast) }
       \eqsepv
      \forall \tslow \in \SlowSpan
       \eqfinp
    \end{align}
  \end{subequations}

\subsubsubsection{Decision and uncertainty spaces, stochastic kernels}
  
For the rest, to the difference with~\S\ref{Structure_of_a_two_time_scale_optimization_problem}, 
we consider decision Borel spaces
$ \nseqa{\CONTROL_{(\tslow,\tfast)}}{(\tslow,\tfast)\in \SlowSpan\times\FastSpan} $,
 uncertainty Borel spaces
 $\nseqa{\UNCERTAIN_{(\tslow,\tfast)}}{(\tslow, \tfast)\in\overline{\SlowSpan\times\FastSpan}}$,
 stochastic kernels
\begin{subequations}
  \label{eq:natural_kernel_bis}
  \begin{align}
    \rho_{(\pre\tslow,\topp\tfast):(\tslow,\bottom\tfast)}
    &\in
      \Delta\np{\UNCERTAIN_{(\tslow,\bottom\tfast)}}
      \eqsepv
      \forall \tslow \in \overline{\SlowSpan}% \SlowSpan\setminus\na{\bottom\tslow}
      \eqfinv
      \label{eq:natural_kernel_slow_bis}
    \\
    \rho_{(\tslow,\pre\tfast):(\tslow,\tfast)}
    &:
   \prod_{\tfast' = \bottom\tfast}^{\pre\tfast} \UNCERTAIN_{(\tslow,\tfast')}
   \to       \Delta\np{ \UNCERTAIN_{(\tslow,\tfast)}}
      \eqsepv
      \forall \tslow \in \SlowSpan
      \eqsepv
      \forall \tfast \in \FastSpan\setminus\na{\bottom\tfast}
      \eqfinp
      \label{eq:natural_kernel_fast_bis}
  \end{align}
  Note that, for a given $\tslow \in \SlowSpan$,
  each stochastic kernel $\rho_{(\tslow,\tfast)}$
  depends at most on the noises of the slow time block
  \( \ClosedIntervalClosed{(\tslow,\bottom{\tfast})}{(\tslow,\topp\tfast)} \)
  in~\eqref{eq:double_time_def}.
  The (constant) assumption~\eqref{eq:natural_kernel_slow_bis} and the
  (single block) assumption~\eqref{eq:natural_kernel_fast_bis} correspond
  to stochastic independence between time blocks, and will be useful in
  the proof of Proposition~\ref{gdp:prop:two_time_scale_decomposition_bis}.
\end{subequations}

\subsubsubsection{History spaces}

In conformity with the unified extended
timeline~\eqref{eq:extended-timeline_bis}
and with the decision and uncertainty spaces, we deduce the history Borel spaces and the histories
for all $(\tslow,\tfast) \in \overline{\SlowSpan{\times}\FastSpan}$
\begin{subequations}
  \begin{align}
  \HISTORY_{(\tslow,\tfast)}
    & =
      \UNCERTAIN_{(\bottom\tslow,\bottom\tfast)} \times
      \prod_{(\bottom\tslow,\nex{\bottom\tfast}) \preceq (\tslow',\tfast') \preceq (\tslow,\tfast)}
      \Bp{\CONTROL_{\pre{(\tslow',\tfast')}}\times\UNCERTAIN_{(\tslow',\tfast')}}
      \eqfinv \\
  \hist[(\tslow,\tfast)]
    & =
      \Bp{\alea[{(\bottom\tslow,\bottom\tfast)}] ,
      \bseqp{
      { \control_{\pre{(\tslow',\tfast')}} ,  \alea[{(\tslow',\tfast')}]}}
      {(\bottom\tslow,\nex{\bottom\tfast}) \preceq (\tslow',\tfast') \preceq (\tslow,\tfast)}}
      \eqfinv
      \label{eq:tts-initial-condition_bis}
  \end{align}
\end{subequations}
and, for suitable indices, the partial history sets and the partial histories
\begin{subequations}
  \begin{align}
    \HISTORY_{(\tslow,\tfast):(\tslow',\tfast')}
    &=
      \prod_{(\tslow,\tfast) \preceq (\tslow'',\tfast'') \preceq (\tslow',\tfast')}
      \np{ \CONTROL_{\pre{(\tslow'',\tfast'')}} \times \UNCERTAIN_{(\tslow'',\tfast'')} }
      % \eqsepv \forall (\tslow,\tfast) \in \SlowSpan\times\FastSpan
      \eqfinv
    \\
    \hist[(\tslow,\tfast):(\tslow',\tfast')]
    & =
      \bp{
      \nseqp{
      { \control_{\pre{(\tslow'',\tfast'')}} ,  \alea[{(\tslow'',\tfast'')}]}}
      {(\tslow,\tfast) \preceq (\tslow'',\tfast'') \preceq (\tslow',\tfast')}}
    % \eqsepv \forall (\tslow,\tfast) \in \SlowSpan\times\FastSpan
      \eqfinp
  \end{align}
\end{subequations}

% As in~\S\ref{Structure_of_a_two_time_scale_optimization_problem},
% we consider a slow time scale represented by a finite totally ordered set~$\np{\SlowSpan,\preceq}$,
% a fast time scale represented by a finite totally ordered       set~$\np{\FastSpan,\preceq}$,
% and deduce a unified timeline~\eqref{eq:timeline}
% and an extended timeline~\eqref{eq:extended-timeline}. 

\subsubsubsection{State reductions and slow time scale dynamics}

Let~$\nseqa{\STATE_{\tslow}}{\tslow\in\overline{\SlowSpan}}$
be a sequence of Borel \emph{(state) spaces}
(where \( \STATE_{\bottom\tslow} = \COST[(\bottom\tslow,\bottom\tfast)][] \)),
$\nseqa{\theta_{\tslow}}{\tslow\in\overline{\SlowSpan}}$
be a sequence of Borel-measurable \emph{reduction mappings}
$ \theta_{\tslow} : \HISTORY_{(\tslow,\bottom\tfast)} \to  \STATE_{\tslow} $
(where \( \theta_{\bottom\tslow} = \textrm{Id} :\HISTORY_{(\bottom\tslow,\bottom\tfast)}
=\COST[(\bottom\tslow,\bottom\tfast)][]
\to\STATE_{\bottom\tslow} = \COST[(\bottom\tslow,\bottom\tfast)][] \)),
  and $\nseqa{ \dynamics_{\tslow} }{\tslow\in\SlowSpan}$
  be a sequence of  Borel-measurable \emph{dynamics}
  \begin{equation}
    \dynamics_{\tslow} : % \Dynamics{(\tslow,\bottom\tfast)}{(\tslow,\topp\tfast)} :
    \STATE_{\tslow} \times \HISTORY_{(\tslow,\nex{\bottom\tfast}):(\nex\tslow,\bottom\tfast)}
    \to \STATE_{\nex\tslow} \eqsepv \forall \tslow\in \SlowSpan
    \eqfinp
  \end{equation}

\subsubsubsection{Cost functions}

\begin{subequations}
  We suppose given a family
  $\nseqa{\instantcost[\tslow]}{\tslow\in\SlowSpan}$
  of \emph{slow time scale nonnegative lower semianalytic cost functions},
  with
  \begin{align}
    \instantcost[\tslow] :
      \STATE_{\tslow} \times
    \HISTORY_{(\tslow,\nex{\bottom\tfast}):(\nex\tslow,\bottom\tfast)}
      \to \bRR \eqsepv \forall \tslow\in \SlowSpan
      \intertext{and a \emph{slow time scale nonnegative lower semianalytic final cost function}~$\instantcost[\nex{\topp\tslow}]$}
      \instantcost[\topp\tslow] :
      \STATE_{\nex{\topp\tslow}}       \to \bRR \eqfinp 
  \end{align}
With the slow time scale cost functions and final cost function, we
make up, by summation, a (nonnegative lower semianalytic)
cost criterion $\criterion:\HIST[{(\nex{\topp\tslow},\bottom\tfast)}] \to \bRR$
given by 
\begin{equation}
  \criterion\np{\history_{(\nex{\topp\tslow},\bottom\tfast)}}
  =  \dsum[\tslow\in\SlowSpan] \instantcost[\tslow]
  \bp{
    \theta_{\tslow}\np{\history_{(\tslow,\bottom\tfast)}},
    \history_{(\tslow,\nex{\bottom\tfast}):(\nex\tslow,\bottom\tfast)}
  }
  + \instantcost[\nex{\topp\tslow}]\bp{\theta_{\nex{\topp\tslow}}\np{\history_{\nex{(\topp\tslow},\bottom\tfast}}}
  \eqfinp
  \label{eq:intertemporal-criterion_bis}
\end{equation}
\end{subequations}

\subsection{Two-time-scale dynamic programming}
\label{Two-time-scale_dynamic_programming_bis}

\begin{proposition}
  Suppose that   the family
  $\big( \nseqa{ \STATE_{\tslow} }{\tslow\in\SlowSpan}$,
  $\nseqa{\theta_{\tslow}}{\tslow\in\SlowSpan}$,
$\nseqa{ \dynamics_{\tslow} }{\tslow\in\SlowSpan}$ 
is a state reduction across the consecutive time blocks~$\ClosedIntervalClosed{(\tslow,\bottom\tfast)}{(\tslow,\topp\tfast)}$,
for $\tslow\in\SlowSpan$
(where we identify \( \Dynamics{(\tslow,\bottom\tfast)}{(\tslow,\topp\tfast)}= 
\dynamics_{\tslow} \) in Definition~\ref{gdp:de:reduction-dynamics_family}). 

  Then, the  multistage (two-time-scale) stochastic optimization problem ---
  formulated like in~\eqref{eq:Vzero} with the data
  in~\S\ref{Data_for_two-time-scale_multistage_optimization_problem} --- 
  has a solution given by a dynamic programming equation at
  the slow scale. More precisely, let $\nseqa{\tilValue[\tslow]}{\tslow \in
    \overline{\SlowSpan}}$
  be a sequence of functions given by
  $\tilValue[\nex{\topp\tslow}] = \instantcost[\nex{\topp\tslow}]$
  and, for $\tslow \in \SlowSpan$, 
  by the backward induction% \footnote{Here again, the formula~\eqref{gdp:eq:overlineBelOp-expression-rephrased}
  % represents a nested sequence of infima of integrals
  % (with respect to different stochastic kernels).}
  \begin{align}
    & \tilValue[\tslow]\np{\state_{\tslow}}
      =
      \inf_{\control_{(\nex\tslow,\bottom\tfast)}\in\CONTROL_{(\nex\tslow,\bottom\tfast)}}
      \int_{\UNCERTAIN_{(\nex\tslow,\nex{\bottom\tfast})}}
      \rho_{(\nex\tslow,{\bottom\tfast}):(\nex\tslow,\nex{\bottom\tfast})}
      \kernelargs{\uncertain_{\nex\tslow}}{\dd\uncertain_{(\nex\tslow,\nex{\bottom\tfast})} }
      \; \cdots
      \nonumber    \\
    & \quad\quad
      \inf_{\control_{(\nex\tslow,\pre{\topp\tfast})}\in\CONTROL_{(\nex\tslow,\pre{\topp\tfast})}}
      \int_{\UNCERTAIN_{(\nex\tslow,\topp\tfast)}}
          \rho_{(\nex\tslow,\pre{\topp\tfast}):(\nex\tslow,{\topp\tfast})}
          \kernelargs{\uncertain_{\nex\tslow},\uncertain_{(\nex\tslow,\nex{\bottom\tfast})},
                      \cdots,\uncertain_{(\nex\tslow,\pre{\topp\tfast})}}
          {\dd\uncertain_{(\nex\tslow,{\topp\tfast})}}
          \nonumber \\
     &
      \hphantom{xxxxxxxxx}
      \Big(
      \instantcost[\tslow]
      \bp{\state_{\tslow},
      \nseqp{
      { \control_{\pre{(\tslow',\tfast')}} ,  \alea[{(\tslow',\tfast')}]}}
      {(\tslow,\nex{\bottom\tfast}) \preceq (\tslow',\tfast') \preceq (\nex\tslow,\bottom\tfast)}}
    \nonumber \\
    & \hphantom{xxxxxxxxxxxxxxxx}
      +
      \tilValue[\nex\tslow]\Bp{
      \dynamics_{\tslow}
      % \Dynamics{(\tslow,\bottom\tfast)}{(\tslow,\topp\tfast)} 
      \bp{\state_{\tslow},
      \nseqp{
      { \control_{\pre{(\tslow',\tfast')}} ,  \alea[{(\tslow',\tfast')}]}}
      {(\tslow,\nex{\bottom\tfast}) \preceq (\tslow',\tfast') \preceq (\nex\tslow,\bottom\tfast)}}}
      \Big)
      \label{gdp:eq:overlineBelOp-expression-rephrased_bis}
      \eqfinp 
  \end{align}
  Then, the value of the optimization problem~\eqref{eq:Vzero}
  is given by
  $\tilValue[\bottom\tslow]\np{\state_{\bottom\tslow}}$,
  where the initial condition $\state_{\bottom\tslow}$
  corresponds to~$\uncertain_{0}$ in~\eqref{eq:Vzero},
  as stated by~\eqref{eq:tts-initial-condition_bis}.
  \label{gdp:prop:two_time_scale_decomposition_bis}
\end{proposition}

\begin{proof}
  The proof is an application of
  Theorem~\ref{cor:instantaneous-costs-in-addition-to-final-cost}
  with the help of
  Remarks~\ref{rem:composed-state-dynamics-as-reduction-mapping}
  and~\ref{rem:block-independent-noises-and-kernel}.
  First, we have framed the multistage (two-time-scale) stochastic
  optimization problem
  in the formalism of~\S\ref{The_Bertsekas-Shreve_setting}
  with the help of the extended timeline~\eqref{eq:extended-timeline_bis}.
  Second, we have by assumption a state reduction
  at times $\nseqa{(\tslow,\bottom\tfast)}{\tslow \in \SlowSpan}$
  by composition of the state dynamics. Moreover,
  as the slow time scale stochastic kernels given by Equation~\eqref{eq:natural_kernel_slow_bis} are constant,
  the state reduction across the slow time scale is compatible
  with the stochastic kernels
  (see Remark~\ref{rem:block-independent-noises-and-kernel}).
  We are thus able to apply
  Theorem~\ref{cor:instantaneous-costs-in-addition-to-final-cost}
  and obtain the slow time scale Bellman
  recursion~\eqref{gdp:eq:overlineBelOp-expression-rephrased_bis}
  as a special case of Equation~\eqref{gdp:eq:overlineBelOp-expression}.
  This ends the proof.
\end{proof}

\makeatletter\def\@biblabel#1{[AM#1]}\makeatother
\ifx\undefined\noopsort\newcommand{\noopsort}[1]{}\fi\ifx\undefined\allcaps\def\allcaps#1{#1}\fi

\preprintstop
\fi

\end{document}